\RequirePackage{fix-cm}
\documentclass[smallextended]{svjour3}       
\smartqed  
\usepackage{float}
\usepackage{graphicx}
\usepackage{url}
\usepackage{amsfonts}
\usepackage{amssymb}
\usepackage{bm}
\newcommand{\Proof}{\vspace{-.1cm}\noindent {\sc Proof}. }

\newcommand{\inc}{\subseteq}
\newcommand{\set}[1]{\{#1\}}
\newcommand{\setc}[2]{\set{#1\!\mid\! #2}}
\newcommand{\union}{\cup}		
\newcommand{\inter}{\cap}		
\newcommand{\Union}{\bigcup}		
\newcommand{\Inter}{\bigcap}		
\newcommand{\Sub}[3]{#1[#2\leftarrow #3]}	
\newcommand{\comp}{\circ}
\newcommand{\nin}{\not\in}

\newcommand{\restrH}[2]{\hyper{#1}\backslash #2}
\newcommand{\hyper}[1]{{\bm #1}}
\renewcommand{\setminus}{\backslash}

\newcommand{\incs}{\subsetneq}
\newcommand{\eqdef}{\stackrel{\Delta}{=}}
\newcommand{\bcdot}{{\bm \cdot}}

\def\mathraise#1#2#3#4{
\mathchoice
{\raisebox{#2}{$\displaystyle{#1}$}}
{\raisebox{#2}{$#1$}}
{\raisebox{#3}{$\scriptstyle{#1}$}}
{\raisebox{#4}{$\scriptscriptstyle{#1}$}}
}
\def\mathraiseord#1#2#3#4{\mathord{\mathraise#1{#2}{#3}{#4}}}
\def\valup{\mathraiseord{\uparrow}{0.29ex}{0.21ex}{0.16ex}}

\newcommand{\Alt}{ \mid\!\!\mid  }

\usepackage{tikz}
\usetikzlibrary{positioning}
\usetikzlibrary{decorations.pathmorphing}
\tikzset{>=stealth}
\newcommand \seq[2]{\shortstack{$#1$ \\ \mbox{}\\
                    \mbox{}\hrulefill\mbox{}\\ \mbox{}\\ $#2$}}
\newcommand \seqtwo[3]{\shortstack{$#1$ \\ $#2$ \\ \mbox{}\\
                    \mbox{}\hrulefill\mbox{}\\ \mbox{}\\ $#3$}}
\newcommand \seqthree[4]{\shortstack{$#1$ \\ $#2$ \\ $#3$ \\ \mbox{}\\
                    \mbox{}\hrulefill\mbox{}\\ \mbox{}\\ $#4$}}

\usepackage{graphicx}
\definecolor{cof}{RGB}{219,10,71}
\definecolor{pur}{RGB}{186,146,162}
\definecolor{greeo}{RGB}{70,63,169}
\definecolor{greet}{RGB}{52,111,72}
 \usepackage{hyperref}
\begin{document}

\title{Syntactic aspects of hypergraph polytopes 
}


\author{Pierre-Louis Curien        \and Jovana Obradovi\'c 
        \and Jelena Ivanovi\'c 
}

\author{Pierre-Louis Curien        \and Jovana Obradovi\'c 
        \and Jelena Ivanovi\'c 
}

\institute{Pierre-Louis Curien \at IRIF, Univ. Paris Diderot, $\pi r^{2}$, INRIA and CNRS \\
              \email{curien@irif.fr}           
           \and
           Jovana Obradovi\'c \at IRIF, Univ. Paris Diderot, $\pi r^{2}$, and INRIA  \\
           \email{jovana@irif.fr}
           \and
    Jelena Ivanovi\'c \at   University of Belgrade\\     
    \email{jelena.ivanovic@arh.bg.ac.rs}
}

\date{Received: date}

\maketitle

\begin{abstract}
This paper introduces an inductively defined  tree notation for all the faces of polytopes arising from a simplex by truncations.  This notation allows us to view inclusion of faces as the process of contracting  tree edges.  Our notation instantiates to the well-known notations for the faces of  associahedra and  permutohedra.    Various authors have independently introduced combinatorial tools for describing such polytopes.
We build  on the particular approach developed by  Do\v sen and Petri\'c, who 
used the formalism of   hypergraphs to describe the  interval of polytopes  from the simplex to the permutohedron.
This interval was further stretched  by Petri\'c to allow truncations of faces that are themselves obtained by truncations, and iteratively so. Our notation applies to all these polytopes.  We illustrate this  by showing that it instantiates  to a notation for the faces of the permutohedron-based associahedra, that consists of  parenthesised words with holes. 
 Do\v sen and Petri\'c have exhibited some families of hypergraph polytopes  (associahedra, permutohedra, and 
hemiassociahedra) describing
 the coherences, and the coherences between coherences etc., arising by weakening
  sequential and parallel associativity of operadic composition.
We complement their work with a criterion allowing us  to recover  the information whether edges of  these ``operadic polytopes"  come from sequential, or from parallel associativity.
We also give  alternative proofs for some of the original results of Do\v sen and Petri\'c.

\end{abstract}

\section{Introduction}\label{intro-section}
Classically, a (convex) polytope  is defined as a bounded intersection of a finite set of half-spaces. More precisely, a polytope $P$ is specified as the set of solutions to a system  $Ax\geq b$ of linear inequalities, where $A$ is an $m\times n$ matrix, $x$ is an $n\times 1$ column vector of variables, and $b$ is an $m\times 1$ column vector of constants. Here,  $n$ is the dimension of the ambient  space containing $P$, and $m$ is the number of half-spaces defining $P$.
 The actual dimension of $P$ is the maximum dimension   of an open ball contained in $P$.

A  {\em face} of $P$ is any intersection of $P$ with one of its bounding hyperplanes (such a hyperplane  intersects $P$ and bounds a closed half-space containing $P$).  
Following usual terminology, the $0$-dimensional (resp. $1$-dimensional) faces of a polytope are called {\em vertices} (resp. {\em edges}), and if $P$ is $n$-dimensional, we call its $(n-1)$-dimensional faces {\em facets}. 
If the definition of a face  is extended to allow
the empty set to be considered as a face, then 
the faces of a convex polytope form a bounded lattice called its {\em face lattice}, the partial ordering being the set containment of faces. The whole polytope (resp. the empty set) 
 is the  maximum (resp. minimum) element of the lattice.
 
\smallskip
 As opposed to the classically (or geometrically) defined polytopes, an {\em abstract polytope} is a structure that captures only the combinatorial properties of the face lattice of a  polytope,  ignoring some of its other properties, particularly measurable ones, such as angles, edge lengths, etc. An abstract polytope is given as a set of faces, together with an order relation  that satisfies certain axioms reflecting the incidence properties
 of polytopes in the classical sense.
 
\smallskip
 In \cite{DP-HP},  Do\v sen and Petri\' c investigate a family of polytopes that may be obtained by truncating the vertices, edges and other faces of simplices of any finite dimension. The permutohedra are limit cases in that family, where all possible truncations have been made. 
 The limit cases at the other end, where no truncation has been made, are simplices.
(Alternatively, one may choose the permutohedron as starting point, and reach the simplex  by successive contractions, see e.g. \cite{T-RAP}).  Other independent  (and even predating)  approaches have been developed for describing polytopes in this family \cite{FK04,FS05,CD-CCGA,D-RGA,PRW08,P09}.
While  the combinatorial description in all these works is essentially the same (with a different terminology), the ways of describing the geometric realisation are quite diverse, as we shall point out.
\\
\indent 
 An easy example of a transition  from a simplex to a permutohedron 
is obtained by truncating all the vertices of a $2$-dimensional simplex (i.e., a triangle) to get a $2$-dimensional permutohedron (i.e., a hexagon):  
\begin{center}
\begin{tikzpicture}[scale=0.7]
    \node[minimum size=0mm, inner sep=0mm] (E) at (-1,0) {};
    \node[minimum size=0mm, inner sep=0mm] (G) at (1,0) {};
    \node[minimum size=0mm, inner sep=0mm] (F) at (0,1.7) {};
    \node[minimum size=0mm, inner sep=0mm] (J) at (-0.5,1.2) {};
    \node[minimum size=0mm, inner sep=0mm] (K) at (0.5,1.2) {};
    \node[minimum size=0mm, inner sep=0mm] (L) at (-0.75,0.6) {};
     \node[minimum size=0mm, inner sep=0mm] (M) at (0.75,0.6) {};
     \node[minimum size=0mm, inner sep=0mm] (N) at (-0.35,-0.15) {};
     \node[minimum size=0mm, inner sep=0mm] (O) at (0.35,-0.15) {};
    \draw[dashed] (J)--(K);
    \draw[dashed] (L)--(N);
    \draw[dashed] (M)--(O);
    \draw[-] (E)--(G) node [midway,above] {};
    \draw[-] (G)--(F) node [midway,right] {};
    \draw[-] (F)--(E) node [midway,left] {};
\end{tikzpicture}
\end{center}
In higher dimensions, the number of possible truncations increases with the number of faces of different dimensions. For example, at dimension $3$ we can truncate not only the vertices of a tetrahedron, but also its edges.  
The connected  subsets of a hypergraph $\hyper{H}$ with $n$ vertices  act as truncating instructions to be  applied to the    simplex of dimension $n-1$.  The polytopes  obtained in this manner are called {\em hypergraph polytopes}.  

\smallskip
In \cite{DP-HP}, the  faces of  hypergraph polytopes are  named by combinatorial objects called {\em constructs}, for which we develop here a new approach.  While they were originally defined in  \cite{DP-HP} as certain sets of connected subsets of a hypergraph, we define them   as decorated trees obtained in an algorithmic manner; this dynamic point of view extends to the definition of the partial order on constructs
 (Section \ref{inductive-section}).  
We show the equivalence with the original definition of Do\v sen and Petri\' c in Section \ref{geometric-realisation-section}, where we also provide an  alternative proof for the main theorem of \cite{DP-HP}, stating the order-isomorphism between the poset of constructs and the poset of faces in the geometric realisation. Unlike the original proof, our proof builds the isomorphism explicitly.

\smallskip
Do\v sen and Petri\'c developed hypergraph polytopes    in connection to their work on the categorification of operads \cite{DP-WCO}.  The coherences arising in this setting display themselves as faces of some hypergraph polytopes.
We complement their work  with a criterion for recognising whether edges in these polytopes arise from sequential or parallel associativity isomorphisms (Section \ref{operadic-application-section}).

\smallskip
Finally, in Section \ref{stretching-section}, we show how to extend our tree notation for constructs to cover iterated truncations, i.e., truncations of faces themselves obtained after (possibly iterated) truncations,   as captured combinatorially  in \cite{P-SISP}, and we illustrate it for the case of the permutohedron-based associahedron (underlying the coherences of symmetric monoidal categories). We present an ad hoc notation for the faces of this polytope (in any finite dimension), based on words with holes and directly suggested by our construct notation.
%

\smallskip
We have tried to give, as much as possible, a self-contained exposition of the material presented.

\smallskip\noindent

{\em Terminological warning:}
Throughout the paper, there will be trees (all rooted), graphs, hypergraphs, and polytopes, sometimes discussed next to  each other. When speaking about ``vertices'' or ``edges'', it should always be clear to which of these structures we are referring. 

We shall  use two notions of subtree. By a {\em subtree of a construct} $T$ we shall  mean a  tree obtained by picking a node of $T$ and taking all its descendants. But in the context of {\em operadic trees} ${\cal T}$ (Section \ref{operadic-application-section}), we shall call subtree {\em any connected subset of} ${\cal T}$.

\vspace{-.2cm}
\section{Hypergraph polytopes and constructs} \label{inductive-section}

\vspace{-.2cm}
In this section, we recall the definition of a hypergraph and some basic related notions. Then we give our own definition of constructs and of the partial ordering between them, postponing to Section \ref{geometric-realisation-section} the proof that these coincide up to isomorphism with the definitions given in \cite{DP-HP}.

\vspace{-.3cm}
\subsection{Hypergraphs}  \label{hypergraph-definition}
\vspace{-.1cm}
A hypergraph   is given by a set  $H$ of vertices (the carrier), and a subset $\hyper{H}\inc {\cal P}(H)\setminus\emptyset$ such that $\Union \hyper{H}=H$. The elements of $\hyper{H}$ are called the {\em hyperedges} of $\hyper{H}$.  We always assume that $\hyper{H}$ is {\em atomic}, by which we mean that $\set{x}\in \hyper{H}$, for all $x\in H$. 
 Identifying $x$ with $\set{x}$, $H$ can be seen as the set of  hyperedges of cardinality $1$, also called {\em vertices}. We shall always use the convention to give the same name to the hypergraph and to its carrier, the former being the bold version of the latter.
A hyperedge of cardinality 2 is called an {\em edge}.  Note that any ordinary graph $(V,E)$ can be viewed as the atomic hypergraph
$\setc{\set{v}}{v\in V} \union \setc{e}{e\in E}$ (with no hyperedge of cardinality $\geq 3$).

\smallskip
If $\hyper{H}$ is a hypergraph, and  if $X\inc H$, we set
\vspace{-.1cm}
$$\hyper{H}_X=\setc{Z}{Z\in \hyper{H}\;\mbox{and}\; Z\inc X}.$$
\vspace{-.5cm}

\noindent
We say that $\hyper{H}$ is {\em connected} if there is no non-trivial partition $H=X_1\union X_2$ such that $\hyper{H}=\hyper{H}_{X_1}\union \hyper{H}_{X_2}$. All our hypergraphs will be finite.  It is easily seen that for each finite hypergraph there exists a partition
$H=X_1\union\ldots\union X_m$ such that each $\hyper{H}_{X_i}$ is connected and $\hyper{H}=\Union(\hyper{H}_{X_i})$.  The $\hyper{H}_{X_i}$'s are called the {\em connected components} of $\hyper{H}$.
We shall also  use  the following notation:
\vspace{-.1cm}
$$\restrH{H}{X}=\hyper{H}_{H\setminus X}.$$
\vspace{-.6cm}

\noindent
As a (standard) abuse of notation, we call
a  non-empty subset $X$ of vertices connected (resp. a connected component) whenever $\hyper{H}_X$ is connected (resp. a connected component).  
We define the  {\em saturation} of  $\hyper{H}$  as the hypergraph
\vspace{-.2cm}
$${\it Sat}(\hyper{H})=\setc{X}{\emptyset\incs X\inc H\;\mbox{and}\;\hyper{H}_X\;\mbox{is connected}}.$$ 
\vspace{-.6cm}

\noindent
A hypergraph is called saturated  when $\hyper{H}={\it Sat}(\hyper{H})$.  Atomic and saturated hypergraphs are called {\em building sets} in the works of Postnikov et al. \cite{PRW08,P09}, and are generalised, with the same name, from the present setting of ${\cal P}(H)$ to that of arbitrary finite lattices in the works of  Feichtner et al. \cite{FK04,FS05}.

The notation
\vspace{-.2cm}
$$\hyper{H},X  \leadsto H_1,\ldots,H_n \quad (\mbox{resp.}\;\; \hyper{H},X  \leadsto \setc{H_i}{i\in I} )$$
\vspace{-.6cm}

\noindent
 will mean that  $H_1,\ldots,H_n\inc H\setminus X$ are  the (resp. $ \setc{H_i}{i\in I}$ is the set of)  connected components of $\restrH{H}{X}$.  
We shall  write
$\hyper{H}_i$  for
$\hyper{H}_{H_i}$. 

\smallskip
We  call a {\em quasi-partition of a set} $X$  a collection of disjoint (possibly empty) subsets whose union is $X$. We shall need the  following   (standard) property.
\begin{lemma}\label{partition-of-i}
Let $\hyper{H}$ be a connected hypergraph, and let $Y\subseteq X\subseteq H$. Let 
$\hyper{H},Y\leadsto \{K_j\,|\,j\in J\}$ and $\hyper{H},X\leadsto \{H_i\,|\,i\in I\}$. 
Then the following two claims hold:
\begin{enumerate}
\item If
$ H_i\cap  K_j\neq\emptyset$, then $H_i\subseteq K_j$.
\item There exists a quasi-partition $\{I_j\,|\, j\in J\}$ of $I$, such that, for each $j\in J$, $K_j \backslash X=\bigcup_{i\in I_j} H_i$. 
Consequently, we have $\hyper{K}_j ,X\leadsto \{H_i\,|\, i\in I_j\}$.
\end{enumerate}
\end{lemma}

\vspace{-.5cm}
\subsection{Constructs and constructions} \label{construct-section}

 A {connected} hypergraph  $\hyper{H}$ gives rise to a partial order of {\em constructs}, which we define below  inductively.
\newpage
\begin{definition} \label{inductive-construct}
Let $\hyper{H}$ be a connected hypergraph and  $Y$ be    an arbitrary non-empty  subset of  $H$:
\begin{itemize}
\item  If $Y = H$, then  the one-node tree decorated with $H$,  written $H$, is a construct of $\hyper{H}$.
\item Otherwise, if   $\hyper{H},Y  \leadsto H_1,\ldots,H_n$, and if  $T_1,\ldots,T_n$ are constructs of $\hyper{H}_1,\ldots,\hyper{H}_n$, respectively, then the
tree whose root is decorated by $Y$, with $n$ outgoing edges on which the respective $T_i\,$'s are grafted, written
$Y(T_1,\ldots,T_n)$, is a construct.  
\end{itemize}
A {\em construction} is a construct  whose nodes are all decorated with singletons.  We shall often use the letter $V$ to denote a construction (since constructions denote  vertices in the geometric realisation, see Section \ref{iso-section}). 
\end{definition}
In $Y(T_1,\ldots,T_n)$, the order  of the constructs $T_1,...,T_n$ is irrelevant.  
We shall  write
$Y\setc{T_i}{i\in I}$  when  the constructs $T_i$ are indexed over some finite set $I$.  When  $I=\emptyset$, we get that    $Y\setc{T_i}{i\in I}$ stands for $Y$,  corresponding to  the base case in Definition \ref{inductive-construct} (note that the only hypergraph with an empty set of connected components is the empty one). It is also convenient to allow ourselves to write $\emptyset\setc{T_i}{i\in I}$, with $I=\set{i_0}$, as a stuttering form of 
$T_{i_0}$ (if $Y=\emptyset$, we are left with building a construct of the original hypergraph).

\medskip
The intuition behind this definition is algorithmic: {\em a construct is built  by picking  a non-empty subset $Y$ of  $H$ and then branching to the connected components of $\restrH{H}{Y}$, and continuing recursively in all the branches}.


\smallskip
The labels of the nodes of a construct of $\hyper{H}$ form a partition of $H$.
We shall  freely confuse the nodes with their labels, since they are a  fortiori all distinct.  For every node $Y$ of $T$, we denote by $\valup_T(Y)$  (or simply $\valup(Y)$) the union of the  labels of the descendants of $Y$ in $T$ (all the way to the leaves), including $Y$. 
For every construct $T$ of $\hyper{H}$ and every node $Z$ of $T$, the subtree of $S$ rooted at $Z$ is a construct of  $\hyper{H}_{\valup_T(Z)}$. 

The notation 
$T:\hyper{H}$  will mean that  $T$ is a construct of $\hyper{H}$.  The following formal system summarises our definition of constructs:
\begin{center}
\resizebox{8.5cm}{!}{\fbox{$
\seq{}{H:\hyper{H}}\quad\quad
\seq{\hyper{H},X\leadsto H_1,\ldots,H_n\quad\quad T_1:\hyper{H}_1,\ldots,T_n:\hyper{H}_n}
{X(T_1,\ldots,T_n):\hyper{H}}
$}}
\end{center}

We note that  while  the inductively-defined  constructions in tree form appear   in \cite{DP-HP}[Section 3] and  in  \cite{PRW08}[Proposition 8.5] exactly like in Definition \ref{inductive-construct}, these authors  did not notice or  exploit the fact that the tree notation could be extended to all constructs simply by replacing singletons with arbitrary subsets. As we shall see, this simple observation gives additional insights. In particular,  it allows us to formulate various equivalent and useful  characterisations of the partial order between constructs. 

The  tree notation for all constructs appears  in \cite{FS05}[Proposition 3.17], but without an inductive characterisation.

\vspace{-.3cm}
\subsection{Ordering constructs} \label{construct-ordering-section}
\vspace{-.2cm}
We next define a partial order between constructs. The algorithmic intuition  is that, given $S$, one can get a larger construct by contracting an edge of $S$, and then merging the decorations of the two nodes related by that edge, as illustrated in the following picture:
\begin{center}
\resizebox{7.75cm}{!}{\begin{tikzpicture}[scale=0.7]
    \node (Y) [rectangle,draw=none,minimum size=0.75cm,inner sep=0.0mm] at (0,-0.9) {\Large $Y$};
    \node (X) [rectangle,draw=none,minimum size=0.75cm,inner sep=0.1mm]  at (-1.5,1) {\Large $X$};
    \node (T11) [circle,draw=none,minimum size=0.75cm,inner sep=0.1mm] at (-2.5,2.5) {\large  $T_{11}$};
   \node (d1) [circle,draw=none,minimum size=0.75cm,inner sep=0.1mm] at (-1.5,2.5) { $\cdots$};
   \node (T1m) [circle,draw=none,minimum size=0.75cm,inner sep=0.1mm] at (-0.5,2.5) {\large $T_{1m}$};
   \node (T2) [circle,draw=none,minimum size=0.75cm,inner sep=0.1mm] at (1,2.5) {\large $T_{2}$};
   \node (d2) [circle,draw=none,minimum size=0.75cm,inner sep=0.1mm] at (2,2.5) {  $\cdots$};
   \node (Tn) [circle,draw=none,minimum size=0.75cm,inner sep=0.1mm] at (3,2.5) {\large  $T_{n}$};
 
  \node (E) [circle,draw=none,minimum size=0.75cm,inner sep=0.1mm] at (4.75,0.5) {\Large $\leq^{\bf H}$};

   \node (Y1) [rectangle,draw=none,minimum size=0.75cm,inner sep=1.5mm] at (9,-0.9) {\Large  $Y\cup X$};
    \node (T111) [circle,draw=none,minimum size=0.75cm,inner sep=0.1mm] at (6,2.5) {\large  $T_{11}$};
   \node (d11) [circle,draw=none,minimum size=0.75cm,inner sep=0.1mm] at (7.25,2.5) { $\cdots$};
   \node (T1m1) [circle,draw=none,minimum size=0.75cm,inner sep=0.1mm] at (8.5,2.5) {\large  $T_{1m}$};
   \node (T21) [circle,draw=none,minimum size=0.75cm,inner sep=0.1mm] at (10,2.5) {\large  $T_{2}$};
   \node (d21) [circle,draw=none,minimum size=0.75cm,inner sep=0.1mm] at (11,2.5) {\large  $\cdots$};
   \node (Tn1) [circle,draw=none,minimum size=0.75cm,inner sep=0.1mm] at (12,2.5) { \large $T_{n}$};
\draw (T111)--(Y1);
\draw (T1m1)--(Y1);
\draw (Y1)--(T21);
\draw (Y1)--(Tn1);
\draw[very thick,red] (Y)--(X);
\draw (T11)--(X);
\draw (T1m)--(X);
\draw (Y)--(T2);
\draw (Y)--(Tn);
   \end{tikzpicture}
 }
\end{center}
Formally, the   partial order $\leq^{\hyper{H}}$ (or simply $\leq$, when $\hyper{H}$ is understood) is defined as the smallest partial order  generated  by the following rules:
\begin{center}
\resizebox{10cm}{!}{\fbox{$\begin{array}{c}
\seqtwo{\hyper{H},Y\leadsto K_1,\ldots, K_n\quad\quad \hyper{K}_1,X\leadsto H_{11},\ldots,H_{1m}}
{T_{11}:\hyper{H}_{11},\ldots,T_{1m}:\hyper{H}_{1m}\quad\quad T_2:\hyper{K}_2,\ldots,T_n:\hyper{K}_n}{Y(X(T_{11},\ldots,T_{1m}),T_2,\ldots,T_n)  \leq^{\hyper{H}}  (Y \union X)(T_{11},\ldots,T_{1m},T_2,\ldots,T_n)} \\\\
\seq{\hyper{H},Y\leadsto H_1,\ldots,H_n\quad\quad T_2:\hyper{H}_2,\ldots,T_n:\hyper{H}_n\quad\quad T_1 \leq^{\hyper{H}_1} T'_1}{ Y(T_1,T_2,\ldots,T_n)   \leq^{\hyper{H}}  Y(T'_1,T_2,\ldots,T_n)}
\end{array}$}}
\end{center}
This definition is well-formed, in the sense that, if $S:\hyper{H}$ and if $S\leq T$ is inferred, then $T:\hyper{H}$ can be inferred.
%
%
The one-node construct $H$ is  maximum, while
%
%
the constructions are the minimal elements (there is no $X\union Y$ to split).

\medskip
The  partial order $\leq$ admits two other equivalent definitions, for which we shall provisionally write $\leq^{\hyper{H}}_2$ and $\leq^{\hyper{H}}_3$ (shortly $\leq_2$ and $\leq_3$, respectively) before we  prove that they define the same relation as $\leq^{\hyper{H}}$.  The formulation $\leq_2$ will allow us to prove the equivalence of our definitions with the original ones of \cite{DP-HP}, while, the formulation $\leq_3$ underlies an algorithm for enumerating  all the vertices  inferior to a given construct (see Section \ref{vertices-of-faces-section}).

\smallskip
The  definition of $\leq_2$ is given by  two  clauses (guided by Lemma \ref{partition-of-i}):
\begin{center}
\resizebox{8.5cm}{!}{\fbox{$\begin{array}{c}
\seq{}{H\leq_2^{\hyper{H}} H} \\\\
\seqthree{Y\inc X\quad\quad \hyper{H},Y\leadsto K_1,\ldots,K_m\quad\quad \hyper{H},X\leadsto H_1,\ldots,H_n}
{S_1:\hyper{K}_1,\ldots,S_m:\hyper{K}_m\quad\quad T_1:\hyper{H}_1,\ldots,T_n:\hyper{H}_n}
{S_j\leq_2^{\hyper{K}_j} (K_j\inter X)\setc{T_i}{H_i\inc K_j}\; \mbox{for all}\;j}
{Y(S_1,\ldots,S_m)\leq_2^{\hyper{H}} X(T_1,\ldots,T_n)}
\end{array}$}}
\end{center}
The relation $\leq_2$  formalises the following intuition. As in the  definition of $\leq$, given a construct $S$, we want to know  which  constructs $T$  (of the same hypergraph) lie above  $S$ in the partial order.  If $S=H$, then $S$ is the maximum construct of $\hyper{H}$, hence $H\leq T$ boils down to $T=H$. Otherwise, the root of $S$ must be a subset of the root of $T$, and  
the task of showing $S\leq T$ is reduced  to that of verifying that each $S_j$ lies lower than an appropriate term.  

\medskip
For the definition of $\leq_3$, we need to introduce a variation of the notion of construct.
We define the {\em partial constructs} of   a connected hypergraph $\hyper{H}$  by adding one clause to the inductive definition of constructs:
\begin{itemize}
\item The single-node tree decorated with $\Omega_H$ is a partial construct of $\hyper{H}$.
\end{itemize}
(and by replacing   ``construct'' with   ``partial construct'' in the original clauses).

\smallskip
To distinguish partial constructs from constructs, we use the font $\mathbb{S},\mathbb{T},\ldots$ for the former.  We summarise   the definition of partial constructs as follows:
\begin{center}
\resizebox{11cm}{!}{\fbox{$
\seq{}{\Omega_H:\hyper{H}}\quad\quad\quad\quad
\seq{}{H:\hyper{H}}\quad\quad\quad\quad
\seq{\hyper{H},X\leadsto H_1,\ldots,H_n\quad\quad \mathbb{T}_1:\hyper{H}_1,\ldots,\mathbb{T}_n:\hyper{H}_n}
{X(\mathbb{T}_1,\ldots,\mathbb{T}_n):\hyper{H}}
$}}
\end{center}
We  define a partial construction to be  a partial construct in which all the non-$\Omega$ nodes are labelled by singletons. 

Note the difference  between decorations  $X$ and $\Omega_X$  in a partial construct: the latter stands for ``undefined'', in the spirit of Scott domain theory.  We shall write $\mathbb{T}[\Omega_X\leftarrow \mathbb{S}]$ for  the partial construct obtained  from $\mathbb{T}$ by replacing $\Omega_X$ with  $\mathbb{S}:\hyper{H}_X$. 

 We shall use the notation  $\mathbb{T} \blacktriangleright^{\hyper{H}} X$ to indicate that $X$ is the union of all   non-$\Omega$-decorations of $\mathbb{T}$, and we shall say that 
$\mathbb{T}$ {\em spans} $X$. 
Formally, this predicate is inductively defined as follows:
\begin{center}
\resizebox{11.75cm}{!}{\fbox{$\begin{array}{c}
\seq{}{\Omega_H\blacktriangleright^{\hyper{H}}\emptyset}\quad \enspace\enspace
\seq{}{H\blacktriangleright^{\hyper{H}} H}\quad  \enspace\enspace
\seq{\hyper{H},X\leadsto H_1,\ldots,H_n\quad\quad\mathbb{T}_1\blacktriangleright^{\hyper{H}_1}X_1,\ldots,
\mathbb{T}_n\blacktriangleright^{\hyper{H}_n}X_n}{X(\mathbb{T}_1,\ldots,\mathbb{T}_n)\blacktriangleright^{\hyper{H}}X\union X_1\union\ldots\union X_n}
\end{array}$}}
\end{center}

\begin{lemma} \label{occurences-of-Omega}
If $\mathbb{T}\blacktriangleright^{\hyper{H}} X$, with $\hyper{H},X\leadsto
 H_1,\ldots,H_n$, then, for each $i\in\set{1,\ldots,n}$, there exists exactly one occurrence of $\Omega_{H_i}$ in $\mathbb{T}$, and these are all the occurrences of an $\Omega$ in $\mathbb{T}$.
\end{lemma}
\Proof  
The  proof is   by structural induction on the proof of  well-formedness of $\mathbb{T}$. The case $\mathbb{T}=Y(\mathbb{T}_1,\ldots,\mathbb{T}_n)$ is settled by appealing to  Lemma \ref{partition-of-i}. \qed

\medskip
It follows  from this lemma that 
the partial constructs (resp. constructions) that span the whole carrier $H$ of a hypergraph $\hyper{H}$ are exactly the constructs (resp. constructions) of $\hyper{H}$.  
%
\begin{lemma}\label{T-X-less} With the notations of Lemma \ref{occurences-of-Omega}, 
 if $S_1,\ldots,S_n$ are constructs of $\hyper{H}_1,\ldots,\hyper{H}_n$, respectively, then, for all $1\leq i\leq n$, $\mathbb{T}[\ldots ,\Omega_{H_i}\leftarrow S_i,\ldots]$ is  a  construct  of $\hyper{H}$, and  $\mathbb{T}[\ldots, \Omega_{H_i}\leftarrow S_i,\ldots]\leq X(S_1,\ldots,S_n)$.
\end{lemma}
\Proof  By structural induction on $\mathbb{T}$.  
We treat  the case  $\mathbb{T}=Y(\mathbb{T}_1,\ldots,\mathbb{T}_m)$, with $\mathbb{T}_j$ spanning $Y_j$ for all $j$. Setting  $I_j=\setc{i}{\Omega_{H_i}\;\mbox{occurs in}\; \mathbb{T}_j}$ and  $T'_j$ to be the result of replacing each $\Omega_{H_i}$ by $S_i$ in $\mathbb{T}_j$  ($i$ ranging over $I_j$), we get by induction that $T'_j\leq  Y_j\setc{S_i}{i\in I_j}$, and we conclude as follows:
$$\begin{array}{llllllll}
\mathbb{T}[\ldots, \Omega_{H_i}\leftarrow S_i,\ldots] & = &Y(T'_1,\ldots,T'_m) \\
& \leq &
Y(Y_1\setc{S_i}{i\in I_1},\ldots,Y_m\setc{S_i}{i\in I_m})
\\& \leq & (Y\union\set{Y_1,\ldots,Y_m})(S_1,\ldots,S_n) \: =\: X(S_1,\ldots,S_n)\;. \mbox{\qed}
\end{array}$$
 
We have now all the prerequisites for  our third presentation of the partial order.  We define $\leq_3^{\hyper{H}}$ by the following two clauses:
\begin{center}
\resizebox{10cm}{!}{\fbox{$\begin{array}{c}
\seq{S:\hyper{H}}{S\leq_3^{\hyper{H}} H} \\\\
\seqtwo{\mathbb{T}\blacktriangleright^{\hyper{H}} X\quad\quad  \hyper{H},X\leadsto H_1,\ldots,H_n}
{S_1:\hyper{H}_1,\ldots,S_n:\hyper{H}_n\quad T_1:\hyper{H}_1,\ldots,T_n:\hyper{H}_n\quad\quad S_i\leq_3^{\hyper{H}_i} T_i\;\;\mbox{for all}\; i }
{\mathbb{T}[\ldots, \Omega_{H_i}\leftarrow S_i,\ldots]\leq_3^{\hyper{H}} X(T_1,\ldots,T_n)}
\end{array}$}}
\end{center}
Unlike for 
$\leq$ and $\leq_2$, the algorithmic reading of $S\leq_3 T$ answers the question of when $S$ lies lower than some fixed $T$. If $T=H$, then any construct of $\hyper{H}$
 lies lower than $T$. Otherwise,  $S$  ``starts by spanning $X$'' (and  recursively so).

%

\begin{proposition}\label{equiv-def-partial-order}
The  relations $\leq$, $\leq_2$, and $\leq_3$ coincide.
\end{proposition}
\Proof  
That $S\leq T$  implies $S\leq_2 T$ is proved by showing that  $\leq_2$ is closed under the rules that define $\leq$, including reflexivity and transitivity. Let us look at transitivity. 
Suppose  that $\hyper{H},Y \leadsto K_1,\ldots,K_m$, $\hyper{H},X \leadsto H_1,\ldots,H_n$, $\hyper{H},Z\leadsto G_1,\ldots,G_k$, and
 $$ Y(S_1,\dots,S_m)\leq_2X(T_1,\dots,T_n) 
\leq_2 Z(U_1,\dots,U_k)\:.
$$
We discuss only the case 
 where $m,n,k\geq 1$.  
We have to show that the two conditions allowing to deduce  $Y(S_1,\dots,S_m)\leq_2 Z(U_1,\dots,U_k)$ hold.
Collecting the  first  conditions in clause 2  of $\leq_2$, relative to our present two assumptions, we have that $Y\inc X$ and $X\inc Z$, and hence $Y\inc Z$. We now show that the second condition holds. Let us fix $j\in\{1,\dots,m\}$ and let $$Y_j=K_j\cap Z,\enspace\enspace Y_j'= K_j\cap X \enspace\mbox{ and } \enspace X_i= H_i\cap Z.$$ 
We have to prove that  $$S_j\leq_2 Y_j(\{U_l\,|\, l\in L_j \}),\;\mbox{where} \;L_j=\{l\in\{1,\dots,k\}\,|\, G_l\subseteq K_j\}\;.$$ 
The second condition for our first assumption gives us that, for all $j$:
$$S_j\leq_2 Y_j'(\{T_i\,|\, i\in I_j\}), \;\mbox{where}\;
I_j=\{i\in\{1,\dots,n\}\,|\, H_i\subseteq  K_j\}\;.$$ 
Now, for each $T_i$, where $i\in I_j$,  the second condition for the second assumption gives us that \begin{equation}
T_i\leq_2 X_i(\{U_m\,|\,m\in M_i\}),
\end{equation} where $M_i=\{m\in\{1,\dots,k\}\,|\, G_m\subseteq  H_i\}$.
Next, we have that
\vspace{-.05cm}
 \begin{equation}
\begin{array}{rcl}
Y_j & = &
 ( K_j\cap X)\cup (K_j\cap (Z\backslash X)) \\
 & = & Y_j' \cup ((K_j\backslash X)\cap Z)\\
 & = & Y_j' \cup ((\bigcup_{i\in I_j}H_i)\cap Z)\\
 & = & Y_j' \cup \bigcup_{i\in I_j}(H_i\cap Z)\\
 & = & Y_j' \cup \bigcup_{i\in I_j}X_i.
\end{array} \end{equation} 
\vspace{-.1cm}
And, lastly, since $K_j\backslash Z=(\bigcup_{i\in I_j} H_i)\backslash Z$, we have that 
\vspace{-.05cm}
\begin{equation}
\begin{array}{rcl}
L_j & = & \{l\in\{1,\dots,k\}\,|\, G_l\subseteq  K_j\} \\
    & = & \{l\in\{1,\dots,k\}\,|\, G_l\subseteq (\bigcup_{i\in I_j}H_i)\}\\
    & = & \bigcup_{i\in I_j} \{l\in\{1,\dots,k\}\,|\, G_l\subseteq H_i\}\\
    & = &\bigcup_{i\in I_j} M_i.
\end{array} \end{equation} 

\vspace{-.15cm}
Finally, (1), (2), (3) and (4), together with the rules from the definition of $\leq$, give us that 
\vspace{-.1cm}
$$\begin{array}{rclr}
S_j & \leq_2 & Y_j'(\{T_i\,|\, i\in I_j\}) & \hspace{0.35cm} (1)\\
 & \leq_2 & Y_j'(\{X_i(\{U_m\,|\,m\in M_i\})\,|\, i\in I_j\}) & \hspace{0.35cm} (2), \mbox{congruence}\\
 & \leq_2 & (Y_j'\cup \bigcup_{i\in I_j}X_i) (\{U_m\,|\,m\in M_i \mbox{ and } i\in I_j\}) & \hspace{0.35cm}  \mbox{axiom of } \leq \\
 & = & Y_j (\{U_m\,|\,m\in \bigcup_{i\in I_j}M_i \}) & \hspace{0.35cm} (3)\\
 & = & Y_j (\{U_l\,|\,l\in L_j \}). & \hspace{0.35cm} (4)
\end{array}$$

\vspace{-.15cm}
Note that this proof is valid provided one has shown beforehand that $\leq_2$ is closed under the other defining clauses of $\leq$.

\smallskip
That $S\leq_2 T$ implies $S\leq_3 T$  (resp. $S\leq_3 T$ implies $S\leq T$) is proved by induction on the proof of $S\leq_2 T$ (resp. of $S\leq_3 T$).  \qed

\subsection{Examples of hypergraphs and constructs} 
\label{construct-examples-section}

In this section, we provide a few examples of hypergraphs and their constructs, conveying  an intuitive understanding of their geometric realisation. 
We shall freely write $x$ instead of $\set{x}$ etc. for the labels of singleton nodes of constructs.

\smallskip As our first example, we describe the $n\!\!-\!\!1$-dimensional simplex:
$$\hyper{H}=\set{\set{x_1},\ldots,\set{x_n}, \set{x_1,\ldots,x_n}}.$$
All of its constructs   have the form $X(y_1,\ldots,y_p),$ 
where $X\inc \set{x_1,\ldots,x_n}$ and $\set{y_1,\ldots,y_p}= \set{x_1,\ldots,x_n}\setminus X$.
Note that $\hyper{H}$ is indeed a hypergraph, and not just the discrete graph with $n$ vertices, because we insist that the hyperedge $ \set{x_1,\ldots,x_n}$ is included.
\begin{itemize}
\item At dimension 2 	and writing $x,y,z$ instead of $x_1,x_2,x_3$,  we have 3 vertices, 3 edges or facets and the maximum face:
 $$\begin{array}{rll}
\mbox{vertices} && x(y,z)\quad y(x,z)\quad z(x,y) \\   \
\mbox{facets} && \set{x,y}(z)\quad \set{y,z}(x)\quad \set{x,z}(y)  \\
\mbox{whole polytope} && \set{x,y,z}.
\end{array}$$
 Note that $x(y,z)\leq \set{x,y}(z)$ and $y(x,z)\leq \set{x,y}(z)$, which says  combinatorially   that the edge $ \set{x,y}(z)$  connects the  vertices $x(y,z)$ and $y(x,z)$.
\item At dimension 3, we get  4 vertices, 6 edges and 4 facets.
\end{itemize}
We illustrate now how the hypergraph structure allows us to make truncations.  The desired effects of truncation will be obtained by adding  hyperedges  to the bare ``simplex hypergraph''. 
\begin{itemize}
\item Truncation of a vertex, say $x(y,z)$, of the 2-dimensional simplex (cf. Section \ref{intro-section}).  We add the hyperedge $\set{y,z}$ 
to the simplex hypergraph:
$$\hyper{H}=\set{\set{x},\set{y},\set{z},\set{y,z}, \set{x,y,z}}\;.$$
Then $x(y,z)$ is not a construction anymore, since $\hyper{H}_{\set{y,z}}$ is now connected.  Instead, we have 3 new constructs (encoding two  vertices and one edge):
$$x(y(z)) \quad x(z(y)) \quad x(\set{y,z})\;.$$
\item Truncation of an edge, say $\set{x,y}(u,z)$, of the 3-dimensional simplex.  Similarly, we add the hyperedge $\set{u,z}$ to the simplex hypergraph:
$$\hyper{H}=\set{\set{x},\set{y},\set{z},\set{u},\set{u,z} ,\set{x,y,z,u}}\; .$$
The edge  $\set{x,y}(u,z)$ and its end vertices are now replaced by a rectangular face (9 new constructs):
$$\begin{array}{ll}
  x(y,u(z))\quad x(y,z(u))\quad y(x,u(z))\quad y(x,z(u))\\
 x(y,\set{u,z})\quad y(x,\set{u,z})\quad  \set{x,y}(u(z))\quad \set{x,y}(z(u))\\
 \set{x,y}(\set{u,z}) \;.
\end{array}$$
\begin{center}
\begin{tikzpicture}[thick,scale=3]
\coordinate (A11) at (0.285,0.8);
\coordinate (A12) at (-0.285,0.8);
\coordinate (A21) at (0.285,-0.8); 
\coordinate (A22) at (-0.285,-0.8); 
\node (a22) at (-0.45,-0.925) {\small $x(y,z(u))$}; 
\node (a21) at (0.45,-0.925) {\small $x(y,u(z))$}; 
\node (a11) at (-0.45,0.925) {\small $y(x,z(u))$}; 
\node (a12) at (0.45,0.925) {\small $y(x,u(z))$}; 
\node (b12) at (-0.55,0.075) {\small $\set{x,y}(z(u))$}; 
\node (b22) at (0.55,-0.075) {\small $\set{x,y}(u(z))$}; 
\node (b23) at (0,0.7) {\small $y(x,\set{u,z})$}; 
\node (b24) at (0,-0.7) {\small $x(y,\set{u,z})$}; 
\node (b25) at (0,0) {\small $\set{x,y}(\!\set{u,z}\!)$}; 

\coordinate (A3) at (-1,0);
\coordinate (A4) at (1,0);
\draw[draw=gray,fill=pur,opacity=0.15] (A11) -- (A12) -- (A3) -- (A22) -- (A21) -- (A4) -- cycle;
\draw[draw=gray] (A11) -- (A12) -- (A3) -- (A22) -- (A21) -- (A4) -- cycle;
\draw[draw=gray,fill=yellow,opacity=0.1] (A11) -- (A12) -- (A22) -- (A21) -- cycle;
\draw[draw=gray] (A12) -- (A22);
\draw[draw=gray] (A11) -- (A21);
\draw[draw=gray,dashed,opacity=0.5] (A3) -- (-0.285,0);
\draw[draw=gray,dashed,opacity=0.5] (A4) -- (0.285,0);
\end{tikzpicture}
\end{center}
\item Truncation of a vertex, say $x(y,z,u)$, of the 3-dimensional simplex.  We achieve this by adding the hyperedge $\set{y,z,u}$ to the simplex hypergraph:
$$\hyper{H}=\set{\set{x},\set{y},\set{z},\set{u},\set{y,z,u} ,\set{x,y,z,u}}\;.$$
 This hypergraph disallows the construction $x(y,z,u)$  since
$\restrH{H}{\{x\}}$ is now connected, and replaces it by 3 vertices, 3 edges, and a facet:
$$\begin{array}{l}
x(y(z,u))\quad x(z(y,u))\quad x(u(y,z))\\
 x(\set{y,z}(u))\quad x(\set{z,u}(y))\quad x(\set{u,y}(z))\\
x(\set{y,z,u})\;.
\end{array}$$
\end{itemize}
\begin{center}
\begin{tikzpicture}[thick,scale=3]
\coordinate (A11) at (0,0.2);
\coordinate (A12) at (-0.45,0.6);
\coordinate (A13) at (0.45,0.6);
\node (a12) at (-0.55,0.7) {\small $x(y(z,u))$};
\node (a13) at (0.6,0.7) {\small $x(z(y,u))$};
\node (a1r1) at (0.25,0.125) {\small $x(u(y,z))$};
\node (b1) at (-0.42,0.285) {\small $x(\set{u,y}(z))$};
\node (b2) at (0.42,0.285) {\small $x(\set{z,u}(y))$};
\node (b3) at (0,0.7) {\small $x(\set{y,z}(u))$};
\node (c) at (0,0.475) {\small $x(\set{y,z,u})$};
\coordinate (A2) at (0,-1); 
\coordinate (A3) at (-1,0);
\coordinate (A4) at (1,0);
\draw[draw=gray,fill=pur,opacity=0.15] (A3) -- (A2) -- (A4) -- (A13) -- (A12) -- (A3) -- cycle;
\draw[draw=gray] (A3) -- (A2) -- (A4) -- (A13) -- (A12) -- (A3) -- cycle;
\draw[draw=gray] (A11) -- (A12) -- (A13) -- cycle;
\draw[draw=gray,fill=yellow,opacity=0.1] (A11) -- (A12) -- (A13) -- cycle;
\draw[draw=gray] (A11) -- (A2) -- cycle;
\draw[draw=gray,dashed,opacity=0.5] (A3) -- (A4);

\end{tikzpicture}
\end{center}

\indent Our next example is the family of associahedra.  One of the  standard labellings  of the faces of the $n$-dimensional associahedron is by  all the (partial or total) parenthesisations of a word of $n+2$ letters.  Here, the idea  is to focus, not on the letters (or the leaves of the corresponding tree), but on the $n+1$  "compositions of these letters" involved.  These compositions are next to each other, as suggested in the following picture for dimension 3 (where $a,b,c,\ldots$ are the letters and $x,y,\ldots$ are the compositions):
\begin{equation}\label{blah}\begin{array}{cccccccccccccccc}
a && b&&c&&d&&e\\
&x&&y&&z&&u
\end{array}\end{equation}
\begin{itemize}
\item  At dimension 2, this suggests to take the following graph (in hypergraph form), expressing  ``$x$ is next to $y$ which is next to $z$'':
$$\hyper{H}=\set{\set{x},\set{y},\set{z},\set{x,y},\set{y,z}}\;.$$
(Note that the hyperedge $\set{x,y,z}$ is no longer necessary to ensure that $\hyper{H}$ is connected.)
The edges $\set{x,y}$ and $\set{y,z}$  are prescriptions for truncating two vertices of a triangle, yielding a pentagon.  The 5 vertices are
$$x(y(z)) \quad x(z(y)) \quad y(x,z)\quad  z(x(y))\quad z(y(x))\;.$$
\item At dimension 3, we take 
$$\hyper{H}=\set{\set{x},\set{y},\set{z},\linebreak\set{u},\set{x,y},\set{y,z},\set{z,u}}\;,$$
which seems like a prescription for truncating  (only)  three edges of the simplex.  But look at what has become of the vertex $u(x,y,z)$. It has been also truncated!  Indeed, it has been split into 5 constructions (with corresponding edges and face):
\vspace{-.05cm}
$$u(x(y(z)))\quad u(x(z(y)))\quad u( y(x,z))\quad u(z(x(y))\quad  u(z(y(x)))\;.$$
\vspace{-.05cm}

To build these constructions, we have used that  $\hyper{H}_{\set{x,y,z}}$ is connected.  In fact,  the truncation prescriptions  are all hyperedges of ${\it Sat}(\hyper{H})$.
\item At dimension $n$, we take
\vspace{-.05cm}
$$\hyper{H}=\set{\set{x_1},\ldots, \set{x_{n+1}},\set{x_1,x_2},\ldots,\set{x_n,x_{n+1}}}\; .$$
\vspace{-.05cm}
\end{itemize}
\vspace{-.3cm}

\indent
Here is the recipe showing how to move  between three equivalent presentations of the faces of the associahedra: partially parenthesised words, rooted  (undecorated) planar trees (see e.g.  \cite{LV-AO}), and constructs.
\begin{itemize}
\item[-] From rooted planar trees to constructs.  Label all the intervals beween the leaves of a tree with $n+2$ leaves  by
$x_1,\ldots,x_{n+1}$ (from left to right). Consider the $x_i$'s as balls and let them fall. Label each node of the tree by the set of balls which fall to that node.  Finally, remove all the leaves. For example:
\vspace{-.05cm}
\begin{center}
\resizebox{9cm}{!}{\begin{tikzpicture}
   \node (Y) [circle,draw=none,minimum size=0.75cm,inner sep=0.1mm] at (-0.5,0) {$a(bc)d$};
   \draw [->,decorate,decoration={snake,amplitude=.4mm,segment length=2mm,post length=1mm}]
(0.5,0) -- (1.5,0);
 \draw [->,decorate,decoration={snake,amplitude=.4mm,segment length=2mm,post length=1mm}]
(6.5,0) -- (7.5,0);

\coordinate (A) at (4,-1);
\coordinate (C) at (4,-0);
\coordinate (B1) at (2,1);
\node at (B1) [yshift=0.15cm] { $a$};
\coordinate (B2) at (3.25,1);
\node at (B2) [yshift=0.15cm] { $b$};
\coordinate (B3) at (4.75,1);
> \node at (B3) [yshift=0.15cm] {$c$};
\coordinate (B4) at (6,1);
\node at (B4) [yshift=0.18cm] { $d$};
\node (x) at (2.85,0.75) { $x$};
\node (z) at (5.15,0.75) {$z$};
\node (y) at (4,0.75) {$y$};
\draw (A)--(B1);
\draw[very thick] (A)--(C);
\draw (C)--(B3);
\draw (C)--(B2);
\draw (A)--(B4);


\node (r) at (8.5,-0.5) { $\{x,z\}$};
\node (l) at (8.5,0.5) {$y$};
\draw (r)--(l);
   \end{tikzpicture}
 }
\end{center}
\vspace{-.05cm}
\item[-] From constructs to parenthesisations. Read a construct from the leaves to the root, and each node as an instruction for building a parenthesis.  If the label is, for example, $\set{x_i,x_{i+2}}$, then the instruction is to
do an unbiased composition  of three  partially parenthesised words $w_1,w_2,w_3$
``above" $x_i$ and $x_{i+2}$ 
 (for example, $a,(bc), d$  are above $x$ and $z$ in (\ref{blah})),  in one shot, resulting in $(w_1w_2w_3)$.
\item[-] From parenthesisations to trees.  This is standard.
\end{itemize}
\noindent For the 2-dimensional associahedron,  the   representation  with  planar rooted trees / constructs is  given on the next picture:
\vspace{-.05cm}
\begin{center}
\resizebox{8cm}{!}{\begin{tikzpicture}
    \node (E) at (0,0) { $y(x,z)$};
    \node (F) at (-2.7,-1.7) { $x(y(z))$};
    \node (A) at (2.7,-1.7) { $z(y(x))$};
    \node (Asubt) at (-1.7,-3.9) {$x(z(y))$};
    \node (P4) at (1.7,-3.9) { $z(x(y))$};
  \node (P41) at (3.8,-1.7) {\resizebox{0.5cm}{!}{\begin{tikzpicture}
\coordinate (L) at (4,-0.25);
\coordinate (M2) at (3.745,0.07);
\coordinate (M1) at (3.835,-0.05);
\coordinate (B1) at (3.6,0.25);
 \coordinate (B2) at (4.05 ,0.25);
 \coordinate (B3) at (3.87,0.25);
 \coordinate (B4) at (4.4,0.25);
 \draw[thick] (M1)--(B2);
\draw[thick] (L)--(B1);
\draw[thick] (M2)--(B3);
\draw[thick] (L)--(B4);\end{tikzpicture}}};

  \node (E1) at (0,0.6)  {\resizebox{0.5cm}{!}{\begin{tikzpicture}
\coordinate (L) at (4,-0.25);
\coordinate (C) at (3.8 ,-0);
\coordinate (M) at (4.2 ,-0);
\coordinate (B1) at (3.6,0.25);
 \coordinate (B2) at (3.95,0.25);
 \coordinate (B3) at (4.05,0.25);
 \coordinate (B4) at (4.4,0.25);
 \draw[thick]  (L)--(B1);
\draw[thick] (L)--(C);
\draw [thick] (C)--(B2);
\draw[thick] (M)--(B3);
\draw [thick] (L)--(B4);\end{tikzpicture}}};
  \node (E1) at (-3.8,-1.7) {\resizebox{0.5cm}{!}{\begin{tikzpicture}
\coordinate (L) at (4,-0.25);
\coordinate (M2) at (4.25,0.05);
\coordinate (M1) at (4.15 ,-0.05);
\coordinate (B1) at (3.6,0.25);
 \coordinate (B2) at (3.925 ,0.25);
 \coordinate (B3) at (4.1,0.25);
 \coordinate (B4) at (4.4,0.25);
 \draw[thick] (M1)--(B2);
\draw[thick] (L)--(B1);
\draw[thick] (M2)--(B3);
\draw[thick] (L)--(B4);\end{tikzpicture}}};
\node (Asubt1) at (2.6,-4.4) {\resizebox{0.5cm}{!}{\begin{tikzpicture}
\coordinate (L) at (4,-0.25);
\coordinate (M2) at (3.93,0.07);
\coordinate (M1) at (3.835,-0.05);
\coordinate (B1) at (3.6,0.25);
 \coordinate (B2) at (4.05 ,0.25);
 \coordinate (B3) at (3.8,0.25);
 \coordinate (B4) at (4.4,0.25);
 \draw[thick] (M1)--(B2);
\draw[thick] (L)--(B1);
\draw[thick] (M2)--(B3);
\draw[thick] (L)--(B4);\end{tikzpicture}}};

  \node (A1d) at (-1.95,-0.3) {\color{violet}\resizebox{0.5cm}{!}{\begin{tikzpicture}
\coordinate (L) at (4,-0.25);
\coordinate (C) at (3.8 ,-0);
\coordinate (M) at (4.2 ,-0);
\coordinate (B1) at (3.6,0.25);
 \coordinate (B2) at (3.95,0.25);
 \coordinate (B3) at (4.05,0.25);
 \coordinate (B4) at (4.4,0.25);
 \draw[thick] (L)--(B2);
\draw [thick] (L)--(B1);
\draw[thick] (M)--(B3);
\draw[thick] (L)--(B4);\end{tikzpicture}}};

  \node (A1dd) at (1.95,-0.3) {\color{violet}\resizebox{0.5cm}{!}{\begin{tikzpicture}
\coordinate (K) at (4,-0.25);
\coordinate (C) at (3.8 ,-0);
\coordinate (B1) at (3.6,0.25);
 \coordinate (B2) at (4,0.25);
 \coordinate (B3) at (4.2,0.25);
 \coordinate (B4) at (4.4,0.25);
 \draw[thick] (K)--(B1);
\draw[thick] (K)--(C);
\draw[thick] (C)--(B2);
\draw[thick] (K)--(B3);
\draw[thick] (K)--(B4);\end{tikzpicture}}};
  \node (fr) at (3.25,-3)  {\color{violet}\resizebox{0.5cm}{!}{\begin{tikzpicture}
\coordinate (Ae) at (4,-0.25);
\coordinate (Ce) at (3.8 ,-0);
\coordinate (B1e) at (3.6,0.25);
 \coordinate (B2e) at (3.8,0.25);
 \coordinate (B3e) at (4,0.25);
 \coordinate (B4e) at (4.4,0.25);
 \draw[thick] (Ae)--(B1e);
\draw[thick] (Ae)--(Ce);
\draw[thick] (Ce)--(B3e);
\draw[thick] (Ce)--(B2e);
\draw[thick] (Ae)--(B4e);\end{tikzpicture}}};

 \node (Cs1) at (0,-2.4) {\color{blue}\resizebox{0.5cm}{!}{\begin{tikzpicture}
\coordinate (Ls) at (4,-0.25);
\coordinate (Cs) at (3.8 ,-0);
\coordinate (Ms) at (4.2 ,-0);
\coordinate (B1s) at (3.6,0.25);
 \coordinate (B2s) at (3.875,0.25);
 \coordinate (B3s) at (4.125,0.25);
 \coordinate (B4s) at (4.4,0.25);
 \draw (Ls)--(B2s);
\draw (Ls)--(B1s);
\draw (Ls)--(B3s);
\draw (Ls)--(B4s);\end{tikzpicture}}};

\node (A1dsddd) at (-3.25,-3)   {\color{violet}\resizebox{0.5cm}{!}{\begin{tikzpicture}
\coordinate (Lq) at (4,-0.25);
\coordinate (Cq) at (3.8 ,-0);
\coordinate (Mq) at (4.2 ,-0);
\coordinate (B1q) at (3.6,0.25);
 \coordinate (B2q) at (4,0.25);
 \coordinate (B3q) at (4.2,0.25);
 \coordinate (B4q) at (4.4,0.25);
 \draw[thick] (Mq)--(B2q);
\draw[thick] (Lq)--(B1q);
\draw[thick] (Mq)--(B3q);
\draw[thick] (Lq)--(B4q);\end{tikzpicture}}};

\node (A1dddd) at (0,-4.4)  { \color{violet}\resizebox{0.5cm}{!}{\begin{tikzpicture}
\coordinate (Aw) at (4,-0.25);
\coordinate (Cw) at (4,-0);
\coordinate (B1w) at (3.6,0.25);
 \coordinate (B2w) at (3.8,0.25);
 \coordinate (B3w) at (4.2,0.25);
 \coordinate (B4w) at (4.4,0.25);
 \draw[thick] (Aw)--(B1w);
\draw[thick] (Aw)--(Cw);
\draw[thick] (Cw)--(B3w);
\draw[thick] (Cw)--(B2w);
\draw[thick] (Aw)--(B4w);\end{tikzpicture}}};

  \node (A1) at (-2.6,-4.4) {\resizebox{0.5cm}{!}{\begin{tikzpicture}
\coordinate (L) at (4,-0.25);
\coordinate (M2) at (4.075,0.05);
\coordinate (M1) at (4.15 ,-0.05);
\coordinate (B1) at (3.6,0.25);
 \coordinate (B2) at (3.925 ,0.25);
 \coordinate (B3) at (4.225,0.25);
 \coordinate (B4) at (4.4,0.25);
 \draw[thick] (M1)--(B2);
\draw[thick] (L)--(B1);
\draw[thick] (M2)--(B3);
\draw[thick] (L)--(B4);\end{tikzpicture}}};
    \draw[-] (E)--(F) node [midway,above,xshift=0.65cm,yshift=-0.5cm] {\color{violet}\scriptsize   $\{y,x\}(z)$};
    \draw[-] (E)--(A) node [midway,above,xshift=-0.5cm,yshift=-0.5cm] {\color{violet}\scriptsize     $\{y,z\}(x)$};
 \draw[-] (F)--(Asubt) node [midway,right] {\scriptsize\color{violet}     $x(\{y,z\})$};
    \draw[-] (A)--(P4) node [midway,left]  {\scriptsize\color{violet}   $z(\{x,y\})$};
    \draw[-] (Asubt)--(P4) node [midway,above] {\scriptsize\color{violet}     $\{x,z\}(y)$};
 \node (C) at (0,-1.95) {\color{blue}\scriptsize $\{x,y,z\}$};
   \end{tikzpicture}}
 \end{center}

Our final example is the family of permutohedra.  Here we take the complete graph on the set of vertices as the hypergraph. We discuss directly the general case at dimension $n$:
$$\hyper{H}=\set{\set{x_1},\ldots, \set{x_{n+1}}} \union\setc{\set{x_i,x_j}}{i,j\in \set{1,\ldots,n+1} \;\mbox{and}\; i\neq j}.$$
Note that  all the constructs of the permutohedra  are filiform, i.e., are trees reduced to a branch.
The faces of the permutohedra have been described in the literature as surjections, and also as planar rooted trees with levels.  The three  representations are related as follows:
\begin{itemize}
\item[-] From trees  with levels to constructs.  Consider again the $x_i$'s as balls  being thrown in the successive intervals between the leaves, and let them fall.  Then we  form the construct $Y_1(Y_2(\ldots(Y_m)\ldots))$, where 
$Y_i$ is the collection of balls that fall to level $i$  (counting levels from the root). For example:
\begin{center}
\resizebox{8cm}{!}{\begin{tikzpicture}
\coordinate (A) at (0,-0);
\coordinate (B1) at (-0.85,1);
\coordinate (B2) at (0.85,1);
\coordinate (C1) at (-1.675,2);
\coordinate (C2) at (-0.35,2);
\coordinate (C3) at (0.35,2);
\coordinate (C4) at (1.675,2);
\draw (B1)--(A)--(B2);
\draw (B1)--(C1);
\draw (B1)--(C2);
\draw (B2)--(C3);
\draw (B2)--(C4);
\draw[dashed] (-1.675,1)--(1.675,1);
\draw[dashed] (-1.675,0)--(1.675,0);
\draw[dashed] (4.3,1)--(5.3,1);
\draw[dashed] (4.3,0)--(5.3,0);
\node (xz) at (5.65,1) { $x,z$};
\node (yr) at (5.45,0) {  $y$};
\node (x) at (-1,2.1) { $x$};
\node (z) at (1,2.1) { $z$};
\node (y) at (0,2.05) { $y$};
\draw [->,decorate,decoration={snake,amplitude=.4mm,segment length=2mm,post length=1mm}]
(2.5,1) -- (3.5,1);
\draw [->,decorate,decoration={snake,amplitude=.4mm,segment length=2mm,post length=1mm}]
(6.6,1) -- (7.6,1);
\node (e) at (9,1) {$y(\{x,z\})$};
   \end{tikzpicture}}
\end{center}
\item[-] A construct $Y_1(Y_2(\ldots(Y_m)\ldots))$ defines a  surjection from $\set{x_1,\ldots,x_{n+1}}$ to $\set{1,\ldots,m}$  mapping each $x$ to $i$, where $i$ is such that $x\in Y_i$.
\item[-] From surjections to trees with levels.  We refer to \cite{LR-P}.
\end{itemize}
For the 2-dimensional permutohedron,  the   representation  with  planar rooted trees with levels / constructs is  given on the next picture: \begin{center}
\resizebox{7.25cm}{!}{\begin{tikzpicture}
    \node (E) at (-1.4,0) {$x(y(z))$};
    \node (G) at (1.4,0) {$y(x(z))$};
    \node (F) at (2.2,-1.7) {$y(z(x))$};
    \node (A) at (-2.2,-1.7) {$x(z(y))$};
    \node (Asubt) at (-1.4,-3.5) {$z(x(y))$};
    \node (P4) at (1.4,-3.5) {$z(y(x))$};

  \node (E1) at (-2.2,0.4) {\resizebox{0.5cm}{!}{\begin{tikzpicture}
\coordinate (L) at (4,-0.25);
\coordinate (M2) at (4.25,0.05);
\coordinate (M1) at (4.15 ,-0.05);
\coordinate (B1) at (3.6,0.25);
 \coordinate (B2) at (3.925 ,0.25);
 \coordinate (B3) at (4.1,0.25);
 \coordinate (B4) at (4.4,0.25);
 \draw[thick] (M1)--(B2);
\draw[thick] (L)--(B1);
\draw[thick] (M2)--(B3);
\draw[thick] (L)--(B4);\end{tikzpicture}}};

 \node (G1) at (2.2,0.4) {\resizebox{0.5cm}{!}{\begin{tikzpicture}
\coordinate (L) at (4,-0.25);
\coordinate (C) at (3.88 ,-0.1);
\coordinate (M) at (4.26 ,0.15);
\coordinate (B1) at (3.725,0.1);
 \coordinate (B2) at (4,0.1);
 \coordinate (B3) at (4.125,0.35);
 \coordinate (B4) at (4.4,0.35);
 \draw[thick] (C)--(B1);
\draw[thick] (L)--(C);
\draw[thick] (C)--(B2);
\draw[thick] (M)--(B3);
\draw[thick] (L)--(B4);\end{tikzpicture}}};

  \node (F1) at (3.2,-1.7) {\resizebox{0.5cm}{!}{\begin{tikzpicture}
\coordinate (L) at (4.19,-0.25);
\coordinate (C) at (3.88 ,0.15);
\coordinate (M) at (4.28 ,-0.1);
\coordinate (B1) at (3.725,0.35);
 \coordinate (B2) at (4,0.35);
 \coordinate (B3) at (4.125,0.1);
 \coordinate (B4) at (4.4,0.1);
 \draw[thick] (C)--(B1);
\draw[thick] (L)--(C);
\draw[thick] (C)--(B2);
\draw[thick] (M)--(B3);
\draw [thick] (L)--(B4);\end{tikzpicture}}};

  \node (A1) at (-3.2,-1.7) {\resizebox{0.5cm}{!}{\begin{tikzpicture}
\coordinate (L) at (4,-0.25);
\coordinate (M2) at (4.075,0.05);
\coordinate (M1) at (4.15 ,-0.05);
\coordinate (B1) at (3.6,0.25);
 \coordinate (B2) at (3.925 ,0.25);
 \coordinate (B3) at (4.225,0.25);
 \coordinate (B4) at (4.4,0.25);
 \draw[thick] (M1)--(B2);
\draw[thick] (L)--(B1);
\draw[thick] (M2)--(B3);
\draw[thick] (L)--(B4);\end{tikzpicture}}};

\node (Asubt1) at (-2.2,-3.9) {\resizebox{0.5cm}{!}{\begin{tikzpicture}
\coordinate (L) at (4,-0.25);
\coordinate (M2) at (3.93,0.07);
\coordinate (M1) at (3.835,-0.05);
\coordinate (B1) at (3.6,0.25);
 \coordinate (B2) at (4.05 ,0.25);
 \coordinate (B3) at (3.8,0.25);
 \coordinate (B4) at (4.4,0.25);
 \draw[thick] (M1)--(B2);
\draw[thick] (L)--(B1);
\draw[thick] (M2)--(B3);
\draw[thick] (L)--(B4);\end{tikzpicture}}};

  \node (P41) at (2.2,-3.9) {\resizebox{0.5cm}{!}{\begin{tikzpicture}
\coordinate (L) at (4,-0.25);
\coordinate (M2) at (3.745,0.07);
\coordinate (M1) at (3.835,-0.05);
\coordinate (B1) at (3.6,0.25);
 \coordinate (B2) at (4.05 ,0.25);
 \coordinate (B3) at (3.87,0.25);
 \coordinate (B4) at (4.4,0.25);
 \draw[thick] (M1)--(B2);
\draw[thick] (L)--(B1);
\draw[thick] (M2)--(B3);
\draw[thick] (L)--(B4);\end{tikzpicture}}};

 \node (C) at (0,-1.5) {\color{blue}\scriptsize $\{x,y,z\}$};
 \node (C1) at (0,-1.95) {\color{blue}\resizebox{0.5cm}{!}{\begin{tikzpicture}
\coordinate (L) at (4,-0.25);
\coordinate (C) at (3.8 ,-0);
\coordinate (M) at (4.2 ,-0);
\coordinate (B1) at (3.6,0.25);
 \coordinate (B2) at (3.875,0.25);
 \coordinate (B3) at (4.125,0.25);
 \coordinate (B4) at (4.4,0.25);
 \draw (L)--(B2);
\draw (L)--(B1);
\draw (L)--(B3);
\draw (L)--(B4);\end{tikzpicture}}};
    \draw (E)--(G) node [midway,below] {\color{violet}\scriptsize   $\{x,y\}(z)$};

 \draw[draw=none] (E1)--(G1) node [midway,above,yshift=-0.25cm] {\color{violet}\resizebox{0.5cm}{!}{\begin{tikzpicture}
\coordinate (L) at (4,-0.25);
\coordinate (C) at (3.8 ,-0);
\coordinate (M) at (4.2 ,-0);
\coordinate (B1) at (3.6,0.25);
 \coordinate (B2) at (3.95,0.25);
 \coordinate (B3) at (4.05,0.25);
 \coordinate (B4) at (4.4,0.25);
 \draw[thick] (L)--(B2);
\draw [thick] (L)--(B1);
\draw[thick] (M)--(B3);
\draw[thick] (L)--(B4);\end{tikzpicture}}};
  
  \draw (G)--(F) node [midway,left] {\color{violet}\scriptsize    $y(\{x,z\})$};

  \draw[draw=none] (G1)--(F1) node [midway,right,xshift=-0.35cm] {\color{violet}\resizebox{0.5cm}{!}{\begin{tikzpicture}
\coordinate (L) at (4,-0.25);
\coordinate (C) at (3.8 ,-0);
\coordinate (M) at (4.2 ,-0);
\coordinate (B1) at (3.6,0.25);
 \coordinate (B2) at (3.95,0.25);
 \coordinate (B3) at (4.05,0.25);
 \coordinate (B4) at (4.4,0.25);
 \draw[thick]  (L)--(B1);
\draw[thick] (L)--(C);
\draw [thick] (C)--(B2);
\draw[thick] (M)--(B3);
\draw [thick] (L)--(B4);\end{tikzpicture}}};
   
 \draw (E)--(A) node [midway,right] {\color{violet}\scriptsize   $x(\{y,z\})$};

 \draw[draw=none] (E1)--(A1) node [midway,left,xshift=0.35cm] {\color{violet}\resizebox{0.5cm}{!}{\begin{tikzpicture}
\coordinate (L) at (4,-0.25);
\coordinate (C) at (3.8 ,-0);
\coordinate (M) at (4.2 ,-0);
\coordinate (B1) at (3.6,0.25);
 \coordinate (B2) at (4,0.25);
 \coordinate (B3) at (4.2,0.25);
 \coordinate (B4) at (4.4,0.25);
 \draw[thick] (M)--(B2);
\draw[thick] (L)--(B1);
\draw[thick] (M)--(B3);
\draw[thick] (L)--(B4);\end{tikzpicture}}};

 \draw (F)--(P4) node [midway,left] {\color{violet}\scriptsize   $\{y,z\}(x)$};

 \draw[draw=none] (F1)--(P41) node [midway,right,xshift=-0.35cm] {\color{violet}\resizebox{0.5cm}{!}{\begin{tikzpicture}
\coordinate (K) at (4,-0.25);
\coordinate (C) at (3.8 ,-0);
\coordinate (B1) at (3.6,0.25);
 \coordinate (B2) at (4,0.25);
 \coordinate (B3) at (4.2,0.25);
 \coordinate (B4) at (4.4,0.25);
 \draw[thick] (K)--(B1);
\draw[thick] (K)--(C);
\draw[thick] (C)--(B2);
\draw[thick] (K)--(B3);
\draw[thick] (K)--(B4);\end{tikzpicture}}};

    \draw (A)--(Asubt) node [midway,right]  {\color{violet}\scriptsize    $\{x,z\}(y)$};
    \draw [draw=none] (A1)--(Asubt1) node [midway,left,xshift=0.35cm]  { \color{violet}\resizebox{0.5cm}{!}{\begin{tikzpicture}
\coordinate (A) at (4,-0.25);
\coordinate (C) at (4,-0);
\coordinate (B1) at (3.6,0.25);
 \coordinate (B2) at (3.8,0.25);
 \coordinate (B3) at (4.2,0.25);
 \coordinate (B4) at (4.4,0.25);
 \draw[thick] (A)--(B1);
\draw[thick] (A)--(C);
\draw[thick] (C)--(B3);
\draw[thick] (C)--(B2);
\draw[thick] (A)--(B4);\end{tikzpicture}}};
    \draw (Asubt)--(P4) node [midway,above] {\color{violet}\scriptsize   $z(\{x,y\})$};
  \draw [draw=none] (Asubt1)--(P41) node [midway,below,yshift=0.2cm] {\color{violet}\resizebox{0.5cm}{!}{\begin{tikzpicture}
\coordinate (A) at (4,-0.25);
\coordinate (C) at (3.8 ,-0);
\coordinate (B1) at (3.6,0.25);
 \coordinate (B2) at (3.8,0.25);
 \coordinate (B3) at (4,0.25);
 \coordinate (B4) at (4.4,0.25);
 \draw[thick] (A)--(B1);
\draw[thick] (A)--(C);
\draw[thick] (C)--(B3);
\draw[thick] (C)--(B2);
\draw[thick] (A)--(B4);\end{tikzpicture}}};

   \end{tikzpicture}

}
\end{center}
\subsection{Vertices of faces}  \label{vertices-of-faces-section}
As a preparation for the following section,  
given a hypergraph $\hyper{H}$ and a construct $T:\hyper{H}$, we give a device for finding all constructions $V$ such that $V\leq^{\hyper{H}} T$.
If $T$ is a construction, then this set is reduced to $T$  itself.  We shall use the notation $V \lessdot^{\hyper{H}}T$  for ``$V \leq^{\hyper{H}} T\;\mbox{and}\; V\;\mbox{is a construction}$''.

First, we notice that, by a straightforward tuning of the definition of $\leq_3$, the predicate $\lessdot$ is defined by the following clauses:
\begin{center}
\resizebox{7cm}{!}{\fbox{$\begin{array}{c} \seq{V\:\mbox{is a construction of}\: \hyper{H}}{V\lessdot^{\hyper{H}} H}\\\\
\seqtwo{\hyper{H},X\leadsto H_1,\ldots,H_n}{V_1 \lessdot^{\hyper{H}_1}  T_1 \quad \ldots \quad V_n \lessdot^{\hyper{H}_n}  T_n \quad\quad \mathbb{V}_0\blacktriangleright^{\hyper{H}} X}
{\mathbb{V}_0[\Omega_{H_1}\leftarrow V_1,\ldots,\Omega_{H_n}\leftarrow V_n]\lessdot^{\hyper{H}}  X(T_1,\ldots,T_n)}
\end{array}$}}
\end{center}
where, in the second clause,  $\mathbb{V}_0$ is a partial construction.

This suggests an algorithm.  For every node $X$ of $T$, we should ``zoom in" and replace it with a partial construction spanning $X$.  Here
 is a formal device for searching all the partial constructions of $\hyper{H}$ spanning a given fixed set $X\inc H$. One  starts from $\Omega_H$, and one performs rewriting  (non-deterministically), until exhaustion of $X$, as follows:
\begin{center}
\resizebox{5.25cm}{!}{\fbox{$\seq{\mathbb{V}\blacktriangleright^{\hyper{H}} Y\quad \quad Y\incs X\quad x\in X\setminus Y}{\mathbb{V}\longrightarrow_X \Sub{\mathbb{V}}{\Omega_K}{x(\Omega_{K_1},\ldots ,\Omega_{K_p})}}$}}
\end{center}
where $K$ is the connected component of $\restrH{H}{ Y}$ to which $x$ belongs and where
$\hyper{K},\set{x} \leadsto K_1,\ldots,K_p$.\\
\indent We write $\longrightarrow_X\!\!{}^*$ for the reflexive and transitive closure of $\longrightarrow_X$. 
We shall say that a partial construction $\mathbb{V}$ is {\em accepted}  if 
$\Omega_H\longrightarrow_X\!\!{}^* \:\mathbb{V}$
 and  there   exists no $\mathbb{V}'$ such that $\mathbb{V}\longrightarrow_X \mathbb{V}'$.
 As immediate observations, we have:
\begin{enumerate}
\item If $\mathbb{V}\longrightarrow_X \mathbb{V}'= \Sub{V}{\Omega_K}{x(\Omega_{K_1},\ldots ,\Omega_{K_p})}$, then
 $\mathbb{V}' \blacktriangleright^{\hyper{H}} Y\cup\{x\}$, i.e., $\mathbb{V}'$ is a partial construction spaning $Y\cup\{x\}$.
\item The  rewriting system $\longrightarrow_X$ is terminating, since the  cardinality of the spanned subset increases by 1 at each step, while remaining a subset of $X$. 
\end{enumerate}
\begin{lemma}\label{rewriting-partial-constructions} 
(1)
 The accepted partial constructions  are precisely the partial constructions spanning $X$.  
(2) For every element $x\in X$, there exists a partial construction spanning $X$  whose root is decorated by $x$.
 \end{lemma}
 \Proof If $\mathbb{V}$ is accepted, then  $\mathbb{V}$ spans $X$ by definition. 
Conversely, we proceed by induction on the cardinality of $X$. 
 We can write   $\mathbb{V}$  as $\mathbb{V}'[\Omega_K \leftarrow x(\Omega_{K_1},\ldots,\Omega_{K_p})]$, since every tree has a node all of whose outgoing edges are leaves, and then apply induction to $\mathbb{V}'$, which spans $X\setminus\set{x}$. 

As for the second claim, given $x\in X$, we can start the rewriting sequence    with $\Omega_H\longrightarrow_X x(\Omega_{K_1},\ldots, \Omega_{K_p})$. Then any continuation of this sequence leads to a partial construction spanning $X$, and has $x$ as a root. \qed

\medskip
Returning  to our  goal of finding all constructions $V$ such that $V\lessdot T$ (for  fixed $T$), we  transform our definition of $\lessdot$ into an algorithmic one  by replacing 
$$\mathbb{V}_0\blacktriangleright^{\hyper{H}} X\quad\mbox{with}\quad
\Omega_H\longrightarrow_X^*\mathbb{V}_0\blacktriangleright^{\hyper{H}} X$$
 in the  second clause: we apply the device repetitively at all nodes of $T$. 
\begin{corollary} \label{choice-of-construction-of-face}
For each construct $X(T_1,\ldots,T_n)$ and each $x\in X$, there exists at least one construction of the form $x(S_1,\ldots,S_m)$,  such that  $x(S_1,\ldots,S_m)\lessdot  X(T_1,\ldots,T_n)$.
\end{corollary}

\section{Constructs as geometric faces} \label{geometric-realisation-section}

In this section, we   recall the geometric realisation of hypergraph polytopes, following Do\v sen and Petri\'c, and we  provide  a new proof of their  theorem stating that the poset of constructs is isomorphic to the poset of geometric faces. The original proof in \cite{DP-HP} relies on Birkhoff's representation theorem, without providing an explicit description of the isomorphism. Our proof  is constructive, in that it  exhibits the isomorphism.

We first  prove the equivalence between our notion of constructs and theirs. Then we recall  the geometric realisation of hypergraph polytopes. Finally, we translate  both formalisations in the  language of simplicial complexes, which provides the environment for exhibiting the desired isomorphism.
\subsection{Non-inductive characterisation of constructs} \label{from-T-to-C-section}
Let $\hyper{H}$ be a finite, atomic and connected hypergraph.
We can define a map $\psi$ from the set of constructs of $\hyper{H}$ to ${\cal P}({\cal P}(H)\setminus\emptyset)\setminus\emptyset$, as follows (with  notation $\valup$ from Section \ref{construct-section}):
$$\psi(T)= \setc{\valup(Y)}{Y\;\mbox{is a (label of a) node of}\;T}\;.$$
We note that the Hasse diagram of $(\psi(T),\supseteq)$ is the same tree as $T$, replacing everywhere $Y$ by   $\valup(Y)$.    We also observe that the old decoration can be recovered from the new one by  noticing that
 $$Y=\valup(Y)\setminus \Union\setc{\valup(Z)}{Z\;\mbox{is a child of}\; Y \mbox{ in } T}\:.$$
From these observations, one can easily conclude that $\psi$ is injective. 

\begin{lemma}
The map $\psi$ is (contravariantly) monotonic and order-reflecting.
\end{lemma}
\Proof  Monotonicity is easy, following the inductive definition of $\leq$. For the second part of the statement,
 we  show that if  $\psi(T')\subseteq \psi(T)$, then  $T\leq_2 T'$, by induction on the size of $T$. 
 In what follows, for an arbitrary construct $T$, we will denote with $\rho(T)$ the root of $T$.\\
\indent If $T=H$, then $\psi(T)=\set{H}$ and  $\set{H}=\valup(\rho(T'))\inc \psi(T')\inc\psi(T)$  implies that also $\psi(T')=\set{H}$, and hence $T'=H$, and we conclude by clause 1 of the definition of $\leq_2$.\\
\indent If  $T=Y(S_1,\dots,S_m)$ ($m\geq 1$), let $T'=X(T_1,\dots,T_n)$ ($n\geq 0$). Since $\psi(T')\inc\psi(T)$, we get   $$\setc{\valup(\rho(T_i))}{1\leq i\leq n}\subseteq\setc{\valup(Z)}{Z\;\mbox{is a node of}\;T}.$$ Denote with $X_1,\dots,X_n$ the nodes of $T$ for which we have  $\valup(\rho(T_i))=\valup(X_i)$, and let, for each $1\leq i\leq n$, $U_i$ be the subtree of $T$ rooted at $X_i$. Note that all $X_i$'s must be different from $Y$. Indeed, if we had that $X_i=Y$ for some $1\leq i\leq n$, i.e., that  $\valup(\rho(T_i))=\valup(Y)=H$, this would imply that $X\subseteq \valup(\rho(T_i))$, which is not possible. We now have $$Y=H\backslash\bigcup_{j=1}^{m}\valup(\rho(S_j))\subseteq H\backslash\bigcup_{i=1}^{n}\valup(X_i)= H\backslash\bigcup_{i=1}^{n}\valup(\rho(T_i))=X.$$ Therefore, the first condition in the second clause defining  $\leq_2$ holds for $T$ and $T'$.  For the second condition, it is enough to establish (for all $j$)
$$\psi((K_j\cap X)\{T_i\,|\, H_i\inc K_j\})\subseteq\psi(S_j)\;,$$
which amounts to proving  $\psi(T_i)\inc \psi(S_j)$, for every $i$ such that $H_i\inc K_j$ . We have, on one hand, $\psi(T_i)\inc\psi(T')\inc\psi(T)$, and, on the other hand, for each element $Z$  of $\psi(T_i)$, $Z\inc H_i\inc K_j$, from which
$\psi(T_i)\inc \psi(S_j)$  follows, since  every non-root node $Z'$ of $S$, other than a node appearing in $S_j$, appears in some other $S_{j'}$, hence is included in $K_{j'}$, and not in $K_j$.    \qed

\medskip
We  now describe the  image of $\psi$.   We shall  characterise the constructs among all possible trees decorated with disjoint subsets of $H$, in a non-inductive way.  
 We note that the definition of $\valup$
makes sense for any such tree.

\smallskip
Recall that an antichain in a poset is a subset of pairwise uncomparable elements. We say that an antichain is proper if its cardinality is at least $2$.
\begin{lemma} \label{non-inductive-constructs-as-trees} Any of the following
properties
 characterises  constructs among  trees $T$ decorated with subsets of $H$:
\begin{enumerate}
\item At every non-leaf node of $T$, 
   $\valup(Y_1)$, \ldots,$\valup(Y_m)$  are the connected components of 
$\hyper{H}_{\valup(Y)\setminus Y}$, where $Y$ is the label of the node, and $Y_1,\ldots,Y_m$ are the labels of its child nodes.
\item The following  three conditions hold:
\begin{itemize} 
\item[A] All labels of the nodes of $T$ are pairwise disjoint and their union is $H$.
\item[B]  At each node $X$, $\valup(X)$ is such that $\hyper{H}_{\valup(X)}$ is connected.
\item[C'] At every non-leaf node $Y$ whose  child nodes are  $Y_1,\ldots,Y_m$, and any subset $I$ of $\set{1,\ldots,m}$  of cardinality at least $2$,  $\hyper{H}_{\Union\setc{\valup(Y_i)\,\,\,}{\,\, i \in I}}$ is not connected.
\end{itemize}
\item Conditions (A) and (B) hold, together with:
\begin{itemize}
\item[C] For each set  $\set{X_1,\ldots,X_m}$ of  labels of $T$ such that $\set{\valup(X_1),\ldots,\valup(X_m)}$
 is a proper antichain,  $\hyper{H}_{\valup(X_1) \union \ldots \union \valup(X_m)}$ is not connected.
\end{itemize}
\end{enumerate}
\end{lemma}
\Proof
(1) is  a  paraphrase of our inductive definition of construct.  We have that (3) obviously implies (2), since (C) a fortiori implies (C'). 

We now prove that (1) implies (3).  (A) and (B)
are obvious through the equivalence of (1) with our definition of inductively defined constructs.
We notice that (C) is vacuously true if $T$ is reduced to one node.    So let  $T=Y(T_1,\ldots,T_p)$, with
$\hyper{H},Y  \leadsto H_1,\ldots,H_p$. Let  $S=\set{X_1,\ldots,X_m}$ be as specified in  the statement, and suppose that $\hyper{H}_{\valup(X_1) \union \ldots \union \valup(X_m)}$ 
is connected.  Then it is included in one of the $H_i$'’s  (note that $(\valup(X_1) \union \ldots \union \valup(X_m)) \inter Y = \emptyset$). But then induction applies and we have a contradiction.

Finally, we prove that (2) implies (1).
By induction, it suffices to   check the property (1)  at the root of $T=X(T_1,\ldots,T_q)$. 
By  (B), we have that every $\valup(X_i)$  ($X_i$ root of $T_i$, $i\in \set{1,\ldots,q}$) is included in some $H_j$.  By (A), we have in fact that each $H_j$ is a union of some $\valup(X_i)$'s. Formally, there exists a non-empty set $I_j$ such that $H_j = \Union \setc{\valup(X_i)}{ i \in I_j}$.  But, by (C'),  $I_j$ must have cardinality $1$ (for every $j$).  Hence, up to permutation, we have $p=q$ and it follows that  (1) holds  at the root of $T$. \qed


%
\begin{proposition}  \label{equivalence-construct}
The map $\psi$ is a (contravariant) order-isomorphism between the set of constructs-as-decorated-trees and the collections of sets $M$ of connected (non-empty) subsets of $H$, containing $H$, and satisfying the following property:
\begin{itemize}
\item[C] For each proper antichain $S= \set{X_1,\ldots,X_m}\inc M$,  $\hyper{H}_{X_1 \union \ldots \union X_m}$ is not connected.
\end{itemize}
\end{proposition}
%
\Proof  By Lemma \ref{non-inductive-constructs-as-trees}, we have that, for any $T$, $\psi(T)$ satisfies (C)  (which we did not even care to rename!).  Conversely, we first show  that the  Hasse diagram of a set $M$ satisfying the conditions of the statement, ordered by reverse inclusion, is a tree.  We  note that if $X, Y$ are in $M$ and neither $X \inc  Y$ nor $Y \inc X$, and thus $\set{X,Y}$ is an antichain,  then, by (C), $\restrH{H}{X \union Y}$ is not connected.  This entails in particular that $X \inter Y$ is empty, as otherwise, since $\restrH{H}{X}$  and $\restrH{H}{Y}$ are connected  by assumption,  $\restrH{H}{X \union Y}$ would be connected. It follows that there cannot be a $Z$ above $X$ and $Y$ in the  Hasse diagram, as this would imply  $Z \inc  X \inter Y$, but  all elements of $M$, and $Z$ in particular, are non-empty: contradiction. Hence this Hasse diagram is a tree, with root $H$. Then, as remarked above,
it is easy to find  a  decoration of the same tree  where at each node the decoration $X$ is such that 
 $\valup(X)$ is  the corresponding element in the Hasse diagram.  Finally, we know from Lemma \ref{non-inductive-constructs-as-trees}  that this tree is indeed a construct $T$, and we have $\psi(T)=M$ by construction.
\qed
\begin{remark}  \label{constructs-nested-tubings}
\begin{enumerate}
\item Sets  as in Proposition \ref{equivalence-construct}
are called {\em nested sets} in  \cite{FS05,PRW08}.  Proposition  \ref{equivalence-construct} thus states that constructs as inductively defined trees are in order-isomorphic correspondence  with nested sets.  In their work, Do\v sen and Petri\'c adopt an intermediate viewpoint: they define constructions inductively, and  they  define constructs as subsets of  $\psi(V)$ containing $H$, for some construction $V$. They prove in \cite{DP-HP}[Proposition 6.13]  that their definition is equivalent to that of nested set.
\item
When $\hyper{H}$ is a graph (or has its set of connected subsets unchanged if restricted to hyperedges of cardinality $\leq 2$), the assumption (C) can be further relaxed to:
\begin{itemize}
\item[$\mbox{C}_g$] For each antichain $S=\set{X_1,X_2}\inc M$ (with respect to inclusion) of cardinality 2, we have that $\hyper{H}_{X_1 \union  X_2}$ is not connected.
\end{itemize}
Indeeed, if  (referring to (C)) $\hyper{H}_{X_1 \union \ldots \union X_m}$ were connected, then since connectedness is path-connectedness in a graph, we would have that  $\hyper{H}_{X_i \union X_j}$ is connected, for every pair of distinct $i,j\in\set{1,\ldots,m}$ (actually, picking just one such pair is enough for proving that ($\mbox{C}_g$) implies (C)). 

\smallskip
 Here is an example of why the stronger condition (C) is needed for general hypergraphs. Consider
$$\hyper{H}=\set{\set{x},\set{y},\set{z},\set{x,y,z}}\;.$$
Then  ($\mbox{C}_g$) holds, but $\set{\set{x},\set{y},\set{z}}$ is a witness that  (C) does not hold.
\item  Going back to graph polytopes, condition ($\mbox{C}_g$) is equivalent to the conditions (1) and (2) 
below:
\begin{itemize}
\item[1]  If $X_1,X_2\in M$  are such that $X_1\inter X_2\neq\emptyset$, then $X_1\subseteq X_2$ or $X_2\subseteq X_1$.
 \item[2] If $X_1,X_2\in M$  are such that $X_1\inter X_2=\emptyset$, then  $\hyper{H}_{X_1 \union  X_2}$ is not connected. 
\end{itemize}
That  (1) and (2) together imply  ($\mbox{C}_g$)  is obvious.  Conversely, we get (the contraposite of) (1) by arguing as in the proof above, and since the  implication $((X_1\subseteq X_2\;\mbox{or}\;X_2\subseteq X_1) \Rightarrow X_1\inter X_2\neq\emptyset)$ holds obviously, (1) is actually an equivalence, through which (2) can be rephrased as
($\mbox{C}_g$).  Conditions (1) and (2) are those given for {\em tubings} in \cite{CD-CCGA}.
 \item We summarise the  terminologies  used in the literature in the following table (see also  Section \ref{geometric-section}):
 \begin{center}
 \begin{tabular}{ccccccccccc}
 Combinatorial &&
 Hypergraphs   & Graphs &  Building sets\\
 && constructs & tubings & nested sets\\\\
 Geometrical && Hypergraph  & Graph  & Nestohedra\\
 && polytopes &  associahedra &
 \end{tabular}
 \end{center}

\end{enumerate}
\end{remark}

\vspace{-0.7cm}
\subsection{Geometric realisation} \label{geometric-section}

\vspace{-0.3cm}
 Following Do\v sen and Petri\'c, given a hypergraph $\hyper{H}$, we show how to associate 
 \begin{itemize}
 \item actual half-spaces and hyperplanes to the connected subsets of ${\hyper H}$ (i.e., to the hyperedges of ${\it Sat}(\hyper{H})$), 
 \item an actual polytope ${\cal G}(\hyper{H})$ to the whole hypergraph and  
 \item actual faces of ${\cal G}(\hyper{H})$  to constructs of $\hyper{H}$.
 \end{itemize}

Let  $H=\set{x_1,\ldots,x_n}$.  For every (non-empty)  $A\inc\set{1,\ldots,n}$, we define two subsets of $\mathbb{R}^n$, as follows: 
$$\begin{array}{lllll}
\pi_A^+=\setc{(x_1,\ldots,x_n)}{\sum_{i\in A} x_i\geq 3^{|A|}} &&&&
\pi_A =\setc{(x_1,\ldots,x_n)}{\sum_{i\in A} x_i = 3^{|A|}}\;.
\end{array}$$
 where $|A|$ is the cardinality of $A$.
Then the polytope associated with $\hyper{H}$ is defined as follows:
$${\cal G}(\hyper{H})= \Inter\setc{\pi^+_Y}{Y \in {\it Sat}(\hyper{H})\setminus\set{H}} \;\inter\; \pi_{H}\;.$$
For an arbitrary $M\inc{\it Sat}(\hyper{H})$, we define
$$\Pi(M) = \Inter\setc{\pi_Y}{Y\in M} \inter {\cal G}(\hyper{H})\;.$$
The definition of ${\cal G}(\hyper{H})$ implements the truncation instructions encoded by the hypergraph $\hyper{H}$.

\smallskip
This construction extends the realisation of the associahedra and of the cyclohedra  originally proposed in \cite{SS94,S-FOPIT}.
In   \cite{CD-CCGA}, graph-associahedra are also realised by means of truncations, although the concrete implementation of truncations is not described
(interestingly, Devadoss gives a more precise realisation in terms of convex hulls in  \cite{D-RGA}, that is also based on powers of $3$).

\smallskip
In the setting of building sets (cf. Remark \ref{constructs-nested-tubings}), 
a realisation that associates a linear inequality to every element of the building set, like in \cite{DP-HP}, can be found in \cite{Z06}. On the other hand, Feichtner et al.  use the elements of a building set as instructions for performing successive stellar subdivisions, starting from the simplex, while Postnikov et al. realise a building set by associating (via a fixed coordinate system) a simplex with each of its elements, and then by  taking the Minkowski sum of these simplices. They call the resulting polytopes nestohedra.


\subsection{Isomorphism between combinatorial and geometric faces}  \label{iso-section}
In this section, we  exhibit   an isomorphism between combinatorial and geometric faces of ${\cal G}(\hyper{H})$, which exploits the fact that ${\cal G}(\hyper{H})$ is a simple polytope.  

\smallskip
We first give an alternative definition of a geometric face. We defined a face  of a polytope  as the intersection of the polytope with a single hyperplane. But by allowing the intersection with several hyperplanes, the choice of those hyperplanes can be restricted, as stated in the following proposition, which is often taken as an alternative definition of geometric face. 
\begin{proposition} \label{faces-of-presented-polytopes}
Each  non-empty face of a polytope $P$ presented by a collection ${\cal S}$ of half-spaces  is defined as the intersection of $P$ with some of the hyperplanes bounding the half-spaces in ${\cal S}$.
\end{proposition}

We next introduce some notation. Given a polytope $P$, we let the letters $F,G$  (resp. $\Phi$) range  over the geometric faces (resp. the facets) of $P$.  We define a map $\phi$ from faces to sets of facets as follows:
$$\phi(F)=\setc{\Phi}{F\inc\Phi}\;.$$
We shall use the following equivalent characterisations of the notion of simple polytope  (which are the item (iii) and a sharpened version of the item  (v) of Proposition 2.16 of \cite{Ziegler}):
\begin{itemize}
\item [S1]  Each vertex of the polytope belongs to exactly $n$ facets of the polytope, where  $n$ is the dimension of the polytope.
\item [S2]
 For every face $F$, the restriction of $\phi$ to $\setc{G}{F\inc G}$ is an order-isomorphism onto 
${\cal P}(\phi(F))$.
\end{itemize}

 We shall also use the  following properties, which are  consequences of Lemmas 9.2, 9.4 and 9.5 of \cite{DP-HP}:
\begin{itemize}
\item [H1]   For every $M\inc \hyper{H}$, if $\Pi(M)$ is non-empty, then $M$ satisfies condition (C) of Proposition \ref{equivalence-construct}. 
\item [H2]
For every  construction $V$, $\Pi(\psi(V))$ is a vertex  $\set{v}$ of ${\cal G}(\hyper{H})$, and for every  $Y\in
 \hyper{H}\setminus\psi(V)$,  we have   $v \nin   \pi_Y$. Conversely, every vertex of ${\cal G}(\hyper{H})$  is obtained as  $\Pi(\psi(V))$ for some construction $V$.  
\end{itemize}
We take three steps in order to come up with the desired isomorphism.
%

\medskip
\noindent
{\bf A}) The poset of (non-empty) faces of a simple polytope is isomorphic to an abstract simplicial complex.

\smallskip
This is well-known, but since we want to express our isomorphisms explicitly, we briefly review here how this goes.
We start by some observations on polytopes (not necessarily simple). In any polytope, we have 
(cf. \cite{Ziegler}[Propositions 2.3 and 2.2]):
\begin{itemize}
\item  Every face of a polytope is the convex hull of its vertices.
\end{itemize}
We shall exploit two consequences of this property.
\begin{enumerate}
\item [P1] The map which associates with a face the set of all vertices that it contains  is monotonic and order-reflecting, and by polarity (cf. \cite{Ziegler}[Section 2.3]), it follows that the map $\phi$ defined above is (contravariantly) monotonic and order-reflecting.
\item [P2] Every non-empty face contains a vertex.
\end{enumerate}
Let $P$ be a polytope, with  vertices $\set{v_1},\ldots,\set{v_n}$. By P2, the lattice of faces (minus the empty face)  can be written as
 $\mathbb{L} = \mathbb{L}_1 \union \ldots \union  \mathbb{L}_n\;,$ 
where $\mathbb{L}_i$ is the set of faces containing $v_i$. 
If the polytope is simple, we know moreover by S2 that  
$\phi$ restricts to an order-isomorphism between 
$\mathbb{L}_i$ and  ${\cal P}(\phi(\set{v_i}))$. 
Let us now define an abstract simplicial complex $\mathbb{N}$ associated with $P$. Recall that a finite abstract simplicial complex (abreviated here as simplicial complex) is given by specifying a set $X$, called the support, and subsets $X_1,\ldots,X_n\inc X$, called the bases, that are  pairwise incomparable (w.r.t. inclusion), and are such that
$X=X_1\union\ldots \union X_n$.  The simplicial complex associated to these data is by definition the set
${\cal P}(X_1)\union\ldots\union{\cal P}(X_n)$, ordered by inclusion. The complex $\mathbb{N}$ is defined as follows:
\begin{itemize}
\item the support $X$ of $\mathbb{N}$ is the set of facets of $P$;
\item we take as bases the sets $\phi(\set{v})$, for all vertices of $P$.
\end{itemize}
Since the local isomorphisms between the $\mathbb{L}_i$'s and  ${\cal P}(\phi(\set{v_i}))$'s are  restrictions of the same function $\phi$,
it follows  that  
$\phi$ is an isomorphism from $\mathbb{L}$ to $\mathbb{N}$.

\medskip
\noindent
{\bf B}) An isomorphism of simplicial complexes.

\smallskip
First, we remark that  the partial order   ${\cal A}(\hyper{H})$  of constructs of $\hyper{H}$ can itself be organised as a simplicial complex, up to the isomorphism identifying each construct $T$ with $\psi(T)\setminus\set{H}$. Under these glasses, ${\cal A}(\hyper{H})$ is  isomorphic to the simplicial complex $\mathbb{M}$
\begin{itemize}
\item  whose support is $\hyper{H}\setminus\set{H}$,
\item and whose bases are the sets $\psi(V)\setminus\set{H}$, where $V$ ranges over the constructions of $\hyper{H}$. 
\end{itemize} This follows  from noting that any subset of a set satisfying condition (C) of Proposition
\ref{equivalence-construct} also satisfies that condition.  

\smallskip
Our goal is to define an isomorphism from $\mathbb{M}$ to $\mathbb{N}$. 
\begin{lemma}
A set $N$ of facets belongs to $\mathbb{N}$ if and only if $\Inter N$ is non-empty.
\end{lemma}
\Proof If $N\in\mathbb{N}$, then $N\inc \phi(\set{v})$, for some $v$, by definition of $\mathbb{N}$. It follows that
$\Inter\phi(\set{v})\inc \Inter N$. But  $v\in \Inter\phi(\set{v})$, by definition of $\phi$, hence $\Inter N$ is not empty.  Conversely, if $\Inter N$ is not empty, then, by (H1),  $\setc{X}{\Pi(\set{X})\in N}=\psi(T)$, for some construct $T$. By Corollary \ref{choice-of-construction-of-face}, we can choose a construction $V$, such that   $V\leq T$. Moreover, by (H2),  $\Pi(\psi(V))$ is a vertex $\set{v}$ of ${\cal G}(\hyper{H})$.  Let now $\chi(X)$ be a facet in $N$. Then $X\in\psi(V)$ since $\psi(T)\inc\psi(V)$, and therefore $\set{v}=\Pi(\psi(V))\inc \chi(X)$, and hence $N\inc \phi(\set{v})$. \qed

\begin{lemma} \label{hint-Zoran}
  If $T$ is a construct of $\hyper{H}$, and if $X\in \hyper{H}\setminus \psi(T)$, then there exists a construction $V \lessdot  T$, such that  $X\in \hyper{H}\setminus \psi(V)$.
\end{lemma}
\Proof  By induction on $T$. Let
$T=Y(T_1,\ldots, T_n)$ (possibly with $n=0$).   We distinguish two cases:
\begin{enumerate}
\item  $X\inter Y=\emptyset$.  Then, $n\geq 1$, and for each $1\leq i\leq n$,  since  $\psi(T_i)\inc \psi(T)$ for each $i$, we can apply induction to  $T_i$'s, and get $V_i\lessdot  T_i$ satisfying the statement relative to $T_i$.  Let then $S'$ be an arbitrary partial construction spanning $Y$.  By grafting the $V_i$'s on the corresponding occurrences of $\Omega$ of $S'$, we get a construction $V\lessdot  T$, which satisfies the statement: this is clear for all nodes $x$ coming from the $V_i$'s, while  all nodes coming from $S'$, being  elements of $Y$, are such that $\valup(x)\inter Y\neq\emptyset$, which implies $\valup(x)\neq X$.
\item  $X\inter Y\neq\emptyset$.  Let $y\in X\inter Y$. By Corollary \ref{choice-of-construction-of-face}, we can choose a construction $V\lessdot  T$ whose root is decorated by $y$.  Then $\psi(V)\setminus\set{H} $ consists only of sets that do not contain $y$, hence none of them can be $X$.  \qed
\end{enumerate}
\begin{lemma} \label{Phi-H}
 The  elements of $\hyper{H}\setminus \{H\}$ are in one-to-one correspondence with the facets of ${\cal G}(\hyper{H})$, through the map $\chi$ defined by  $\chi(X) = \Pi(\set{X})$.
\end{lemma}
\Proof 
We need to show that $\chi$ is both bijective and well-defined, in the sense that $\Pi(\set{X})$ is actually a facet.  We take the following steps.
\begin{enumerate}
\item
For all $X, Y\in\hyper{H}\setminus\set{H}$, if  $X\neq Y$, then  $\chi(X)$ is not included in $\chi(Y)$ (this a fortiori implies that $\chi$ is injective).  Since $X\neq Y$, we have
$Y\nin\psi((H\setminus X)(X))$. Then, by Lemma \ref{hint-Zoran}, there exists a construction $V$ such that  $V\leq (H\setminus X)(X)$ and  $Y\nin\psi(V)$. By (H2), we have $\Pi(\psi(V))=\set{v}$ for some $v$ such that $v\nin\pi_Y$, and therefore $v\nin \chi(Y)$. On the other hand, $V\leq (H\setminus X)(X)$ implies  $v\in \chi(X)$, which proves the claim.
\item  $\chi(X)$ is a facet, for all $X$.
Suppose that $\chi(X) \incs  F$ for some face $F$ of ${\cal G}(\hyper{H})$.  It follows from Proposition \ref{faces-of-presented-polytopes} that every face is included in some $\chi(Y)$  (just pick one of  the hyperplanes in the statement).  So we have $F\inc \chi(Y)$ for some $Y$, and a fortiori $\chi(X)\inc \chi(Y)$, from which we deduce $X=Y$ by (1).  But this forces $\chi(X)=F$, contradicting our assumption.
\item
$\chi$ is surjective.  We already observed that every face is included in some $\chi(Y)$, from which surjectivity follows. \qed
\end{enumerate}

Then the claimed isomorphism  from $\mathbb{M}$ to $\mathbb{N}$ is defined through the map $\chi$ of Lemma \ref{Phi-H}, using the following easy fact.
\begin{itemize}
\item If  
$\chi$ is a bijection from the support  of $\mathbb{M}$ to the support of $\mathbb{N}$ whose extension  to subsets (notation $\chi[M]=\setc{\chi(X)}{X\in M}$) is such that, for all subsets $M,N$ of the respective supports, we have 
$M\in \mathbb{M}\Leftrightarrow \chi[M]\in\mathbb{N}$,
then it defines an order-isomorphism between $\mathbb{M}$ and $\mathbb{N}$.
\end{itemize}

\noindent
{\bf C}) ${\cal G}(\hyper{H})$ is simple.

\smallskip
First, we establish the dimension of {${\cal G}(\hyper{H})$.
\begin{lemma} \label{dimension-of-HP}
If $H$ has cardinality $n+1$, then ${\cal G}(\hyper{H})$ has dimension  $n$.
\end{lemma}
\Proof It is enough to prove the statement in the case of the permutohedron, since ${\cal G}(\hyper{H})$ contains the  permutohedron defined by the complete graph on $H$. Simple calculations prove that  the point
$(3^{n+1}/n+1,\ldots , 3^{n+1}/n+1)$ lies in $\pi_H$ and in the interior of $\pi_Y^+$ for all  non-empty  $Y\incs H$, from which one concludes easily. \qed

\smallskip
We  prove simplicity via condition S1, as follows. First, the dimension of  $ {\cal G}(\hyper{H})$ is  $|H|-1$, by Lemma \ref{dimension-of-HP}.
Second, we note that a construction $V$ has always exactly $|H|$ nodes, hence $\psi(V)$ has   exactly $|H|-1$ elements different from $H$. Since, by (H2), every vertex can be written as  $\set{v}=\Pi(\psi(V))$, for some construction $V$,  we conclude by observing that $\Pi(\psi(V))$ is included by definition  in  all of the $|H|-1$ facets $\chi(X)$, for $X$ ranging over $\psi(V)\setminus\set{H}$, and in no other facet, by  (H2) and Lemma \ref{Phi-H}.
 
\smallskip
Thus we can combine steps (B) and (A).
\begin{theorem} \label{iso-theorem}
The map $\Pi\comp\psi$, where $\psi$ and $\Pi$ are defined in Sections \ref{from-T-to-C-section} and  \ref{geometric-section},  is an order-isomorphism.
\end{theorem}
\Proof  Our analysis gives us  the isomorphism $\phi^{-1}\comp\chi$. It can be shown that ``taking the intersection'' is inverse to $\phi$, which  allows us to reformulate the  isomorphism as follows 
$\phi^{-1}(\chi[\psi(T)\setminus\set{H}])=\Pi(\psi(T))$. \qed

\section{Operadic coherences}  \label{operadic-application-section}
In  \cite{DP-WCO}, Do\v sen and Petri\'c have used hypergraph polytopes in the study   of coherences arising when categorifying the notion of operad \cite{LV-AO}, i.e., when the axioms of sequential and parallel associativity are turned into coherent isomorphisms $\beta$ and $\theta$, the coherence conditions being   naturally associated with suitable  polytopes.
We shall not need the precise definition of an operad, and shall rely instead on simple graphical intuitions.
\subsection{Weak Cat-operads} \label{weak-cat-operads-section}
%

In monoidal categories,  a coherence condition is  imposed on the associator $\alpha_{A,B,C}:(A\otimes B)\otimes C\rightarrow A\otimes(B\otimes C)$, ensuring that all the diagrams made of instances of $\alpha$ (possibly whiskered by identites), and their inverses,  commute. This condition is expressed by the commutation of Mac Lane's  pentagon (see diagram (2) on the next page).
In an operad, the role of the objects $A,B$ of a monoidal category is now played by operations labelling the nodes of a rooted (non-planar) tree.  We call such a tree, acting as a pasting scheme, an {\em operadic tree}. Any two neigbouring operations in an operadic tree may be composed (imagine that the edge connecting them is  contracted in the process), and then composed with a neigbouring operation, etc.  The axioms of operads guarantee that the overall composition of the operations in the tree does not depend on the order of compositions.  Consider the two trees with three nodes:
\begin{center}
  \resizebox{0.4cm}{!}{\begin{tikzpicture}[scale=0.7]
    \node (E)[circle,draw=black,minimum size=4mm,inner sep=0.1mm] at (0,0) {\small $a$};
    \node (F) [circle,draw=black,minimum size=4mm,inner sep=0.1mm] at (0,1) {\small $b$};
    \node (A) [circle,draw=black,minimum size=4mm,inner sep=0.1mm] at (0,2) {\small $c$};
    \draw[-] (E)--(F) node  {};
    \draw[-] (F)--(A) node  {};
   \end{tikzpicture}}
   \enspace\enspace\enspace\enspace\enspace\enspace\enspace  \begin{tikzpicture}
    \node (E)[circle,draw=none,minimum size=4mm,inner sep=0.1mm] at (0,1) {and};
    \node (F) [circle,draw=none,minimum size=4mm,inner sep=0.1mm] at (0,0) {};
   \end{tikzpicture}  \enspace\enspace\enspace\enspace\enspace\enspace\enspace
\resizebox{1.43cm}{!}{\begin{tikzpicture}
    \node (E)[circle,draw=black,minimum size=4mm,inner sep=0.1mm] at (0,0) {\small $a$};
    \node (F) [circle,draw=black,minimum size=4mm,inner sep=0.1mm] at (-0.5,1) {\small $b$};
    \node (A) [circle,draw=black,minimum size=4mm,inner sep=0.1mm] at (0.5,1) {\small $c$};
    \draw[-] (E)--(F) node  {};
    \draw[-] (E)--(A) node  {};
    \end{tikzpicture}}
\end{center}
The axiom of sequential (resp. parallel) associativity says that the two ways to build the tree on the left (resp. on the right) by means of grafting
and to perform compositions accordingly, yield the same operation:  first compose $a$ with $b$, and then compose with (or {\em insert}) $c$, or first compose $b$ with $c$ and then insert $a$  (resp. first compose $a$ with $b$, or first compose $a$ with $c$).  In a weak Cat-operad, these identifications are turned into isomorphisms
$$\beta:(ab)c \rightarrow a(bc)\quad \mbox{and}\quad\theta:(ab)c \rightarrow (ac)b$$
(writing composition as juxtaposition).  
 
 \smallskip
To synthesise the coherence conditions that $\beta$ and $\theta$ have to satisfy, we need to consider the four possible shapes of operadic trees with four nodes:
\begin{center}
\resizebox{9.75cm}{!}{\begin{tikzpicture}[scale=0.7]
\node (E1)[circle,draw=none,minimum size=4mm,inner sep=0.1mm] at (0,-0.75) { (1)};
    \node (E)[circle,draw=black,minimum size=4mm,inner sep=0.1mm] at (0,0) {\scriptsize  $a$};
    \node (F) [circle,draw=black,minimum size=4mm,inner sep=0.1mm] at (-1,1) {\scriptsize  $b$};
    \node (A) [circle,draw=black,minimum size=4mm,inner sep=0.1mm] at (0,1) {\scriptsize  $c$};
    \node (Asubt) [circle,draw=black,minimum size=4mm,inner sep=0.1mm] at (1,1) {\scriptsize  $d$};
    \draw[-] (E)--(F) node  {};
    \draw[-] (E)--(A) node  {};
 \draw[-] (E)--(Asubt) node {};
   \end{tikzpicture}
   \enspace\enspace\enspace\enspace\enspace\enspace\enspace
   \begin{tikzpicture}
\node (E1)[circle,draw=none,minimum size=4mm,inner sep=0.1mm] at (0,-0.75) { (2)};
    \node (E)[circle,draw=black,minimum size=4mm,inner sep=0.1mm] at (0,0) {\scriptsize  $a$};
    \node (F) [circle,draw=black,minimum size=4mm,inner sep=0.1mm] at (0,1) {\scriptsize $b$};
    \node (A) [circle,draw=black,minimum size=4mm,inner sep=0.1mm] at (0,2) {\scriptsize   $c$};
    \node (Asubt) [circle,draw=black,minimum size=4mm,inner sep=0.1mm] at (0,3) {\scriptsize $d$};
    \draw[-] (E)--(F) node  {};
    \draw[-] (F)--(A) node  {};
 \draw[-] (A)--(Asubt) node {};
   \end{tikzpicture}
      \enspace\enspace\enspace\enspace\enspace\enspace\enspace
    \begin{tikzpicture}
\node (E1)[circle,draw=none,minimum size=4mm,inner sep=0.1mm] at (0,-0.75) { (3)};
    \node (E)[circle,draw=black,minimum size=4mm,inner sep=0.1mm] at (0,0) {\scriptsize  $a$};
    \node (F) [circle,draw=black,minimum size=4mm,inner sep=0.1mm] at (0,1) {\scriptsize  $b$};
    \node (A) [circle,draw=black,minimum size=4mm,inner sep=0.1mm] at (-0.5,2) {\scriptsize $c$};
    \node (Asubt) [circle,draw=black,minimum size=4mm,inner sep=0.1mm] at (0.5,2) {\scriptsize $d$};
    \draw[-] (E)--(F) node  {};
    \draw[-] (F)--(A) node  {};
 \draw[-] (F)--(Asubt) node {};
   \end{tikzpicture}
 \enspace\enspace\enspace\enspace\enspace\enspace\enspace
   \begin{tikzpicture}
\node (E1)[circle,draw=none,minimum size=4mm,inner sep=0.1mm] at (0,-0.75) { (4)};
    \node (E)[circle,draw=black,minimum size=4mm,inner sep=0.1mm] at (0,0) {\scriptsize $a$};
    \node (F) [circle,draw=black,minimum size=4mm,inner sep=0.1mm] at (-0.5,1) {\scriptsize $b$};
    \node (A) [circle,draw=black,minimum size=4mm,inner sep=0.1mm] at (0.5,1) {\scriptsize $c$};
    \node (Asubt) [circle,draw=black,minimum size=4mm,inner sep=0.1mm] at (-0.5,2) {\scriptsize $d$};
    \draw[-] (E)--(F) node  {};
    \draw[-] (E)--(A) node  {};
 \draw[-] (F)--(Asubt) node {};
   \end{tikzpicture}
}
\end{center}
 Each of these trees guides the interpretation of
 parenthesised  words such as  $((ab)c)d$
  as   sequences of insertions,
and each of the diagrams below (one for each tree) features the resolution of the critical pair (or overlapping)
$$(\underbrace{(ab)c})d  \quad\quad \underbrace{((ab)c)d}\;,$$
interpreted diversely according to whether the associativities are parallel or sequential, as prescribed by the respective trees.
\begin{center}
\resizebox{10.5cm}{!}{\begin{tabular}{ccc}
\begin{tikzpicture}
\node (E1)[circle,draw=none,minimum size=4mm,inner sep=0.1mm] at (0,-1.7) {\small (1)};
    \node (E) at (-1.4,0) {$((ab)c)d$};
    \node (G) at (1.4,0) {$((ac)b)d$};
    \node (F) at (2.2,-1.7) {$((ac)d)b$};
    \node (A) at (-2.2,-1.7) {$((ab)d)c$};
    \node (Asubt) at (-1.4,-3.5) {$((ad)b)c$};
    \node (P4) at (1.4,-3.5) {$((ad)c)b$};
    \draw[->] (E)--(G) node [midway,above] {\footnotesize   $\theta$};
    \draw[->] (G)--(F) node [midway,right] {\footnotesize   $\theta$};
    \draw[->] (E)--(A) node [midway,left] {\footnotesize  $\theta$};
 \draw[->] (F)--(P4) node [midway,right] {\footnotesize  $\theta$};
    \draw[->] (A)--(Asubt) node [midway,left]  {\footnotesize   $\theta$};
    \draw[->] (Asubt)--(P4) node [midway,above] {\footnotesize  $\theta$};
   \end{tikzpicture}

&&

\begin{tikzpicture}
\node (E1)[circle,draw=none,minimum size=4mm,inner sep=0.1mm] at (0,-1.7) {\small (2)};
    \node (E) at (0,0) {$((ab)c)d$};
    \node (F) at (-2,-1.7) {$(a(bc))d$};
    \node (A) at (2,-1.7) {$(ab)(cd)$};
    \node (Asubt) at (-1.5,-3.5) {$a((bc)d)$};
    \node (P4) at (1.5,-3.5) {$a(b(cd))$};
    \draw[->] (E)--(F) node [midway,above] {\footnotesize   $\beta$};
    \draw[->] (E)--(A) node [midway,above] {\footnotesize  $\beta$};
 \draw[->] (F)--(Asubt) node [midway,left] {\footnotesize  $\beta$};
    \draw[->] (A)--(P4) node [midway,right]  {\footnotesize   $\beta$};
    \draw[->] (Asubt)--(P4) node [midway,above] {\footnotesize  $\beta$};
   \end{tikzpicture}

\end{tabular}}
\end{center}

\begin{center}
\resizebox{10.5cm}{!}{\begin{tabular}{ccc}
\begin{tikzpicture}
\node (E1)[circle,draw=none,minimum size=4mm,inner sep=0.1mm] at (0,-1.7) {\small (3)};
    \node (E) at (-1.4,0) {$((ab)c)d$};
    \node (G) at (1.4,0) {$((ab)d)c$};
    \node (F) at (2.2,-1.7) {$(a(bd))c$};
    \node (A) at (-2.2,-1.7) {$(a(bc))d$};
    \node (Asubt) at (-1.4,-3.5) {$a((bc)d)$};
    \node (P4) at (1.4,-3.5) {$a((bd)c)$};
    \draw[->] (E)--(G) node [midway,above] {\footnotesize   $\theta$};
    \draw[->] (G)--(F) node [midway,right] {\footnotesize   $\beta$};
    \draw[->] (E)--(A) node [midway,left] {\footnotesize  $\beta$};
 \draw[->] (F)--(P4) node [midway,right] {\footnotesize  $\beta$};
    \draw[->] (A)--(Asubt) node [midway,left]  {\footnotesize   $\beta$};
    \draw[->] (Asubt)--(P4) node [midway,above] {\footnotesize  $\theta$};
   \end{tikzpicture}

&&

\begin{tikzpicture}
\node (E1)[circle,draw=none,minimum size=4mm,inner sep=0.1mm] at (0,-1.7) {\small (4)};
    \node (E) at (0,0) {$((ab)c)d$};
    \node (F) at (-2,-1.7) {$((ac)b)d$};
    \node (A) at (2,-1.7) {$((ab)d)c$};
    \node (Asubt) at (-1.5,-3.5) {$(ac)(bd)$};
    \node (P4) at (1.5,-3.5) {$(a(bd))c$};
    \draw[->] (E)--(F) node [midway,above] {\footnotesize   $\theta$};
    \draw[->] (E)--(A) node [midway,above] {\footnotesize  $\theta$};
 \draw[->] (F)--(Asubt) node [midway,left] {\footnotesize  $\beta$};
    \draw[->] (A)--(P4) node [midway,right]  {\footnotesize   $\beta$};
    \draw[->] (Asubt)--(P4) node [midway,above] {\footnotesize  $\theta$};
   \end{tikzpicture}

\end{tabular}}
\end{center}
\begin{remark}
 Before asking the question of distinguishing $\beta$ and $\theta$ edges in these ``operadic polytopes" (and, in general, in operadic polytopes of arbitrary dimension), one must be able to systematically assign them labels.
In Section  \ref{construct-examples-section}, we have seen various ways to label all the faces of the pentagon and of the hexagon, that would work here for the pentagon made of $\beta$-arrows only and the hexagon made of $\theta$-arrows only.
However, it is a priori not  clear how we could do this for the other two mixed $\beta/\theta$-diagrams.  This question will be addressed in Section \ref{graph-of-tree-section}.
\end{remark}

By ``lifting'' the methodology of coherence chasing to the 3-dimensional setting, i.e., by considering trees with 5 nodes, we find  $9$ possible configurations. We shall draw only three of them:
\begin{center}
\resizebox{9.5cm}{!}{\begin{tabular}{ccccccccccccccc}
\begin{tikzpicture}[scale=0.7]
    \node (E)[circle,draw=black,minimum size=4mm,inner sep=0.1mm] at (0,0) {\scriptsize $a$};
    \node (F) [circle,draw=black,minimum size=4mm,inner sep=0.1mm] at (-0.35,1) {\scriptsize $c$};
    \node (A) [circle,draw=black,minimum size=4mm,inner sep=0.1mm] at (-1,1) { \scriptsize $b$};
    \node (Asubt) [circle,draw=black,minimum size=4mm,inner sep=0.1mm] at (0.35,1) {\scriptsize $d$};
    \node (P) [circle,draw=black,minimum size=4mm,inner sep=0.1mm] at (1,1) { \scriptsize $e$};
    \draw[-] (E)--(F) node  {};
    \draw[-] (E)--(A) node  {};
 \draw[-] (E)--(Asubt) node {};
 \draw[-] (E)--(P) node {};
   \end{tikzpicture}
&& &&&

   \begin{tikzpicture}
    \node (E)[circle,draw=black,minimum size=4mm,inner sep=0.1mm] at (0,0) {\scriptsize $a$};
    \node (F) [circle,draw=black,minimum size=4mm,inner sep=0.1mm] at (0,0.75) {\scriptsize $b$};
    \node (A) [circle,draw=black,minimum size=4mm,inner sep=0.1mm] at (0,1.5) {\scriptsize $c$};
    \node (Asubt) [circle,draw=black,minimum size=4mm,inner sep=0.1mm] at (0,2.25) {\scriptsize $d$};
    \node (P) [circle,draw=black,minimum size=4mm,inner sep=0.1mm] at (0,3) {\scriptsize $e$};
    \draw[-] (E)--(F) node  {};
    \draw[-] (F)--(A) node  {};
 \draw[-] (A)--(Asubt) node {};
 \draw[-] (Asubt)--(P) node {};
   \end{tikzpicture}
&&&&&
\begin{tikzpicture}
    \node (E)[circle,draw=black,minimum size=4mm,inner sep=0.1mm] at (0,0) {\scriptsize $a$};
    \node (F) [circle,draw=black,minimum size=4mm,inner sep=0.1mm] at (-0.5,1) {\scriptsize $b$};
    \node (A) [circle,draw=black,minimum size=4mm,inner sep=0.1mm] at (0.5,1) {\scriptsize $e$};
    \node (Asubt) [circle,draw=black,minimum size=4mm,inner sep=0.1mm] at (-1,2) {\scriptsize $c$};
    \node (P) [circle,draw=black,minimum size=4mm,inner sep=0.1mm] at (0,2) {\scriptsize $d$};
    \draw[-] (E)--(F) node  {};
    \draw[-] (E)--(A) node  {};
 \draw[-] (F)--(Asubt) node {};
 \draw[-] (F)--(P) node {};
   \end{tikzpicture}
\end{tabular}}
\end{center}
The expression $(((ab)c)d)e$ is now subject to a three-fold overlapping,
$$\begin{array} {ccc}
((\underbrace{(ab)c})d)e & (\underbrace{((ab)c)d})e & \underbrace{(((ab)c)d)e}
\end{array}$$
which is resolved  differently for  each of the 9 trees, leading to 9 ``coherence conditions between coherences''  (in a framework  where   the coherence equations would not hold up to equality),
each described by a suitable 3-dimensional polytope.

 For the first  two  trees above,
 we get  the $3$-dimensional permutohedron and  associahedron, respectively, whose edges all stand for $\theta$-arrows in the first case, and $\beta$-arrows in the second. For the third  one,
we get a polytope called the {\em hemiassociahedron}, which, as we shall see,  also belongs to the familly of hypergraph polytopes. 
In Figure \ref{hemi-assoc}, we labelled some of the vertices of this polytope, matching them with decompositions of our example tree (this matching will be spelled out  in Proposition \ref{constructions-insertions}).
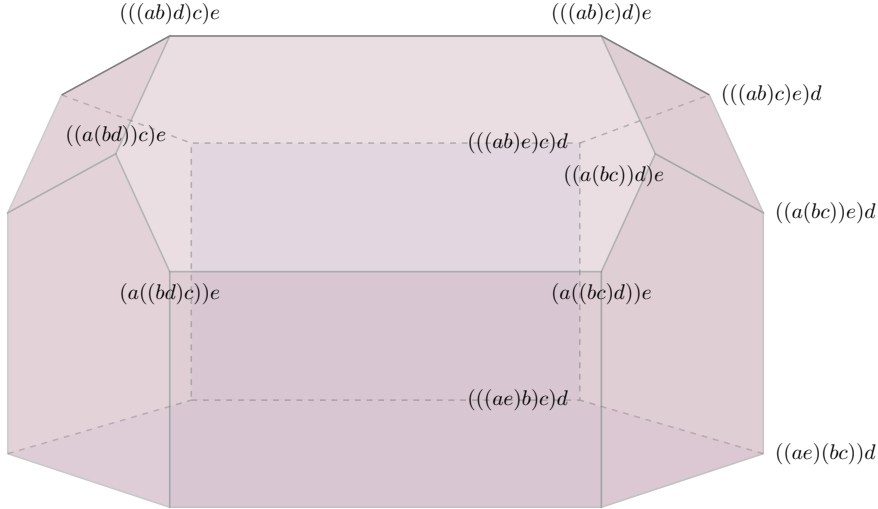
\begin{figure}[H] 
\resizebox{11.5cm}{!}{\begin{tikzpicture}[thick,scale=19.5]
\coordinate (A1) at (-0.2,0);
\coordinate (A2) at (0.2,0); 
\coordinate (D1) at (-0.35,0.05);
\coordinate (D2) at (0.35,0.05);
\coordinate (D3) at (-0.18,0.1);
\coordinate (D4)  at (0.18,0.1);
\coordinate (E3) at (-0.18,0.34);
\coordinate (E4) at (0.18,0.34);
\coordinate (A3) at (-0.2,0.22);
\coordinate (A4) at (0.2,0.22);
 \coordinate (C1) at (-0.25,0.33);
\coordinate (C2) at (0.25,0.33);
 \coordinate (F1) at (-0.35,0.275);
\coordinate (F2) at (0.35,0.275);
 \coordinate (G1) at (-0.3,0.385);
\coordinate (G2) at (0.3,0.385);
 \coordinate (B3) at (-0.2,0.44);
\node at (B3) [above = 1mm of B3] {\large $(((ab)d)c)e$};
\coordinate (B4) at (0.2,0.44);
\node at (B4) [above = 1mm of B4] {\large $(((ab)c)d)e$};
\draw[draw=gray,fill=pur,opacity=0.45] (A1) -- (A2) -- (A4) -- (A3) -- cycle;
\draw[draw=gray,fill=pur,opacity=0.3] (A3) -- (A4) -- (C2) -- (B4) -- (B3) -- (C1)-- cycle;
\draw[draw=gray,fill=greeo,opacity=0.05] (A2) -- (A1) -- (D1) -- (D3) -- (D4) -- (D2) --  cycle;
\draw[draw=gray,fill=greeo,opacity=0.05,dashed] (D3) -- (D4) -- (E4) -- (E3) -- cycle;
\draw[draw=gray,dashed,opacity=0.5] (D3) -- (D4) -- (E4) -- (E3) -- cycle;
\draw[draw=gray,dashed,opacity=0.5] (D1)   -- (D3);
\draw[draw=gray,dashed,opacity=0.5] (D2)   -- (D4);
\draw[draw=gray] (G1) -- (B3) -- (B4) -- (G2);
\draw[draw=gray,dashed,opacity=0.5]   (G1) -- (E3); 
\draw[draw=gray,dashed,opacity=0.5] (E4) -- (G2);
\draw[draw=gray,fill=pur,opacity=0.4]  (A1) -- (D1) -- (F1) -- (C1) -- (A3) --  cycle;
\draw[draw=gray,fill=pur,opacity=0.45]  (A2) -- (D2) -- (F2) -- (C2) -- (A4) --  cycle;
\draw[draw=gray,fill=pur,opacity=0.4]  (F1) -- (G1) -- (B3) -- (C1)  --  cycle;
\draw[draw=gray,fill=pur,opacity=0.45] (F2) -- (G2) -- (B4) -- (C2)  --  cycle;
 \node at (C2) [left= 0.3mm of C2,yshift=-4mm,xshift=0.3cm] {\large $((a(bc))d)e$};
 \node at (C1) [above= 0.3mm  of C1] {\large $((a(bd))c)e$};
 \node at (A4) [below = 1mm of A4] {\large $(a((bc)d))e$};
 \node at (A3) [below= 1mm of A3] {\large  $(a((bd)c))e$};
 \node at (G2) [right= 1mm of G2] {\large  $(((ab)c)e)d$};
  \node at (F2) [right= 1mm of F2] {\large  $((a(bc))e)d$};
  \node at (D2) [right= 1mm of D2] {\large $((ae)(bc))d$};
  \node at (D4) [left= 1mm of D4] {\large  $(((ae)b)c)d$};
\node at (E4) [left= 1mm of E4] {\large $(((ab)e)c)d$};





\end{tikzpicture}}
\caption{The hemiassociahedron}
\label{hemi-assoc}
\end{figure}

\subsection{Graphs associated with operadic trees} \label{graph-of-tree-section}

To every rooted tree ${\cal T}$ representing a pasting scheme for operadic operations, Do\v sen and Petri\'c associate a graph $\mathbb{G}({\cal T})$, obtained as follows.  Its vertices are the edges of ${\cal T}\!$, and two vertices are connected whenever as edges of ${\cal T}$ they share a common vertex.  

 It is clear that one can identify the edges of ${\cal T}$ with the non-root nodes of ${\cal T}$  (for example, in Figure \ref{G(T)-picture},  there is a bijection mapping $x$ to $c$, $y$ to $d$, $z$ to $b$, and $u$ to $e$). By this identification,  seeing now the nodes of $\mathbb{G}({\cal T})$ as the non-root nodes of ${\cal T}$, all edges of ${\cal T}$, apart from those stemming from the root, are in $\mathbb{G}({\cal T})$. All the other edges of $\mathbb{G}({\cal T})$ are edges witnessing  that  two edges of ${\cal T}$ are siblings. We record the latter (resp. the former) by representing them with  a dashed  (resp.  solid) line. 

The graph   $\mathbb{G}({\cal T})$ is  connected and can be represented itself as a  tree with some horizontal dashed edges such that, by construction, each dashed horizontal zone is a complete graph all of whose nodes are connected to their father node (if it exists)  by solid edges.  The nodes of $\mathbb{G}({\cal T})$ are thus  organised in levels.  We say that $\mathbb{G}({\cal T})$ has a root when there is no horizontal dashed layer at the bottom of ${\mathbb G}({\cal T})$.

 Figure \ref{G(T)-picture} shows the  graph associated to the third tree considered at the end of Section \ref{weak-cat-operads-section} ($z,u$ are at level $1$, and $x,y$ are at level $2$).

\begin{figure}[t]  
\begin{tabular}{ccc}
\resizebox{2cm}{!}{\begin{tikzpicture}[scale=0.8]
    \node (E)[circle,draw=black,minimum size=4mm,inner sep=0.1mm] at (0,0) {\scriptsize $a$};
    \node (F) [circle,draw=black,minimum size=4mm,inner sep=0.1mm] at (-0.5,1) { \scriptsize $b$};
    \node (A) [circle,draw=black,minimum size=4mm,inner sep=0.1mm] at (0.5,1) {\scriptsize $e$};
    \node (Asubt) [circle,draw=black,minimum size=4mm,inner sep=0.1mm] at (-1,2) {\scriptsize  $c$};
    \node (P) [circle,draw=black,minimum size=4mm,inner sep=0.1mm] at (0,2) {\scriptsize $d$};
    \draw[-] (E)--(F) node  [midway,left] {\scriptsize $z$};
    \draw[-] (E)--(A) node  [midway,right] {\scriptsize $u$};
 \draw[-] (F)--(Asubt) node [midway,left] {\scriptsize $x$};
 \draw[-] (F)--(P) node [midway,right] {\scriptsize $y$};
   \end{tikzpicture}}

&&
\resizebox{2cm}{!}{
\begin{tikzpicture}
    \node (Z)[] at (-0.5,0) {$z$};
    \node (U)[]  at (0.5,0) {$u$};
    \node (X)[]  at (-1,1) {$x$};
    \node (Y)[]  at (0,1) {$y$};
    \draw[dashed] (Z)--(U) node  {};
 \draw[-] (Z)--(X) node  {};
 \draw[-] (Z)--(Y) node {};
 \draw[dashed] (X)--(Y) node {};
   \end{tikzpicture}}
\end{tabular}
\caption{The ${\mathbb G}({\cal T})$ construction} \label{G(T)-picture}
\end{figure}
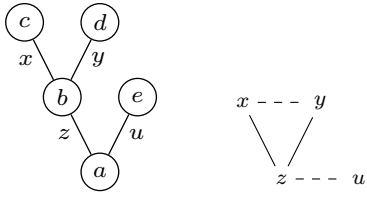
We insist that the   dashed/solid informations on the edges of $\mathbb{G}({\cal T})$  are
not part of the graph structure $\mathbb{G}({\cal T})$: they are additional data that we shall use
to derive both the type ($\beta$ or $\theta$) and (in the case of $\beta$) the orientation of all edges of the corresponding polytope (as dictated by ${\cal T}$). 

Recall that, in the language of constructs, vertices   are trees whose nodes are all labelled with singletons. An  edge $E$  is a  tree whose   nodes are all singletons, {\em except one}, which is a two-element set $\set{u_E,v_E}$. 
We will show that
$\mathbb{G}({\cal T})$,  together with its bipartition of dashed and solid edges, determines the type (and the orientation) of $E$. Let us call a min-path of a graph a path of minimum length  between  two vertices (we will show that in $\mathbb{G}({\cal T})$ min-paths are always unique).  Our criterion is the following:
\begin{itemize}
\item[$\dagger$]  If the min-path between $u_E$ and $v_E$ in $\mathbb{G}({\cal T})$ is made only of solid edges,  $E$ corresponds to a $\beta$-arrow, oriented towards the vertex of $E$ 
in which the label  $u_E$ appears below the label $v_E$  if and only if the level of  $u_E$  is inferior to the level of $v_E $ in $\mathbb{G}({\cal T})$. Otherwise, $E$ witnesses a $\theta$-arrow.
\end{itemize}

\smallskip
As an example, let us derive the edge information for the mixed pentagon (4), out of the associated graph:
\begin{center}
\begin{tabular}{cccc}
\resizebox{1.75cm}{!}{\begin{tikzpicture}
    \node (Z)[] at (0,1.5) {$z$};
    \node (U)[]  at (1.5,0) {$y$};
    \node (Y)[]  at (0,0) {$x$};
    \draw[dashed] (Y)--(U) node  {};
 \draw[-] (Z)--(Y) node  {};
   \end{tikzpicture}}
&&&
\resizebox{5.5cm}{!}{\begin{tikzpicture}
    \node (E) at (0,0) {$z(y(x))$};
    \node (F) at (-2,-1.7) {$z(x(y))$};
    \node (A) at (2,-1.7) {$y(z(x))$};
    \node (Asubt) at (-1.5,-3.5) {$x(y,z)$};
    \node (P4) at (1.5,-3.5) {$y(x(z))$};
\node (M) at (0,-1.7) {$\set{x,y,z}$};
    \draw[-] (E)--(F) node [xshift=-0.55cm,midway,above] {\footnotesize   $z(\set{x,y})$};
    \draw[-] (E)--(F) node [xshift=0.1cm,yshift=0.1cm,midway,below] {\footnotesize   $\theta$};
    \draw[-] (E)--(A) node [xshift=0.55cm,midway,above] {\footnotesize  $\set{y,z}(x)$};
    \draw[-] (E)--(A) node [xshift=-0.15cm,yshift=0.1cm,midway,below] {\footnotesize  $\theta$};
 \draw[-] (F)--(Asubt) node [xshift=0.3cm,midway,left] {\footnotesize  $\set{x,z}(y)\quad$};
 \draw[->] (F)--(Asubt) node [midway,right] {\footnotesize  $\beta$};
    \draw[-] (A)--(P4) node [xshift=-0.3cm,midway,right]  {\footnotesize   $\quad y(\set{x,z})$};
     \draw[->] (A)--(P4) node [midway,left]  {\footnotesize   $\beta$};
    \draw[-] (Asubt)--(P4) node [midway,above] {\footnotesize  $\set{x,y}(z)$};
     \draw[-] (Asubt)--(P4) node [midway,below] {\footnotesize  $\theta$};
   \end{tikzpicture}}
\end{tabular}
\end{center}
According to the criterion,   the  orientation of, say,   the $\beta$ edge connecting $z(x(y))$ and $x(y,z)$ is dictated by the fact that $x$ is below $z$ in $\mathbb{G}({\cal T})$.
The orientation of  the $\theta$ edges is then determined after   choosing a starting vertex (one of the three upper vertices).

\smallskip
 We now embark on the proof of  soundness and completeness of this criterion. We shall formulate the criterion in different   ways, and we shall exhibit the   relationship between  the connectedness properties of ${\cal T}$ and of $\mathbb{G}({\cal T})$.

\smallskip
We first observe that for any two distinct vertices $u,v$ of $\mathbb{G}({\cal T})$, exactly one of the following two situations occurs
(referring to $u$, $v$ as edges of ${\cal T}$):
\begin{itemize}
\item [$\bullet$] Type I: $u$ is above $v$ or conversely.
\item [$\bullet$] Type II: $u$ and $v$ are situated in disjoint branches of a subtree of ${\cal T}$.  We will denote by  ${\it meet}(u,v)$ the node of ${\cal T}$ at which the two branches diverge.
\end{itemize}
We can reformulate these two situations in  $\mathbb{G}({\cal T})$, without reference to ${\cal T}$:
\begin{itemize}
\item [$\bullet$]  Type I: There is a descending path of solid edges (i.e., the level decreases by 1 at each node in the path)   from $u$ to $v$ or from $v$ to $u$ (such a path will be called of type I);
\item [$\bullet$]  Type II: There exists a path  $p=p_1,u',v',p_2$ from $u$ to $v$ whose parts $p_1,u'$ and $v',p_2$  are descending and ascending, respectively (and therefore are made of solid edges only) and which is such that $(u',v')$ is a dashed edge (such a path will\ be called of type II).
\end{itemize}
That this indeed is a reformulation is obvious for type I, while for type II, the desired path in $\mathbb{G}({\cal T})$ is obtained by going down in ${\cal T}$ from (the child vertex of) $u$ all the way down to $u'$ whose father node is  ${\it meet}(u,v)$,  then through a dashed arrow to the branch carrying $v'$, and then all the way up to (the child vertex of) $v$.  Conversely, transcribing the path $p_1,u',v',p_2$ in the language of ${\cal T}$, we find a configuration of type II there.  

\smallskip
In the next lemma, we show how to transform any path into a path of type I or II with the same end nodes.  The transformations are specified by the following picture:
\begin{center}
\resizebox{11cm}{!}{\begin{tabular}{cccccccccccccccccccc}

\begin{tikzpicture}[]
    \node (Z)[] at (1,1) {$z$};
    \node (X)[]  at (-1,1) {$x$};
    \node (Y)[]  at (0,1) {$y$};

 \draw[dashed] (X)--(Y) node {};
 \draw[dashed] (Y)--(Z) node {};
   \end{tikzpicture}
 &&&&

\begin{tikzpicture}
    \node (Y)[] at (-0.5,0) {$y$};
    \node (X)[]  at (-1,1) {$x$};
    \node (Z)[]  at (0,1) {$z$};

 \draw[-] (Y)--(X) node  {};
 \draw[-] (Z)--(Y) node {};
 
   \end{tikzpicture}
&& &&
\begin{tikzpicture}
    \node (Z)[] at (-0.5,0) {$z$};
    \node (X)[]  at (-1,1) {$x$};
    \node (Y)[]  at (0,1) {$y$};
 \draw[-] (Z)--(Y) node {};
 \draw[dashed] (X)--(Y) node {};
   \end{tikzpicture}

 &&&&
\begin{tikzpicture}
    \node (X)[] at (-0.5,0) {$x$};
    \node (Y)[]  at (-1,1) {$y$};
    \node (Z)[]  at (0,1) {$z$};
 \draw[-] (X)--(Y) node  {};
 \draw[dashed] (Z)--(Y) node {};
   \end{tikzpicture}
\\[0.4cm]
$\downarrow$ &&&& $\downarrow$ &&&& $\downarrow$ &&&& $\downarrow$ 
\\[0.4cm]

\begin{tikzpicture}
    \node (X)[]  at (-1,1) {$x$};
    \node (Z)[]  at (0,1) {$z$};
    \node (a)[]  at (0,0) {};

 \draw[dashed] (X)--(Z) node {};
 
   \end{tikzpicture}

&&&&
\begin{tikzpicture}
    \node (X)[]  at (-1,1) {$x$};
    \node (Z)[]  at (0,1) {$z$};
 \node (a)[]  at (0,0) {};

 \draw[dashed] (X)--(Z) node {};
 
   \end{tikzpicture}

&&&&

\begin{tikzpicture}
    \node (Z)[] at (-0.5,0) {$z$};
    \node (X)[]  at (-1,1) {$x$};
 \draw[-] (Z)--(X) node  {};
   \end{tikzpicture}

&&&&
\begin{tikzpicture}
    \node (X)[] at (-0.5,0) {$x$};
    \node (Z)[]  at (0,1) {$z$};
 \draw[-] (X)--(Z) node {};
   \end{tikzpicture}

\end{tabular}}
\end{center}
This specification is then used to define a rewriting system:
$$p_1,x,y,z,p_2 \longrightarrow  p_1,x,z,p_2\vspace{-0.1cm}$$
when $x,y,z$ are in one of the four configurations at the top of the picture.

\begin{lemma}  \label{equivalent-criteria}
This rewriting system is confluent and terminating. It is complete in the sense that any two paths between the same pair of end points are provably equal by a zigzag of such rewritings, and sound in the sense that any such zigzag always relates two paths with the same endpoints.
 The  normal forms of the rewriting system are  the paths of type I or II, and  are the min-paths.
\end{lemma}
\Proof  Termination is obvious, since  the length decreases by 1 at each step.  As for confluence, we list the critical pairs, which all admit immediate solutions (note that the sequence $(x,y)$ solid, $(y,z)$ dashed, $(z,u)$ solid is excluded since one would then have $x=u$, which contradicts the definition of a path):
\begin{center}

\resizebox{3.35cm}{!}{\begin{tikzpicture}
    \node (Z)[] at (1,1) {$z$};
    \node (X)[]  at (-1,1) {$x$};
    \node (Y)[]  at (0,1) {$y$};
    \node (U)[]  at (2,1) {$u$};
\draw[dashed] (Z)--(U) node {};
 \draw[dashed] (X)--(Y) node {};
 \draw[dashed] (Y)--(Z) node {};
   \end{tikzpicture}}
\vspace{0.25cm}

\resizebox{11.4cm}{!}{\begin{tabular}{cccc}
\begin{tikzpicture}
    \node (Y)[] at (-0.5,0) {$y$};
    \node (X)[]  at (-1,1) {$x$};
    \node (Z)[]  at (0,1) {$z$};
 \node (U)[]  at (1,1) {$u$};
\draw[dashed] (Z)--(U) node {};
 \draw[-] (Y)--(X) node  {};
 \draw[-] (Z)--(Y) node {};
 
   \end{tikzpicture}\enspace
& \begin{tikzpicture}
    \node (Z)[] at (-0.5,0) {$z$};
    \node (Y)[]  at (-1,1) {$y$};
    \node (U)[]  at (0,1) {$u$};
 \node (X)[]  at (-2,1) {$x$};
\draw[dashed] (X)--(Y) node {};
 \draw[-] (Y)--(Z) node  {};
 \draw[-] (Z)--(U) node {};
 
   \end{tikzpicture}\enspace &
\begin{tikzpicture}
    \node (U)[] at (-0.5,0) {$u$};
    \node (Y)[]  at (-1,1) {$y$};
    \node (Z)[]  at (0,1) {$z$};
    \node (X)[]  at (-2,1) {$x$};
 \draw[-] (U)--(Z) node {};
 \draw[dashed] (Y)--(Z) node {};
 \draw[dashed] (X)--(Y) node {};
   \end{tikzpicture}\enspace
& \begin{tikzpicture}
    \node (X)[] at (-0.5,0) {$x$};
    \node (Y)[]  at (-1,1) {$y$};
    \node (Z)[]  at (0,1) {$z$};
    \node (U)[]  at (1,1) {$u$};
 \draw[-] (X)--(Y) node  {};
 \draw[dashed] (Z)--(Y) node {};
 \draw[dashed] (Z)--(U) node {};
   \end{tikzpicture}

\end{tabular}}
\end{center}
That the paths of type I and II are in normal form is also immediate (there is no matching for the left hand sides of our rewriting rules).  It remains to check that all normal forms are indeed of one of these two shapes. We proceed by induction on the length of the normal form $p$.  Every path of length 1 is indeed of type I or II.   Let now $p=u,v,p_1$. We can assume by induction that $p_1$ is of type I or II.  There are three cases:
\begin{itemize}
\item $(u,v)$ is solid with $v$ one level up from $u$.  Then $p_1$ cannot start with a solid edge going down, because then $p$ would visit $u$ twice, nor with a dashed edge, because $p$ would then not be a normal form. Hence $p_1$ is of type I, and morevoer goes up (again because otherwise $p$ would not be a path).  Then adding $(u,v)$ in front still results in a path of type I.
\item $(u,v)$ is solid with $v$ one down  from $u$.  Then $p_1$ cannot start with a solid edge going up, since $p$ would not be in normal form.  Hence prefixing $p_1$ with $(u,v)$ yields a path of type I (resp. II) if $p_1$ was of type I (resp. II).
\item $(u,v)$ is dashed. Then $p_1$ cannot start with a dashed edge nor a solid edge going down, as $p$ would then not be in normal form. Hence $p_1$ has to be of type I, going up, which makes $p$ a path of type II. 
\end{itemize}
We now prove completeness. We have already observed the uniqueness of the paths of type I or II. Since we have established that the normal forms are the paths of type I or II,  it  follows that all paths in normal form from $u$ to $v$ coincide  (notation $\equiv$), and we have, for any two paths $p_1,p_2$ from $u$ to $v$, and writing ${\it nf}(p)$ for the normal form of a path $p$\vspace{-0.1cm}
$$ p_1 \longrightarrow^* {\it nf}(p_1)\equiv {\it nf}(p_2) {}^*\longleftarrow p_2\vspace{-0.1cm}$$
Conversely, the rewriting system leaves the endpoints of the path unchanged at each step, and hence any zigzag maintains this inviariant, which establishes soundness.

That every minpath is normal is clear, since any rewriting step decreases the length of a path.  For the converse, we use completeness. Suppose that $p$ is normal, but is not a min-path, and let $p_1$ be a min-path with the same endpoints as $p$.   By completeness, there exists a zigzag between $p$ and $p_1$, or equivalently, by confluence, $p$ and $p_1$ have the same normal form. But ${\it nf}(p_1)$ has a fortiori a length strictly smaller than $p={\it nf}(p)$: contradiction.   \qed

\medskip
Summing up, the following are characterisations of ``being of type I or II", for  two distinct vertices $u,v$ of $\mathbb{G}({\cal T})$ (or equivalently, two edges $u,v$ of ${\cal T}$):
\begin{center}
\begin{tabular}{ll}
Type I & Type II\\\\
$u,v$ are one above the other in ${\cal T}$  &  $u,v$ are on disjoint  branches \\
& of a subtree of ${\cal T}$\\\\
$u,v$ connected by a path of type I  & $u,v$ connected by a path of type II\\\\
 min-path between $u,v$ is of type I &  min-path between $u,v$ is of type II\\\\
 min-path between $u,v$ contains  &  min-path between $u,v$  contains \\
  only solid edges &  at least one dashed edge
\end{tabular}
\end{center}
Indeed, by Lemma \ref{equivalent-criteria}, we know that the min-paths are exactly the paths of type I or II, and crossing or not a dashed edge is what distinguishes among min-paths those that are of type II or I, respectively.

\begin{lemma}\label{dds} There is a one-to-one correspondence between the subtrees of ${\cal T}$ and the connected subsets of $\mathbb{G}({\cal T})$.
\end{lemma}
\Proof The connected subset of $\mathbb{G}({\cal T})$ corresponding to a subtree ${\cal T}'$ of ${\cal T}$ is precisely $\mathbb{G}({\cal T}')$. In the other direction, let
  $K$ be  a connected subset of $\mathbb{G}({\cal T})$. By connectedness, for any $u,v$ in $K$ there exists a path $p$ from $u$ to $v$ that is included in $K$.  By Lemma \ref{equivalent-criteria}, we know that ${\it nf}(p)$ is  also included in $K$.  By this observation,  through the transcription in ${\cal T}$  of paths of type I or II of $\mathbb{G}({\cal T})$,  we conclude that the subgraph of ${\cal T}$  whose edges are precisely the vertices of $K$ is a subtree.  \qed


\medskip
 In what follows, in the context of operadic trees, we shall  say that a tree is {\em non-Empty} if it contains at least one edge (whence the capital ``E").
 Clearly, all operadic trees relevant for describing operadic laws are non-Empty.
\begin{lemma}\label{removing_n_edges}
If $x_1,\dots,x_n$ are arbitrary distinct edges of ${\cal T}$, then the following  claims hold.
\begin{enumerate}
\item By removing $x_1,\ldots,x_n$ from ${\cal T}$, we obtain exactly $n+1$ subtrees of ${\cal T}$.
\item The number $k$ of non-Empty subtrees of ${\cal T}$ obtained in this way is equal to the number of connected components of $\mathbb{G}({\cal T})$ obtained by removing the vertices $x_1,\dots,x_n$, and $k\in\{0,\dots,n+1\}$.
\item Let ${\cal T}'$ be one of the non-Empty subtrees of ${\cal T}$ obtained by removing  $x_1,\dots,$ $x_n$, and let $K$ be the  connected subset of $\mathbb{G}({\cal T})$ associated with ${\cal T}$ by (2). Then, if $y$ is an edge of ${\cal T}'$ and a vertex of $K$, we have that $\mathbb{G}({\cal T}')=K$.
\end{enumerate} 
\end{lemma}
\Proof We consider only the case $n=1$ (the general case  follows easily by induction), and we write $x$ for $x_1$.
Let $a$ and $b$ be the vertices adjacent to $x$, with $a$ being the child vertex for $b$.

The first claim is standard: the subtrees obtained after the removal of $x$ are  the subtree ${\cal T}_1$ rooted at $a$ and containing all descendants of $a$, and the subtree ${\cal T}_2$ obtained from ${\cal T}$ by removing all of  ${\cal T}_1$.
Note that $b$ is a leaf of ${\cal T}_2$.

We prove the other two claims in parallel, by looking at the possible configurations of ${\cal T}$.
 If $x$ is the only edge of ${\cal T}$, then ${\cal T}_1$ and  ${\cal T}_{2}$ are the vertex $a$ and the vertex $b$, respectively, and, therefore, $k=0$. Suppose that $x$ is not the only edge of ${\cal T}$. Then, if $x$ is on the highest level in ${\cal T}$,  ${\cal T}_1$ is just  the vertex $a$, while $
{\cal T}_2$ is clearly non-Empty, and, hence, $k=1$. We also get $k=1$ when $x$ is the unique edge on the first level of ${\cal T}$, in which case ${\cal T}_1$ is non-Empty and ${\cal T}_2$ is just the vertex $b$. In all other situations, we have $k=2$.

\smallskip
Let us now prove that $k$ is also the number of connected components of $\mathbb{G}({\cal T})$ obtained by removing the vertex $x$. We examine only the case  $k=2$. Let $K_1=\mathbb{G}({\cal T}_1)$ and $K_2=\mathbb{G}({\cal T}_2)$.
Since $K_1$ and $K_2$ are connected and disjoint and  $\mathbb{G}({\cal T})\backslash\{x\}=K_1\cup K_2$, we only have to show that the set of edges of $\mathbb{G}({\cal T})\backslash\{x\}$ is also the (disjoint) union of sets of edges of  $K_1$ and $K_2$. For this, we use the fact that  the removal of $x$ from $\mathbb{G}({\cal T})$  involves the removal of all edges of $\mathbb{G}({\cal T})$ that have $x$ as one of its adjacent vertices. Let $e$ be  an edge of $\mathbb{G}({\cal T})$, with $y$ and $z$ being its adjacent vertices. 

Suppose first that $e$ is an edge of $\mathbb{G}({\cal T})\backslash\{x\}$. We then know that both $y$ and $z$ are different from $x$, and share a common vertex $v$ when considered as edges of ${\cal T}$. Since ${\cal T}_1$ and ${\cal T}_2$ form a partition of the set of vertices of ${\cal T}$, let us  assume, say, that $v$ is a vertex of ${\cal T}_1$. If $v\neq a$, we can immediately conclude that $y$ and $z$ are edges of ${\cal T}_1$, and, if $v=a$, then, since both $y$ and $z$ are different from $x$, it must be the case that $v$ is the parent vertex for both $y$ and $z$, which also implies that $y$ and $z$ are edges of ${\cal T}_1$. Therefore, $y$ and $z$ are both vertices of $K_1$, and, hence, $e$ is an edge of $K_1$.  

Conversely, if $e$ is an edge of $K_1$, then $y$ and $z$ are edges of ${\cal T}_1$, and therefore must both be different from $x$. Since they share a common vertex in ${\cal T}_1$, and, hence, in ${\cal T}$, we  conclude that $e$ is an edge of $\mathbb{G}({\cal T})\backslash\{x\}$. \qed




\medskip 

\medskip
The following  proposition is only implicit in \cite{DP-WCO}.
\begin{proposition}\label{constructions-insertions}
For every operadic tree ${\cal T}$, the constructions of $\mathbb{G}({\cal T})$ (considered as hypergraph) are in one-to-one correspondence with the (fully) parenthesised words that denote  decompositions of ${\cal T}$.
\end{proposition}
\Proof  To every decomposition/parenthesisation of ${\cal T}$, one can associate a tree each of whose nodes is decorated by an edge of ${\cal T}$: one proceeds from the most internal parentheses to the most external ones, recording each insertion on the way.  

Formally, the fullly parenthesised words  are declared by  the syntax $w:: a\Alt (ww)$, where $a$ ranges over the nodes of ${\cal T}$ (all named with different letters).

 Not all words correspond to decompositions of ${\cal T}$. When this is the case, we say that $w$ is admissible for ${\cal T}$ (the precise definition of admissibility can be easily reconstructed from the inductive construction below).
 
\smallskip
Since we deal with non-Empty trees, our base case is that of a word $(ab)$ corresponding to a single edge operadic tree connecting $a$ and $b$.  Then there is only one decomposition and one construction, hence the statement holds. 

Otherwise, we have a word 
$w=(w_1w_2)$, where  at least one of  the words $w_1$ or $w_2$ is not reduced to a letter. We proceed by structural induction on $w$, providing both the decorated tree and the proof that is indeed a construction.
Let us  call ${\cal T}_1$, ${\cal T}_2$ the trees decomposed by $w_1$, $w_2$, respectively (cf. Lemma \ref{removing_n_edges}). Let $x$ be the edge on which  ${\cal T}_1$ is grafted on the tree ${\cal T}_2$. 
We distinguish three cases.
\begin{enumerate}
\item  If  neither $w_1$ nor $w_2$ are  reduced to a letter, then $\mathbb{G}({\cal T}_1)$ and $\mathbb{G}({\cal T}_2)$ are both non-empty, and 
are the connected components of 
$\mathbb{G}({\cal T})\setminus \set{x}$. 
We can thus apply induction: if $V_1$ and $V_2$ are the constructions associated with 
$w_1$ and $w_2$,  then we associate  $x(V_1,V_2)$   with $(w_1w_2)$, which is  a construction.
\item If $w_2=a$ is a reduced to a letter and $w_1$ is not reduced to a letter, then $\mathbb{G}({\cal T}_2)$ is empty, and $x$ is a leaf of $\mathbb{G}({\cal T})$. We conclude  by induction that  the tree $x(V_1)$ associated with $((w_1)a)$ is a construction.
\item If $w_1=a$ is reduced to a letter, then ${\cal T}$ is of the form $a({\cal T}_2)$, i.e. $a$ is the root and has only one child which is the root of ${\cal T}_2$.  
We  conclude  by induction that the tree $x(V_2)$ associated with $(a(w_2))$ is a construction.
\end{enumerate}
Note that case 1  (resp. cases 2 and 3) correspond to the situation where $k=2$ (resp. $k=1$), while the base case is the case where $k=0$ (in the terminology of Lemma \ref{removing_n_edges}).

\smallskip
The converse mapping is defined much in the same way. We observe that, for any ${\cal T}$ (with at least 3 nodes), constructions of $\mathbb{G}({\cal T})$ can only be of the form
$x(V)$ or $x(V_1,V_2)$.  They are of the first (resp. second)  form when the node $x$ is either a leaf or the root of $\mathbb{G}({\cal T})$, (resp. when $x$ is any other node).   We can deploy induction on the number of nodes of {\cal T}  and map constructions back to parenthesised words.  By induction too, we can show that these are inverse transformations. \qed

\smallskip As an illustration, referring to Figure \ref{G(T)-picture}, 
 $(ae)((bd)c)$ is  mapped to $z(x(y),u)$, obtained as follows:  the leaf $u$ records  $ae$, while in parallel the leaf $y$ records $bd$ and then $x(y)$ encodes $(bd)c$ and, finally, the last insertion is along $z$.  (Note that, following common practice, in examples, we do not write the most external parentheses.)


\medskip
We are now in a position to conclude.
\begin{theorem}  The criterion $(\dagger)$ is sound and complete.
\end{theorem}

\Proof  We write $u,v$ for $u_E,v_E$. Let $E'$ be the subtree of $E$ whose root is $\{u,v\}$ and let $K$ be the connected subset of $\mathbb{G}({\cal T})$ out of which $E'$ arises as a construct. Let ${\cal T}'$ be the subtree of ${\cal T}$ that corresponds to $K$ by Lemma \ref{dds}.\\
\indent The  number of constructions grafted to $\{u,v\}$ in $E'$ is  the number of connected components of $\restrH{K}{\{u,v\}}$. By  Lemma \ref{removing_n_edges}, it is also the number of non-Empty subtrees of ${\cal T}'$ obtained by removing the edges $u$ and $v$. Moreover, there can be at most $3$ such subtrees. Let us now introduce some names. \\
\indent  Let ${\cal T}'_1, {\cal T}'_2$ and ${\cal T}'_3$ be  the subtrees of ${\cal T}'$ obtained by removing $u$ and $v$. Let $I\subseteq\{1,2,3\}$ be  such that $i\in I$ if and only if ${\cal T}'_i$ is non-Empty, and let $J=\{1,2,3\}\backslash I$. Let, for all $i\in I$, $K_i$ be the  connected component of $\restrH{K}{\{u,v\}}$ corresponding to ${\cal T}'_i$, $V_i$ be the construction of $K_i$ that is grafted to $\{u,v\}$ in $E'$, and $w_i$ be the decomposition of ${\cal T}'_i$  corresponding to   $V_i$ according to  Proposition \ref{constructions-insertions}. On the other hand, for all $j\in J$, ${\cal T}'_j$ is a vertex $a_j$, and let $w_j$ be precisely  $a_j$.

By analysing the case $n=2$ of Lemma \ref{removing_n_edges}, we get that  $E'$
determines an (incomplete) decomposition $W$ of ${\cal T}'$, in which the insertions of  $u$ and $v$ are the only ones not yet performed, and which has one of the following two shapes:
\begin{center}
\begin{tabular}{ccc}
\resizebox{0.85cm}{!}{ \begin{tikzpicture}
    \node (E)[circle,draw=black,minimum size=7mm,inner sep=0.1mm] at (0,0) { \footnotesize $w_{k_1}$};
    \node (F) [circle,draw=black,minimum size=7mm,inner sep=0.1mm] at (0,1.25) {\footnotesize $w_{k_2}$};
    \node (A) [circle,draw=black,minimum size=7mm,inner sep=0.1mm] at (0,2.5) {\footnotesize $w_{k_3}$};
    \draw[-] (E)--(F) node  [midway,left] {\footnotesize $u$};
    \draw[-] (F)--(A) node  [midway,left] {\footnotesize $v$};
\end{tikzpicture}}
&  \enspace \enspace \enspace
\begin{tikzpicture}
 \node (E)[circle,minimum size=8mm,inner sep=0.1mm] at (0,0.2) {};
    \node (F) [circle,minimum size=8mm,inner sep=0.1mm] at (0,1.4) {$\mbox{and}$};
\end{tikzpicture}
 & \enspace \enspace
\resizebox{2.8cm}{!}{ \begin{tikzpicture}
\node (E)[circle,minimum size=10mm,inner sep=0.1mm] at (0,0.52) {};
    \node (E)[circle,draw=black,minimum size=7mm,inner sep=0.1mm] at (0,0.9) {\footnotesize $w_{k_2}$};
    \node (F) [circle,draw=black,minimum size=7mm,inner sep=0.1mm] at (-1,2.4) {\footnotesize $w_{k_1}$};
    \node (A) [circle,draw=black,minimum size=7mm,inner sep=0.1mm] at (1,2.4) {\footnotesize $w_{k_3}$};
    \draw[-] (E)--(F) node  [midway,left] {\footnotesize $u$};
    \draw[-] (A)--(E) node  [midway,right] {\footnotesize $v$};
\end{tikzpicture}}
\end{tabular}
\end{center}
where $\{k_1,k_2,k_3\}=I\cup J$, and where the words $w_{k_i}$ are as defined above. The  shape on the left arises in the case when there exists a sequence $u=x_0,\dots,x_n=v$ of edges in ${\cal T}'$ such that the child vertex of $x_{i-1}$ is a parent vertex of $x_i$, for all $1\leq i \leq n$, and the  one on the right when there exists a subtree of ${\cal T}'$ that has $u$ and $v$ on different branches.
We observe that $\set{u,v}$ is of type I  (resp. of type II) in ${\cal T}'$ (and hence in ${\cal T}$) if $W$ has the shape on the left (resp. on the right).

Now, if $V_1$ and $V_2$ are the vertices of $P_{\cal T}$ adjacent to $E$, then, in order to get complete decompositions of ${\cal T}'$ corresponding to $V_1$ and $V_2$, it remains  to add $u$ and $v$ (in a way dictated by $V_1$ and $V_2$, respectively) in the sequence of insertions obtained previously from $E'$.
More precisely, if we assume that    
$u$ is the child of $v$  (resp. $v$ is the child of $u$) in
 $V_1$ (resp. $V_2$),  then in the decomposition of ${\cal T}'$ guided by $V_1$ (resp. $V_2$), the insertion of $v$ (resp. $u$) will be applied last.  We then conclude by examining the two  possible shapes of $W$.
\begin{itemize} 
\item In the type I case,  $V_1$ and $V_2$ differ only by the subwords 
$(w_{k_1} w_{k_2})w_{k_3}$ and  $w_{k_1} (w_{k_2}w_{k_3})$, respectively. Hence $E'$ features a $\beta$-arrow.
Moreover, the orientation prescribed in the statement of the criterion  tells us that the edge should be oriented from $V_1$ to $V_2$, given our (arbitrary) choice of placing $u$ under $v$ in our drawing on the left.
 \item In the type II case, $V_1$ and $V_2$ differ only by the subwords
$(w_{k_1} w_{k_2})w_{k_3}$ and  $(w_{k_1} w_{k_3})w_{k_2}$, respectively. Hence $E'$ features a $\theta$-arrow. \qed
\end{itemize}
We illustrate the constructions of the proof below, with $E',K,{\cal T}'$ as follows:
\begin{center}
\begin{tabular}{ccccc}
\begin{tikzpicture}
    \node (Z)[] at (0,0) { \footnotesize $\{u,v\}$};
      \node (V)[]  at (-0.5,1) {\footnotesize $x$};
        \node (y)[]  at (0.5,1) {\footnotesize $y$};
\node (z)[]  at (0,2) {\footnotesize $z$};
\node (w)[]  at (1,2) {\footnotesize $w$};
    \node (A)[] at (-1.2,1) {$E'=$};
\node (f)[] at (-1.2,-1.03cm) { };
 \draw[-] (Z)--(y) node  {};
 \draw[-] (Z)--(V) node  {};
 \draw[-] (y)--(w) node  {};
  \draw[-] (z)--(y) node  {};
\end{tikzpicture} &    \begin{tikzpicture}
    \node (Z)[] at (-2,0) {\footnotesize $u$};
      \node (V)[]  at (-2.5,1) {\footnotesize $x$};
        \node (y)[]  at (-1.5,1) {\footnotesize $v$};
\node (z)[]  at (-2,2) {\footnotesize $z$};
\node (w)[]  at (-1,2) {\footnotesize $y$};
\node (r)[]  at (-1,3) {\footnotesize $w$};
    \node (A)[] at (-3.2,1.5) {$K=$};
 \node (f)[] at (-1.2,-0.55cm) { };
 \draw[-] (Z)--(y) node  {};
 \draw[-] (Z)--(V) node  {};
 \draw[dashed] (y)--(V) node  {};
 \draw[-] (z)--(y) node  {};
 \draw[-] (w)--(y) node  {};
\draw[-] (w)--(r) node  {};
  \draw[dashed] (z)--(w) node  {};
\end{tikzpicture}
&
 \begin{tikzpicture}
    \node (A)[circle,draw=black,minimum size=4mm,inner sep=0.1mm] at (0,-1) {\footnotesize $a$};
    \node (B)[circle,draw=black,minimum size=4mm,inner sep=0.1mm] at (0,0) {\footnotesize $b$};
      \node (V)[circle,draw=black,minimum size=4mm,inner sep=0.1mm]  at (-0.5,1) {\footnotesize $c$};
        \node (y)[circle,draw=black,minimum size=4mm,inner sep=0.1mm]  at (0.5,1) {\footnotesize $d$};
\node (z)[circle,draw=black,minimum size=4mm,inner sep=0.1mm]  at (0,2) {\footnotesize $e$};
\node (w)[circle,draw=black,minimum size=4mm,inner sep=0.1mm]  at (1,2) {\footnotesize $f$};
\node (r)[circle,draw=black,minimum size=4mm,inner sep=0.1mm]  at (1,3) {\footnotesize $g$};
    \node (L)[] at (-1.5,1) {${\cal T}'=$};
 \draw[-] (B)--(y) node  [midway,right]   {\footnotesize $v$};
 \draw[-] (B)--(V) node  [midway,left]  {\footnotesize $x$};
 \draw[-] (z)--(y) node [midway,left]  {\footnotesize $z$};
 \draw[-] (w)--(y) node [midway,right]   {\footnotesize $y$};
\draw[-] (w)--(r) node [midway,right] {\footnotesize $w$};
\draw[-] (A)--(B) node  [midway,right]   {\footnotesize $u$};
\end{tikzpicture}

\end{tabular}
\end{center}
The subtrees we get after removing $u$ and $v$ from ${\cal T}'$ are
\begin{center}
\begin{tabular}{ccc}
\begin{tikzpicture}[scale=0.8]
    \node (L)[] at (-1,0) {${\cal T}'_1=$};
 \node (f)[] at (-1,-1.3) { };
    \node (A)[circle,draw=black,minimum size=4mm,inner sep=0.1mm] at (0,0) {\footnotesize $a$};
\end{tikzpicture} \enspace\enspace & \enspace\enspace
 \begin{tikzpicture}
    \node (B)[circle,draw=black,minimum size=4mm,inner sep=0.1mm] at (0,0) { \footnotesize $b$};
      \node (V)[circle,draw=black,minimum size=4mm,inner sep=0.1mm]  at (0,1) {\footnotesize $c$};
       \node (L)[] at (-1.2,0.5) {${\cal T}'_2=$};
\node (f)[] at (-1,-0.55) { };
 \draw[-] (B)--(V) node  [midway,left]  {\footnotesize $x$};
\end{tikzpicture} \enspace\enspace
 &  \enspace\enspace
 \begin{tikzpicture}
        \node (y)[circle,draw=black,minimum size=4mm,inner sep=0.1mm]  at (0,0) {\footnotesize $d$};
\node (z)[circle,draw=black,minimum size=4mm,inner sep=0.1mm]  at (-0.5,1) {\footnotesize $e$};
\node (w)[circle,draw=black,minimum size=4mm,inner sep=0.1mm]  at (0.5,1) {\footnotesize $f$};
\node (r)[circle,draw=black,minimum size=4mm,inner sep=0.1mm]  at (0.5,2) {\footnotesize $g$};
    \node (L)[] at (-1.5,1) {${\cal T}'_3=$};
 \draw[-] (z)--(y) node [midway,left]  {\footnotesize $z$};
 \draw[-] (w)--(y) node [midway,right]   {\footnotesize $y$};
\draw[-] (w)--(r) node [midway,right] {\footnotesize $w$};
\end{tikzpicture}
\end{tabular}
\end{center}
and the corresponding decompositions are $$w_1=a\enspace\enspace\enspace \enspace\enspace\enspace\enspace\enspace\enspace w_2= bc\enspace\enspace\enspace\enspace\enspace\enspace\enspace\enspace\enspace w_3=(de)(fg)$$
Hence, $E'$ corresponds to the following decomposition of ${\cal T}'$:
\begin{center}
 \begin{tikzpicture}
   \node (E)[circle,draw=black,minimum size=4mm,inner sep=0.1mm] at (0,0.5) {\footnotesize $a$};
   \node (F) [circle,draw=black,minimum size=6mm,inner sep=0.1mm] at (0,1.55) {\footnotesize $bc$};
  \node (A) [circle,draw=black,minimum size=10mm,inner sep=0.1mm] at (0,3) {\footnotesize $(de)(fg)$};
    \node (B) [circle,minimum size=10mm,inner sep=0.1mm] at (-1.5,1.5) {$W=$};
    \draw[-] (E)--(F) node  [midway,left] {\footnotesize $u$};
    \draw[-] (F)--(A) node  [midway,left] {\footnotesize $v$};
\end{tikzpicture}
\end{center}
For this example, the vertices  $V_1$ and $V_2$  adjacent to $E$ induce decompositions $(a(bc))((de)(fg))$   and  $a((bc)((de)(fg)))$,  respectively, and $E$ features a $\beta$-arrow from $V_1$ to $V_2$.

\vspace{-.4cm}

 \section{Iterated truncations} \label{stretching-section}

In this section, we recast the iterated truncations of \cite{P-SISP} in our setting. Hypergraph polytopes allow us to describe all the polytopes in the interval between the simplex and the permutohedron. The hypergraph specifies at once all truncations to be made to reach a particular polytope in this interval. But what about truncating a new face that was not present in the original simplex, i.e. a face already obtained as a result of a truncation?
We shall build a whole ``tree'' of polytopes,  each polytope in the tree giving rise to a whole interval of truncations which are all its child nodes in the tree. A polytope at distance $n$ of the root  will be obtained through $n$ runs of truncations.  The root is occupied by the simplex.   Our tree notation for constructs extends  to this setting.


\vspace{-.3cm}
\subsection{Successive rounds of truncations}

\vspace{-.1cm}

Let ${\cal X}$ be a set (whose elements stand for the facets of the initial simplex).
All the work will be carried out within ${\cal M}^f({\cal X})$, the set of finite multisets of elements of ${\cal X}$, which gives rise to a monad.
At each round of truncation, we are given 
\begin{itemize}
\item 
a non-empty set $H\inc {\cal M}^f({\cal X})$ (whose elements  stand for the facets of the polytope that is to be truncated);
\item
a hypergraph   $\hyper{H}^v$  (whose hyperedges stand for the \boldmath{$v$}ertices of the same polytope); 
\item 
an atomic and connected hypergraph  $\hyper{H}^t$
(whose connected subsets give  instructions for the \boldmath{$t$}runcations to be performed at this round).
\end{itemize}
We require that
$(\union \hyper{H}^v) = H = (\union \hyper{H}^t).$ We also require  $\hyper{H}^v$ to satisfy the following property::\\[0.2cm]
 \enspace $(P)$ \quad $\forall x\in H,\:\exists V\in\hyper{H}^v,\:  x\in V\;.$\\[0.2cm]
\indent 
%
%
The intuition is that $\hyper{H}^v$ serves to tame the constructs of $\hyper{H}^t$. Indeed, the polytopes that we are building in this way are simple (this is a consequence of \cite{P-SISP}[Proposition 9.3]), so from every vertex at round $n$, the local view is that of a simplex -- a property that  is often taken as the definition of simplicity --, which makes those polytopes liable to the machinery of hypergraph polytopes.

We modify the  definition of construct (and construction) as follows.  Constructs are defined exactly as in Section \ref{inductive-section}, except for the root, for which one has to pick, not an arbitrary non-empty subset $Y$ of $H$, but  one which contains the complement of some hyperedge $V$ of $\hyper{H}^v$.  
Such a construct will be called a {\em construct of $\hyper{H}^t$ rel to $\hyper{H}^v$}, and we shall say that it is  {\em tamed} by $V$.  Constructs in the ``old sense'' will be called plain constructs.
 Here is the full definition.  Pick an arbitrary subset  $Y\inc H$ such that $(H\setminus V)\inc Y$ for some $V\in \hyper{H}^v$.
\begin{itemize}
\item 
If $Y = H$, then  the one node tree decorated with $Y$, and written $Y$, is a construct of $\hyper{H}^t$ rel to $\hyper{H}^v$.
\item  Otherwise, if   $\hyper{H}^t ,Y  \leadsto H_1,\ldots,H_n$, and if  $T_1,\ldots,T_n$ are plain constructs of $\hyper{H}_1,\ldots,\hyper{H}_n$, respectively, then
$Y(T_1,\ldots,T_n)$, is a construct of $\hyper{H}^t$ rel to $\hyper{H}^v$. 
\end{itemize}
Note that the taming is only performed at the root. 
We denote by  ${\cal A}_{\hyper{H}^v}(\hyper{H}^t)$  the set of  constructs  rel to $\hyper{H}^v$.  

The definition of construction is also slightly modified: it is a tree where all non-root nodes are decorated by singletons while the root is decorated exactly  by the complement of some hyperedge of $\hyper{H}^v$.

\smallskip
The initial round of truncations is along the simplex-permutohedron interval.
We take:
\begin{itemize}
\item $H_1={\cal X}$ (identifying an element $x$ of ${\cal X}$ with the associated one element multiset);
\item $\hyper{H}_1^v=\setc{{\cal X}\setminus\set{x}}{x\in{\cal X}}$;
\item  $\hyper{H}_1^t$ is any atomic connected hypergraph on $H_1$.
\end{itemize}
Note that the constructs of $\hyper{H}_1^t$ rel to $\hyper{H}_1^v$ are all the constructs of $\hyper{H}_1^t$   (no taming yet).

\smallskip
We explain now how round $n+1$ is prepared from round $n$.
From $H_n\inc {\cal M}^f({\cal X})$, $\hyper{H}_n^v$, $\hyper{H}_n^t$, we generate 
${\cal A}_{\hyper{H}_n^v}(\hyper{H}_n^t)$, which induces  a set  $H_{n+1}\inc {\cal M}^f({\cal X})$ and a hypergraph $ \hyper{H}_{n+1}^v$ on $H_{n+1}$,  as follows:
\begin{itemize}
\item  The maximal elements of  ${\cal A}_{\hyper{H}_n^v}(\hyper{H}_n^t) \setminus\set{H_n}$,
which we shall   call  {\em constrs}, 
are all of the form   $X(Y)$, 
where $X \union Y = H_n$  (by definition of constructs), and hence are entirely characterized by $Y$.  We set  
$$H_{n+1}= \{(\mu_{{\cal X}}\circ\sigma_{{\cal M}^f({\cal X})})(Y)\,|\, (H\setminus Y)(Y) \;\mbox{is a constr of}\; \hyper{H}_n^t\; \mbox{rel to}\; \hyper{H}_n^v \},$$
where $\sigma$  turns a set into the formal sum of its elements, 
and where $\mu$ is the multiplication of the monad ${\cal M}^f$.
\item  $\hyper{H}_{n+1}^v$ is in bijection with  the set of constructions of $\hyper{H}_n^t$ rel to $\hyper{H}_n^v$:
%
 $$ \begin{array}{rcl}
\hyper{H}_{n+1}^v&   = &  
\{ {\cal P}(\mu_{{\cal X}}\circ\sigma_{{\cal M}^f({\cal X})})(\psi(T)\setminus\set{H_n}) \, |\, 
 \\
 && 
T \;\mbox{is a construction of}\; \hyper{H}_n^t \;\mbox{rel to}\; \hyper{H}_n^v\} \;.
\end{array}$$
\end{itemize}

\begin{proposition}  \begin{enumerate}
\item  $\hyper{H}_{n+1}^v$ is indeed a subset of ${\cal P}(H_{n+1})$, and 
$\mu_{{\cal X}}\circ\sigma_{{\cal M}^f({\cal X})}$ is bijective on the subsets to which it is applied.
\item At every round, $\hyper{H}_n^v$ satisfies property $(P)$.
\item  At every round, we have $H_n\inc H_{n+1}$.
\end{enumerate}
\end{proposition}
\Proof  For (1), we refer (mutatis mutandis) to  \cite{P-SISP}[section 7].
We prove (3) first.  In terms of constructs, we have to show that for each $y\in H_n$, 
$(H_n\setminus\set{y})(y)$ is a construct rel to $\hyper{H}_n^v$. Unrolling what it means, we see that this is the 
case if there exists $V$ such that $(H_n\setminus \set{y})\supseteq (H\setminus V)$ which is statement (2) (at round $n$).  We now prove (2) at round $n+1$.  Let $Y$ be such that $(H_{n+1}\setminus Y)(Y)$ is a constr of $\hyper{H}_n^t$ rel to $\hyper{H}_n^v$.  By an easy adaptation of the devices described in Section \ref{vertices-of-faces-section}, we can obtain at least one construction $V$ of $\hyper{H}_n^t$ rel to $\hyper{H}_n^v$  such that $V\leq (H\setminus Y)(Y)$, which entails $Y\in \psi(V)$. \qed

\medskip
As an illustration, here is how to recast the example  of \cite{P-SISP}[p. 11]. We take ${\cal X}=\set{x,y,z,u}$. We set:
$$\begin{array}{l}
H_1=\set{x,y,z,u}  \quad (\mbox{considered as a subset of}\; {\cal M}^f({\cal X}))\\
\hyper{H}_1^v=\set{\set{x,y,z} , \set{y,z,u} , \set{z,u,x} , \set{u,x,y}}\\
\hyper{H}_1^t = \set{\set{x} , \set{y} , \set{z} , \set{u} , \set{x,y} , \set{x,y,z,u}}
\end{array}$$
resulting in the following truncation of the 3-dimensional simplex (decorating  the vertices as in $\hyper{H}_2^v$, by anticipation):
\begin{center}
\begin{center}
\begin{tikzpicture}[thick,scale=2.7]
\coordinate (A11) at (0.285,0.8);
\coordinate (A12) at (-0.285,0.8);
\coordinate (A21) at (0.285,-0.8); 
\coordinate (A22) at (-0.285,-0.8); 
\node (dl) at (-1.25,0) {\small $\set{x,z,u}$};
\node (dr) at (1.25,0) {\small $\set{y,z,u}$};
\node (a22) at (-0.45,-0.925) {\small $\set{x,x+y,z}$}; 
\node (a21) at (0.45,-0.925) {\small $\set{y,x+y,z}$}; 
\node (a11) at (-0.45,0.925) {\small $\set{x,x+y,u}$}; 
\node (a12) at (0.45,0.925) {\small $\set{y,x+y,u}$}; 
\coordinate (A3) at (-1,0);
\coordinate (A4) at (1,0);
\draw[draw=gray,fill=pur,opacity=0.15] (A11) -- (A12) -- (A3) -- (A22) -- (A21) -- (A4) -- cycle;
\draw[draw=gray] (A11) -- (A12) -- (A3) -- (A22) -- (A21) -- (A4) -- cycle;
\draw[draw=gray,fill=yellow,opacity=0.1] (A11) -- (A12) -- (A22) -- (A21) -- cycle;
\draw[draw=gray] (A12) -- (A22);
\draw[draw=gray] (A11) -- (A21);
\draw[draw=gray,dashed,opacity=0.5] (A3) -- (A4);
\end{tikzpicture}
\end{center}
\end{center}
For example,  $\set{x,x+y,u}$ is obtained from the construction $V= z(y(x),u)$ as prescribed by the specification of $\hyper{H}_2^v$: slowly, we have 
$$\psi(V)\setminus H_1=\set{\set{x},\set{x,y},\set{u}}\:,$$
from which we  get  $\set{x,x+y,u}$ by applying $\sigma$ elementwise (no $\mu$ to perform here).
%

\smallskip
The first round induces
$$\begin{array}{rcl}
H_2 &=& \set{x , y , z,u, x+y}\\
 \hyper{H}_2^v & = & \set{\set{y,z,u}, \set{x,z,u}, \set{y,x+y,z}, \set{x,x+y,z},\set{y,x+y,u}, \set{x,x+y,u}}\;,
\end{array}$$
and let the second round be instructed by 
$$\hyper{H}_2^t=\set{\set{u} , \set{x} , \set{y} , \set{z} , \set{x+y} ,  \set{x,x+y} , \set{u,x,y,z,x+y}}.$$
resulting in the following polytope:
\begin{center}
\begin{tikzpicture}[thick,scale=2.7]
\coordinate (A11) at (0.285,0.8);
\coordinate (A12) at (-0.285,0.8);
\coordinate (B12) at (-0.1,0.8);
\coordinate (B22) at (-0.42,0.65);
\coordinate (B32) at (-0.1,-0.8);
\coordinate (B42) at (-0.42,-0.65);
\coordinate (A21) at (0.285,-0.8); 
\coordinate (A22) at (-0.285,-0.8); 
\node (dl) at (-1.4,0) {\small $\set{y,x+y}(x,z,u)$};
\node (dr) at (1.4,0) {\small $ \set{x,x+y}(y,z,u)$};
\node (a22) at (-0.9,-0.65) {\small $\set{y,u}((x+y)(x),z)$}; 
\node (a21) at (0.76,-0.8) {\small $\set{x,u}(x+y,y,z)$}; 
\node (a11) at (-0.9,0.65) {\small $\set{y,z}((x+y)(x),u)$}; 
\node (a12) at (0.76,0.8) {\small $\set{x,z}(x+y,y,u)$}; 
\node (c12) at (-0.1,0.925) {\small $\set{y,z}(x(x+y),u)$}; 
\node (d12) at (-0.1,-0.925) {\small $\set{y,u}(x(x+y),z)$}; 
\coordinate (A3) at (-1,0);
\coordinate (A4) at (1,0);
\draw[draw=gray,fill=pur,opacity=0.15]  (A21) -- (A11) -- (A4) -- cycle;
\draw[draw=gray,fill=pur,opacity=0.15] (A3) -- (B22) -- (B42) -- cycle;
\draw[draw=gray]  (A21) -- (A11) -- (A4) -- cycle;
\draw[draw=gray] (A3) -- (B22) -- (B42) -- cycle;
 \draw[draw=gray,fill=yellow,opacity=0.1] (B12) -- (B32) -- (A21) -- (A11) -- cycle;
 \draw[draw=gray] (B12) -- (B32) -- (A21) -- (A11) -- cycle;
\draw[draw=gray,fill=yellow,opacity=0.3] (B12) -- (B22) -- (B42) -- (B32) -- cycle;
\draw[draw=gray] (B12) -- (B22) -- (B42) -- (B32) -- cycle;
\draw[draw=gray,dashed,opacity=0.5] (A3) -- (A4);
\end{tikzpicture}
\end{center}
in which the new edge between $x$ and $x+y$ (created after the first round) has been itself truncated.
This induces
 $$\begin{array}{rcl}
H_3 &=& \set{x , y ,  z, u, x+y, 2x+y}\\
\hyper{H}_3^v &=&\{ \set{x,z,u} ,  \set{y,z,u} , \set{x+y,y,z} , \set{x+y,y,u} , \set{x+y,2x+y,z} ,   \\ 
 &&\;\:\set{x,2x+y,z} , \set{x+y,2x+y,u} , \set{x,2x+y,u}\}\end{array}$$
 Here, say, $\set{x+y,2x+y,u}$ corresponds to $\set{y,z}(x(x+y),u)$  (note the use of $\mu$ on $x+(x+y)$).

\smallskip
One could go on on this example:  we could truncate the new edge between the faces $x$ and $2x+y$, 
and create the new face $3x+y$, etc.   It can be shown (see \cite{P-SISP}[Section 6]) that
the flattening from $x+(x+y)$  to $2x+y$ (or of
$x+(2x+y)$ to $3x+y$) incurs no loss of information, provided the traces of the rounds (i.e., the successive pairs of hypergraphs) are recorded. In this way, an untractable combinatorial explosion in the description of iterated truncations is avoided.

\subsection{The permutohedron-based associahedron}  \label{permutohedron-based-associahedron-section}

As a more sophisticated example of iterated   truncations, we now describe the combinatorics of the family of permutohedron-based associahedra, which are polytopes describing the coherences of symmetric monoidal categories (see Figures \ref{pba-3} and \ref{3-net}). They were introduced in \cite{P-SISP}, and further studied in \cite{BIP-SPA}. These polytopes are different from the  permutoassociahedra, which were introduced for the same purpose in
 \cite{K-Permuto}, and which are not simple polytopes. The reason for this diversity is that different choices of generating isomorphisms lead to different combinatorial / geometrical interpretations of the same coherence theorem.

\smallskip
We take ${\cal X}=\set{x_1,\ldots,x_{n+1}}$.  The first round of truncation is that producing the permutohedron, with as truncation hypergraph the complete graph over ${\cal X}$:
$$\begin{array}{l}
H_1={\cal X}\\
\hyper{H}_1^v=\setc{{\cal X}\setminus\set{x}}{x\in {\cal X}}\\
\hyper{H}_1^t =\setc{\set{x_i}}{i\in[1,n+1]}
\union\setc{\set{x_i,x_j}}{i,j\in [1,n+1],\; i\neq j}\;.
\end{array}$$
Recall from Section \ref{construct-examples-section} that all  constructs of the permuohedra are filiform.  Thus constrs are in bijection with proper subsets of ${\cal X}$  (different from ${\cal X}$), and the set of constructions is in bijection with the symmetric group $S_{n+1}$.  More precisely, a construction  
$x_{\sigma(n+1)}(\ldots (x_{\sigma(2)}(x_{\sigma(1)}))\ldots)$
is encoded as 
\begin{figure}[b]
\includegraphics[scale=0.5]{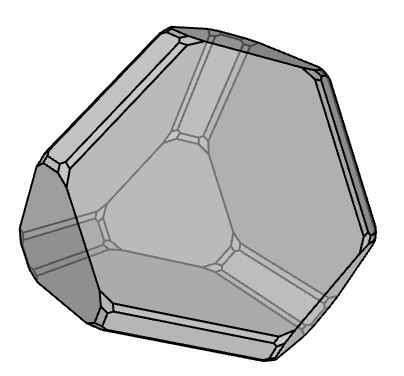}
\caption{The 3-dimensional permutohedron-based associahedron}
\label{pba-3}
\end{figure}
$$x_\sigma \eqdef \set{x_{\sigma(1)}\:,\:x_{\sigma(1)}+x_{\sigma(2)}\:,\:\ldots\:,\: x_{\sigma(1)}+\ldots+x_{\sigma(n)}}\; ,$$ where $\sigma\in S_{n+1}.$
This leads us to
 $$\begin{array}{l}
H_2=\setc{\setc{\Sigma_{i\in I} x_i}{i\in I}}{ \emptyset\neq I\incs{\cal X}}\\
\hyper{H}_2^v=\setc{x_\sigma}{\sigma\in S_{n+1}}\;.
\end{array}$$
We shall write $x_I=\Sigma_{i\in I}x_i$.
The next round of truncations is defined by
$$\hyper{H}_2^t=\setc{\set{x_I}}{\emptyset\neq I\incs{\cal X}}\union
\setc{\set{x_I,x_J}}{I\inc J\:\mbox{and}\:J\setminus I\: \mbox{is a singleton}}\;.$$
We note that, for each $V\in \hyper{H}_2^v$, $(\hyper{H}_2^t)_V$ is the hypergraph specifying the associahedron.
Now, recall  from Section \ref{construct-examples-section} that in the setting of associahedra we have a  bijective correspondence between (fully) parenthesized words and constructions. In the present case, it is guided by the following picture:
\vspace{-.1cm}
\begin{center}
{\resizebox{12cm}{!}{\begin{tikzpicture}
    \node (A) at (-0.5,0.5) {$x_{\sigma(1)}$};
    \node (B) at (0.5,0) {$x_{\sigma(1)}$};
    \node (C) at (1.7,0.5) {$x_{\sigma(2)}$};
    \node (D) at (3.4,0) {$x_{\sigma(1)}+x_{\sigma(2)}$};
\node (E) at (6.1,0) {$\ldots$};
\node (F) at (9.3,0) {$x_{\sigma(1)}+\ldots,+x_{\sigma(n)}$};
    \node (G) at (11.7,0.5) {$x_{\sigma(n+1)}$};
  \node (H) at (5.1,0.5) {$x_{\sigma(3)}$};
  \node (I) at (7.1,0.5) {$x_{\sigma(n)}$};
    \draw[-] (B)--(D) ;
  \draw[-] (D)--(E) ;
  \draw[-] (E)--(F) ;
      \end{tikzpicture}}}
\end{center}
\vspace{-.1cm}
We get that the set of constructs of $\hyper{H}_2^t$ tamed by 
$x_\sigma$ is in one-to-one correspondence with the set of  parenthesized words over ${\cal X}$, in which the order of the letters from left to right is the one that we adopted in the  definition of $x_\sigma$.  
This takes care of all 24 pentagons (corresponding to all possible permutations $\sigma$) of  Figure \ref{pba-3}..
We next show how to name the remaining edges and faces.
\begin{itemize}
\item Let us set $a=x_1, b=x_1+x_2, c= x_1+x_2+x_3, d=x_1+x_3$.
Then 
\vspace{-.05cm}
$$\left.\begin{array}{r}(H_2\setminus\set{a,b,c})(b(a,c))
\\(H_2\setminus\set{a,c,d})(d(a,c))
\end{array}\right\}
\; \mbox{correspond to}\;
\left\{\begin{array}{l}
(x_1x_2)(x_3x_4)\\
(x_1x_3)(x_2x_4)
\end{array}\right.$$
\vspace{-.2cm}

There is an edge between these two vertices, named by $(H_2\setminus\set{a,c})(a,c)$. Here is a way to name it in the style of parenthesized words:  
\vspace{-.05cm}
$$ ((x_1\bcdot_1)(\bcdot_1x_4)\:,\: (\bcdot_1\mapsto \set{x_2,x_3}))$$
\vspace{-.35cm}

The notation here is a way  to formalise the  surjection that maps  $x_1$ to $x_1$, $x_2,x_3$ to $\bcdot_1$, and $x_4$ to $x_4$.  After all,  we are seeking a mix of the notation for associahedra and permutohedra, hence a mix of parenthisations and surjections!

In this way, we  account for all single edges relating two pentagons.

\item 
We now account for parallel edges between two pentagons, and the corresponding rectangular faces:
$$\begin{array}{l}
\left.\begin{array}{r}
(\bcdot_1(\bcdot_1(x_3x_4)),(\bcdot_1\!\mapsto\! \set{x_1,\!x_2}))\\
(\bcdot_1((\bcdot_1x_3)x_4),(\bcdot_1\!\mapsto\! \set{x_1,\!x_2}))
\end{array}\right\}
\; \mbox{for}\;
\left\{\begin{array}{l}
x_1(x_2(x_3x_4)) -x_2(x_1(x_3x_4))\\
x_1((x_2x_3)x_4) - x_2((x_1x_3)x_4)
\end{array}\right.
\\\\
\quad\bcdot_1(\bcdot_1x_3x_4) \;\;\mbox{for}\;\;
(H_2\setminus\set{b,c})(\set{b,c})
\end{array}$$
and
$$\begin{array}{l}
\left.\begin{array}{r}
((x_1(x_2\bcdot_1))\bcdot_1,(\bcdot_1\!\mapsto\! \set{x_3,\!x_4}))\\
(((x_1x_2)\bcdot_1)\bcdot_1,(\bcdot_1\!\mapsto\! \set{x_3,\!x_4}))
\end{array}\right\}
\; \mbox{for}\;
\left\{\begin{array}{l}
(x_1(x_2x_3))x_4 - (x_1(x_2x_4))x_3\\
 ((x_1x_2)x_3)x_4 -  ((x_1x_2)x_4)x_3
\end{array}\right.
\\\\
\quad(x_1x_2\bcdot_1)\bcdot_1 \;\;\mbox{for}\;\;
(H_2\setminus\set{a,b})(\set{a,b})
\end{array}$$
\item 
We are left with the remaining faces.
The eight dodecagons are named by
$$((x_i\bcdot_1)\bcdot_1\bcdot_1,(\bcdot_1\mapsto (H_2\setminus\set{x_i}))\enspace\mbox{and} \enspace
(\bcdot_1\bcdot_1(\bcdot_1 x_i) , (\bcdot_1\mapsto (H_2\setminus\set{x_i}))$$
standing for  $(H_2\setminus\set{x_i})(x_i)$ and $(H_2\setminus\set{\sum_{j\neq i}x_j})(\sum_{j\neq i}x_j)$, respectively,
and the 6 octagons by, say:
$$(\bcdot_1(\bcdot_1\bcdot_2)\bcdot_2,(\bcdot_1\mapsto\set{x_1,x_2}\,,\, \bcdot_2\mapsto\set{x_3,x_4})).$$
Indeed, this octagon should contain the following  four edges (which are sides of four pentagons), for each of which we give the corresponding construct:
$$\begin{array}{lll}
x_1(x_2x_3)x_4 &&(H_2\setminus\set{a,b,c})(\set{a,c})(b))\\
x_1(x_2x_4)x_3 &&(H_2\setminus\set{a,b,f})(\set{a,f})(b))\\
x_2(x_1x_3)x_4 &&(H_2\setminus\set{e,b,c})(\set{e,c})(b))\\
x_2(x_1x_4)x_3 &&(H_2\setminus\set{e,b,f})(\set{e,f})(b))
\end{array}$$
where  $a= x_1$,  $b=x_1+x_2$, $c=x_1+x_2+x_3$, $e=x_2$,  $f=x_1+x_2+x_4$.  
The least upper bound of these constructs is $(H_2\setminus\set{b})(b)$, and all what this construct specifies is that  we should do the operation $b$ as innermost operation.  It is a ``Mastermind" kind of 
partial information:
\begin{center}
\resizebox{8cm}{!}{\begin{tikzpicture}
    \node (A) at (-0.5,0.5) {$\bcdot_1$};
    \node (B) at (0.5,0) {$?_1$};
    \node (C) at (1.7,0.5) {$\bcdot_1$};
    \node (D) at (3.4,0) {$x_1+x_2$};
\node (F) at (6.3,0) {$?_2$};
  \node (H) at (5.1,0.5) {$\bcdot_2$};
  \node (I) at (7.4,0.5) {$\bcdot_2$};
    \draw[-] (B)--(D) ;
  \draw[-] (D)--(F) ;
      \end{tikzpicture}}
\end{center}
Note that $b$, being the sum of two letters, has to be the central node, and that it being the sum of $x_1$ and $x_2$  entails that  $?_1$ is  $x_1$ or $x_2$, and $?_2$ is $x_1+x_2+x_3$ or $x_1+x_2+x_4$.  The same information is carried out by our encoding.
\end{itemize}

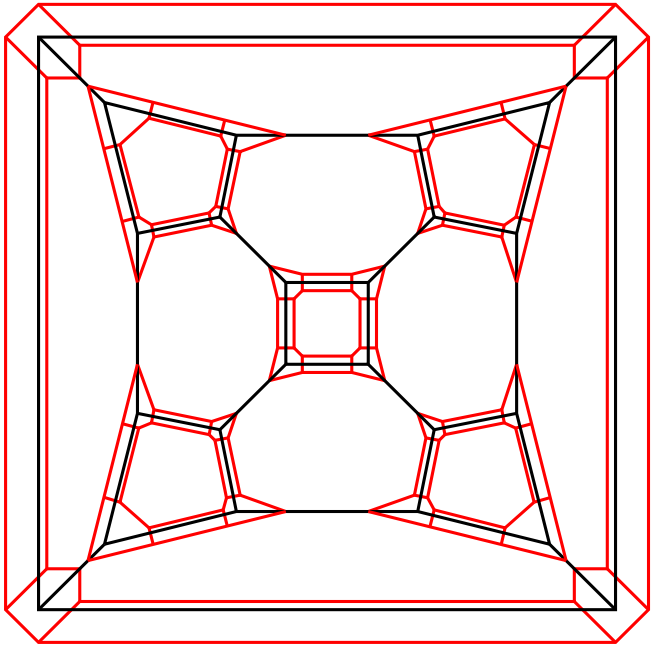
\begin{figure}[b] 
\resizebox{8.5cm}{!}{\begin{tikzpicture}[thick,scale=0.3]
\coordinate (A1) at (0,0);
\coordinate (A2) at (17.5,0);
\coordinate (A3) at (18.5,1);
\coordinate (A4) at (18.5,18.5);
 \coordinate (A5) at (17.5,19.5);
 \coordinate (A6) at (0,19.5);
 \coordinate (A7) at (-1,18.5);
 \coordinate (A8) at (-1,1);
 \coordinate (B1) at (0,1);
 \coordinate (B2) at (17.5,1);
 \coordinate (B3) at (17.5,18.5);
 \coordinate (B4) at (0,18.5);
 \coordinate (P1) at (1.25,1.25);
 \coordinate (P2) at (1.25,2.25);
 \coordinate (P3) at (0.25,2.25);
  \coordinate (Q1) at (16.25,1.25);
  \coordinate (Q2) at (16.25,2.25);
  \coordinate (Q3) at (17.25,2.25);
  \coordinate (R1) at (17.25,17.25);
  \coordinate (R2) at (16.25,17.25);
  \coordinate (R3) at (16.25,18.25);
 \coordinate (S1) at (1.25,18.25);
 \coordinate (S2) at (1.25,17.25); 
\coordinate (S3) at (0.25,17.25);
 \coordinate (K1) at (2,3);
 \coordinate (K2) at (15.5,3);
  \coordinate (K3) at (15.5,16.5);
  \coordinate (K4) at (2,16.5);
  \coordinate (F1) at (6,4);
 \coordinate (F2) at (11.5,4);
 \coordinate (F3) at (14.5,7);
 \coordinate (F4) at (14.5,12.5);
 \coordinate (F5) at (11.5,15.5);
 \coordinate (F6) at (6,15.5);
 \coordinate (F7) at (3,12.5);
 \coordinate (F8) at (3,7);
  \coordinate (G1) at (5.5,6.5);
 \coordinate (G2) at (12,6.5);
  \coordinate (G3) at (12,13);
 \coordinate (G4) at (5.5,13);
  \coordinate (J1) at (7.5,8.5);
 \coordinate (J2) at (10,8.5);
  \coordinate (J3) at (10,11);
 \coordinate (J4) at (7.5,11);
  \coordinate (J1') at (8,8.75);
 \coordinate (J2') at (9.5,8.75);
 \coordinate (J2'') at (9.75,9);
  \coordinate (J3') at (9.75,10.5);
  \coordinate (J3'') at (9.5,10.75);
 \coordinate (J4') at (8,10.75);
 \coordinate (J4'') at (7.75,10.5);

  \coordinate (J1'') at (7.75,9);
  \coordinate (J'1) at (7,8);
 \coordinate (J'2) at (10.5,8);
  \coordinate (J'3) at (10.5,11.5);
 \coordinate (J'4) at (7,11.5);
 \coordinate (M1) at (1.5,2.5);
 \coordinate (M1') at (1.5,17);
 \coordinate (M2) at (7.5,4);
 \coordinate (M3) at (10,4);
 \coordinate (M4) at (16,2.5);
 \coordinate (M5) at (3,8.5);
 \coordinate (M6) at (3,11);
 \coordinate (M7) at (14.5,8.5);
 \coordinate (M8) at (14.5,11);
 \coordinate (M9) at (16,17);
 \coordinate (M10) at (10,15.5);
 \coordinate (M11) at (7.5,15.5);
  \coordinate (J1'-) at (8,8.25);
 \coordinate (J2'-) at (9.5,8.25);
 \coordinate (J2''-) at (10.25,9);
  \coordinate (J3'-) at (10.25,10.5);
 \coordinate (J3a'') at (9.5,11.25);
 \coordinate (J4a') at (8,11.25);
 \coordinate (J4''-) at (7.25,10.5);
  \coordinate (J1''-) at (7.25,9);

\draw[thick,red] (J'1) -- (J1'-) -- (J2'-) -- (J'2) -- (J2''-) -- (J3'-) -- (J'3) -- (J3a'') -- (J4a') -- (J'4) -- (J4''-) -- (J1''-) -- (J'1);
\draw[thick,red] (J1') -- (J1'-);
\draw[thick,red] (J2') -- (J2'-);
\draw[thick,red] (J2'') -- (J2''-);
\draw[thick,red] (J3') -- (J3'-);
\draw[thick,red] (J3'') -- (J3a'');
\draw[thick,red] (J4') -- (J4a');
\draw[thick,red] (J4'') -- (J4''-); 
\draw[thick,red] (J1'') -- (J1''-); 

  \coordinate (G1') at (6,7);
 \coordinate (G2') at (11.5,7);
  \coordinate (G3') at (11.5,12.5);
 \coordinate (G4') at (6,12.5);

 \coordinate (G1'-) at (5.25,6.75);
 \coordinate (G1'+) at (5.75,6.25);
 \coordinate (G1'--) at (5.175,6.35);
 \coordinate (G1'++) at (5.35,6.175);

 \coordinate (G1'---) at (3.425,6.71);
 \coordinate (W1) at (2.47,4.285);
 \coordinate (W2) at (3.35,3.5);

 \coordinate (G2'-) at (11.75,6.25);
 \coordinate (G2'+) at (12.25,6.75);

  \coordinate (G3'-) at (11.75,13.25);
  \coordinate (G3'+) at (12.25,12.75);

 \coordinate (G4'-) at (5.25,12.75);
 \coordinate (G4'+) at (5.75,13.25);

\draw[thick,red] (M5) -- (3.5,7.11) -- (G1'-) -- (G1') -- (G1'+) -- (6.115,4.5) -- (M2);
\draw[thick,red] (G1'--) -- (G1'-);
\draw[thick,red] (3.5,7.11) -- (G1'---) -- (3.04,6.55) -- (W1) -- (W2) -- (5.6,4.05) -- (5.7,4.42) --(6.115,4.5);
\draw[thick,red] (3.04,6.55) -- (2.53,6.687);
\draw[thick,red]   (W1) -- (2,4.43);
\draw[thick,red] (G1'--) -- (G1'---);
\draw[thick,red] (G1'++) -- (G1'+);
\draw[thick,red] (G1'--) -- (G1'++);
\draw[thick,red] (5.7,4.42) -- (G1'++);
\draw[thick,red] (5.6,4.05) -- (5.726,3.55);
 \draw[thick,red] (W2) -- (3.476,3);

\draw[thick,red] (M7) -- (14.05,7.11) -- (G2'+) -- (G2') -- (G2'-) -- (11.4 ,4.5) -- (M3);
 \coordinate (W1') at (15.03,4.285);
 \coordinate (W2') at (14.15,3.5);
\coordinate (G2'++) at (12.325,6.35);
 \coordinate (G2'+++) at (14.125,6.71);
 \coordinate (G2'--) at (12.15,6.175);
\draw[thick,red]  (G2'++) -- (G2'--);
\draw[thick,red]  (G2'+++) -- (G2'++) -- (G2'+);
 \draw[thick,red] (W1')--(W2')-- (12,4.05)-- (11.8,4.42);
 \draw[thick,red] (12,4.05)-- (11.875,3.55);
\draw[thick,red] (G2'--) -- (G2'-);
\draw[thick,red] (11.4 ,4.5) -- (11.8,4.42) -- (G2'--);
\draw[thick,red] (14.05,7.11)-- (G2'+++)--(14.47,6.55) -- (14.97,6.687);
 \draw[thick,red] (14.47,6.55) -- (W1');
\draw[thick,red] (W2') -- (14.035,3);
 \draw[thick,red]   (W1') -- (15.5,4.42);

\draw[thick,red] (M8) -- (14.05,12.39) -- (G3'+) -- (G3') -- (G3'-) -- (11.4 ,15) -- (M10);
 \coordinate (W1") at (15.04,15.215);
 \coordinate (W2") at (14.15,16);

\draw[thick,red] (W1")--(W2");
\coordinate (G3'++) at (12.325,13.125);
 \coordinate (G3'+++) at (14.125,12.75);
 \coordinate (G3'--) at (12.15,13.325);
\draw[thick,red]  (G3'++) -- (G3'--);
\draw[thick,red]  (G3'+++) -- (G3'++) -- (G3'+);
 \draw[thick,red] (11.4 ,15.00) -- (11.8,15.075) -- (G3'--);
  \draw[thick,red] (11.8,15.075) -- (12,15.473) -- (W2");
\draw[thick,red] (12,15.473)-- (11.875,15.973);
\draw[thick,red] (W2")--(14.03,16.5);
\draw[thick,red] (W1")--(15.525,15.093);
\draw[thick,red] (W1")--(14.48,13) -- (G3'+++)-- (14.05,12.39);
\draw[thick,red] (14.48,13) -- (14.975,12.87);
\draw[thick,red] (G3'--) -- (G3'-);

\draw[thick,red] (M11) -- (6.115,15) -- (G4'+) -- (G4') -- (G4'-) -- (3.5,12.39) -- (M6);
 \coordinate (W1") at (2.47,15.215);
 \coordinate (W2") at (3.35,16.02);
\draw[thick,red] (W1")--(W2");
\coordinate (G4'--) at (5.17,13.125);
 \coordinate (G4'---) at (3.425,12.765);
 \coordinate (G4'++) at (5.375,13.325);
\draw[thick,red] (G4'--) -- (G4'-);
\draw[thick,red] (G4'++) -- (G4'+);
\draw[thick,red]  (G4'++) -- (G4'--);
 \draw[thick,red]  (6.115,15) -- (5.725,15.075) -- (G4'++);
  \draw[thick,red] (5.725,15.08) -- (5.525,15.475) -- (W2");
\draw[thick,red] (5.525,15.473)-- (5.65,15.973);
\draw[thick,red] (W2")--(3.47,16.5);
\draw[thick,red] (W1")--(2,15.093);
\draw[thick,red] (W1")--(3.045,13)-- (G4'---) -- (G4'--);
\draw[thick,red] (G4'---) -- (3.5,12.39);
\draw[thick,red] (2.55,12.868) -- (3.045,13);

\draw[thick,red] (A1) -- (A2) -- (A3) -- (A4) -- (A5) -- (A6) -- (A7) -- (A8) -- cycle;
 \draw[thick,red]  (A1) -- (P1) -- (P2) -- (P3) -- (A8);
  \draw[thick,red]  (A2) -- (Q1) -- (Q2) -- (Q3) -- (A3);
  \draw[thick,red]  (A4) -- (R1) -- (R2) -- (R3) -- (A5);
  \draw[thick,red]  (A6) -- (S1) -- (S2) -- (S3) -- (A7);
 \draw[thick,red]   (S3) -- (P3);
 \draw[thick,red]   (P1) -- (Q1);
 \draw[thick,red]   (Q3) -- (R1);
 \draw[thick,red]   (S1) -- (R3);
 \draw[thick] (B1) -- (B2) -- (B3) -- (B4) -- cycle;
  \draw[thick] (K1) -- (F1) -- (F2) -- (K2) -- (F3) -- (F4) -- (K3) -- (F5) -- (F6) -- (K4) -- (F7) -- (F8) -- (K1);
  \draw[thick] (F1) -- (G1) -- (F8);
  \draw[thick] (F2) -- (G2) -- (F3);
 \draw[thick] (F4) -- (G3) -- (F5);
 \draw[thick] (F6) -- (G4) -- (F7);
 \draw[thick] (B1) -- (K1);
 \draw[thick] (B2) -- (K2);
 \draw[thick] (B3) -- (K3);
 \draw[thick] (B4) -- (K4);
 \draw[thick] (J1) -- (J2) -- (J3) -- (J4) -- cycle;
 \draw[thick,red] (J1') -- (J2') -- (J2'') -- (J3') -- (J3'') -- (J4') -- (J4'') -- (J1'') -- cycle;
\draw[thick] (G1) -- (J1);
 \draw[thick] (G2) -- (J2);
 \draw[thick] (G3) -- (J3);
 \draw[thick] (G4) -- (J4);
 \draw[thick,red] (M1) -- (M2);
 \draw[thick,red] (M3) -- (M4);
 \draw[thick,red] (M7) -- (M4);
 \draw[thick,red] (M1) -- (M5);
 \draw[thick,red] (M1') -- (M6);
 \draw[thick,red] (M8) -- (M9);
 \draw[thick,red] (M10) -- (M9);
 \draw[thick,red] (M1') -- (M11);
\end{tikzpicture}}
\caption{Plane projection of the  permutohedron-based associahedron of dimension 3}
\label{3-net}
\end{figure}
The   notation can be systematised in any finite dimension. We shall describe an algorithm transforming any 
 construct $T=(H_2\setminus Y)(T_1,\ldots, T_p)$ of $\hyper{H}_2^t$ rel to
$\hyper{H}_2^v$ into a pair $(W,\pi)$ of the kind discovered above. 
We can write $$Y = \set{\sum_{i\in I_1}x_i,\ldots, \sum_{i\in I_k}x_i}\quad\quad\mbox{with}\; I_1\incs I_2\incs\ldots\incs I_k\;.$$
We can encode the information provided by $Y$ through the  following map  $\tilde\pi_Y$:
$$\begin{array}{llllllllll}
\tilde\pi_Y(i1)=\setc{x_i}{i\in I_1}, &&&& \tilde\pi_Y(2)=\setc{x_i}{i\in I_2\setminus I_1},\ldots ,
\\
\tilde\pi_Y(k)=\setc{x_i}{i\in I_k\setminus I_{k-1}} &&&& \tilde\pi_Y(k+1)=\setc{x_i}{i\in [1,n+1]\setminus I_k}
\;.
\end{array}$$
We associate with $\tilde\pi_Y$ the word $\tilde w_Y$ starting with $|I_1|$ occurrences of the letter  $\bcdot_1$, followed by 
$|I_2\setminus I_1|$ occurrences of the letter $\bcdot_2$,\ldots, ending with $|[1,n+1]\setminus I_k|$ occurrences of $\bcdot_{k+1}$.  We then do a bit of ``making up'':
we replace in $\tilde w_Y$ all letters $\bcdot_i$  occurring only once by the unique element of $\tilde\pi_Y(i)$; we also renumber the remaining letters $\bcdot_j$, and we reindex $\tilde\pi_Y$ accordingly. We denote the new word by $w_Y$ and the new map by $\pi_Y$
  We call the $\cdot_j$'s and the $x_i$'s   holes and determined letters, respectively.
For example, if  
$$\begin{array}{llllll}
\tilde\pi_Y(1)= \set{x_9}\ &&\tilde\pi_Y(3)=\set{x_3}\ && \tilde\pi_Y(5)=\set{x_6}\\
\tilde\pi_Y(2)=\set{x_2,x_4,x_8} && \tilde\pi_Y(4)=\set{x_1,x_7} && \tilde\pi_Y(6)=\set{x_5,x_7}\;,
\end{array}$$	
then $\bcdot_1\bcdot_2\bcdot_2\bcdot_2\bcdot_3\bcdot_4\bcdot_4\bcdot_5\bcdot_6\bcdot_6$ becomes 
$$x_9 \bcdot_1\bcdot_1\bcdot_1x_3\bcdot_2\bcdot_2 x_6 \bcdot_3\bcdot_3\;,$$ and we have  
$$\pi_Y(1) =\set{x_2,x_4,x_8}\quad \pi_Y(2) \set{x_1,x_7}\quad \pi_Y(3)=\set{x_5,x_7}\;.$$ 
We complete the making up by placing parentheses  in $w_Y$: \begin{itemize}
\item  around every  subword 
$\bcdot_j  x_{i_1}\ldots x_{i_l} \bcdot_{j+1}$ ($l$ may be $0$), and, if this applies,
\item   around  the prefix of $w_Y$ of the form   $x_{i_1}\ldots x_{i_l} \bcdot_1$ ($l>0$),
\item
   and around the suffix of $w_Y$ of the form   $\bcdot_k  x_{i_1}\ldots x_{i_l} $ ($l>0$).
\end{itemize}
(with all letters in the $\ldots$'s determined). We denote the obtained parenthesised word  by $(w_Y\!)^{\it st}$ (for  {\em st}andardisartion). We use square brackets for writing the parentheses  in $(w_Y\!)^{\it st}$, to distinguish them (visually only)  from further parentheses that will be induced by $T_1,\ldots,T_p$.
For our example, we get:
$$[x_9 \bcdot_1]\bcdot_1[\bcdot_1x_3\bcdot_2][\bcdot_2 x_6 \bcdot_3]\bcdot_3$$
We now examine how to encode the information provided by $T_1,\ldots,T_p$.
The square brackets in $(w_Y\!)^{\it st}$ delimit the zones of $w_Y$ that correspond to the connected components of $(\hyper{H}_2^t)_Y$, which have the form 
$$\set{\sum_{i\in I_m}x_i,\sum_{i\in I_{m+1}}\!\! x_i,\ldots, \sum_{i\in I_{m+q}}\!\! x_i}\; \mbox{with}\left\{\begin{array}{l}
I_m\inc I_{m+1}\inc\ldots\inc I_{m+q}\\
|I_{m+1}\setminus I_m|=\ldots=|I_{m+q}\setminus I_{m+q-1}|=1
.
\end{array}\right.$$
 We  
examine first the two  degenerate cases:
\begin{itemize}
\item $Y=\emptyset$. Then $k=0$, and we set by convention $I_0=\emptyset$, so that $\pi(1)=\set{1,\ldots,n+1}\setminus I_0=\set{1,\ldots,n+1}$. Then $(w_Y\!)^{\it st}= w_Y=\bcdot_1\ldots\bcdot_1$ (with
  length $n+1$), which encodes the maximum face, i.e., the entire polytope.
\item $Y=x_\sigma$ for some $\sigma$. Then all sets $I_1$, $I_2\setminus I_1$,\ldots, $\set{1,\ldots,n+1}\setminus I_k$ are singletons, and the construct $T$  is of the form $(H\setminus x_\sigma)(S)$, where $S$ is a construct of the associahedron generated by the hypergraph 
$$ \hyper{H}_\sigma \eqdef  \{\set{x_{\sigma(1)}},\ldots,\set{x_{\sigma(n+1)}}, \set{x_{\sigma(1)},x_{\sigma(2)}},\ldots ,\set{x_{\sigma(n)},x_{\sigma(n+1)}}\}.$$
It follows that $(w_{x_\sigma}\!)^{\it st}= w_{x_\sigma}=x_{\sigma(1)}\:x_{\sigma(2)}\:\ldots\: x_{\sigma(n+1)}$. 
Then  $S$ determines a parenthesisation of this word.
\end{itemize}
In the non-degenerate cases, if we fix a permutation $\sigma$ such that $\emptyset\incs Y\incs x_\sigma$,  we can show that $T_\sigma=(x_\sigma\setminus Y)(T_1,\ldots,T_p)$ is a construct of $\hyper{H}_{x_\sigma}$, hence the data
$T_1$,\ldots, $T_p$ amount to giving parenthesisations in the $p$ zones delimited by the square parentheses, resulting in a parenthesised word which we denote by $W_T$. 
It can be shown easily that the synthesis of $W_T$ does not depend on the choice of $\sigma$ such that $Y\incs x_\sigma$.

Our analysis also identifies the target of the translation associating $W_T$ and $\pi_Y$ to $T=(H_2\setminus Y)(T_1,\ldots, T_p)$.  It consists of all pairs 
$(W,\pi)$, where 
\begin{itemize}
\item $\pi$ is a map from $\set{1,\ldots,q}$ to the set of subsets  of ${\cal X}$ of cardinality at least $2$, for some $q$,
\item  $W$ is a parenthesised word  over
${\cal X}\union\set{\bcdot_1,\ldots,\bcdot_{q}}$
such that, writing $\overline{W}$ for the word obtained by removing all parentheses from $W$:
\begin{itemize}
\item each letter $x_i$ appears at most once in $\overline{W}$;
\item for all $i\in\set{1,\ldots,q}$, all the occurrences of  $\bcdot_i$ appear as a block of length $|\pi(i)|$ in $\overline{W}$, and before any occurrence of $\bcdot_{i+1}$ (if $i<q$).
\item
the sets $\pi(r)$  (for $r$ ranging over $\set{1,\ldots,q}$) and the singletons $\set{x_i}$ such that $x_i$ appears in $\overline{W}$  form   a partition of ${\cal X}$;
\item  all the parentheses of $W$ are within the scope of some parentheses of $(\overline{W})^{{\it st}}$ (as defined above), and $W$ carries all the parentheses of  $(\overline{W})^{{\it st}}$.
\end{itemize}
As an example, in reference to the avove example of standardisation, 
$$[x_9 \bcdot_1]\bcdot_1[\bcdot_1(x_3\bcdot_2)][\bcdot_2 x_6 \bcdot_3]\bcdot_3$$
respects the scoping condition. As a prototypical counter-example, the word
$(x_1x_2)(\bcdot_1 \bcdot_1)$ is not accepted, since $(x_1x_2\bcdot_1 \bcdot_1)^{\it st}=  [x_1x_2\bcdot_1]\bcdot_1$.

\end{itemize}

We leave to the reader the proof that the translation is well defined and bijective, and that the following description of a partial order 
makes it actually an isomorphism:
\begin{itemize}
\item
 set $(W,\pi)\leq (W',\pi)$ if  $W'$ has one pair of parentheses (other than the standard ones) less than $W$;
\item
set $(W,\pi)\leq (W',\pi')$, if  $W$ inherits the parentheses of  $W'$, and if $\pi$ is an elementary refinement of $\pi'$, i.e.,  
$\pi(1)=\pi'(1),\ldots,\pi(i-1)=\pi'(i-1)$, $\pi(i)\union \pi(i+1)=\pi'(i)$,  
$\pi(i+2)=\pi'(i+1),\ldots$, up to ``making up'';
\item close by reflexivity and transitivity.
\end{itemize}
For example, we have:
$$\begin{array}{lllllll}
\bcdot_1((\bcdot_1 x_3)x_4) <  \bcdot_1 (\bcdot_1 x_3x_4) &&  (\pi(1)=\set{x_1,x_2}) &&\mbox{by the first rule}\\
\bcdot_1((\bcdot_1 x_3)x_4) <  \bcdot_1(\bcdot_1\bcdot_2)\bcdot_2  && (\pi(1)=\set{x_1,x_2}) && \mbox{by the second rule}\\
(x_1\bcdot_1)\bcdot_1\bcdot_1 < \bcdot_1 \bcdot_1\bcdot_1\bcdot_1 && (\pi(1)=\set{x_2,x_3,x_4}) && \mbox{by the second rule} \;.
\end{array}$$
We detail the derivation of $\bcdot_1((\bcdot_1 x_3)x_4) <  \bcdot_1(\bcdot_1\bcdot_2)\bcdot_2$: the refinement splits $\pi'(2)=\set{x_3,x_4}$, yielding the standardised word 
$\bcdot_1(\bcdot_1x_3 x_4)$, and because parentheses are inherited, we indeed get $\bcdot_1((\bcdot_1 x_3)x_4) $ as a predecessor.

\vspace{-.4cm}
 
\section{Directions for future work}

\vspace{-.2cm}
We plan  to apply hypergraph polytopes to  study  other coherence problems. In recent work, the first two authors have identified the coherence conditions for categorified cyclic operads, but it is not yet clear what the relevant polytopes 
are  in this setting.
 The third author is working on  giving precise geometric realisations of the  polytopes obtained by iterated truncations. 
 The case of the permutohedron-based associahedron has already  been settled in  \cite{BIP-SPA}. 



\vspace{-.3cm}
\begin{acknowledgements}
The authors wish to thank Kosta Do\v sen and Zoran Petri\'c for  enlightening discussions.
\end{acknowledgements}


\vspace{-.55cm}
\begin{thebibliography}{99}
\bibitem{BIP-SPA} D. Barali\'c, J. Ivanovi\'c, Z. Petri\'c, A simple permutoassociahedron, arXiv:1708.02482.
\bibitem{CD-CCGA} M. Carr and S. Devadoss, Coxeter complexes and graph-associahedra, Topology and its Applications 153, (12),  2155--2168 (2006).
\bibitem{D-RGA} S. Devadoss, A realization of graph-associahedra, Discrete Mathematics 309, 2009, pp. 271--276.
\bibitem{DP-HP}  K. Do\v sen and Z. Petri\'c, Hypergraph polytopes, Topology and its Applications 158(2011), pp. 1405--1444 (arXiv:1010.5477). 
\bibitem{DP-WCO}  K. Do\v sen and Z. Petri\'c, Weak Cat-operads, Logical Methods in Computer Science 11(2015), issue 1, paper 10, pp. 1--23 (arXiv:1005.4633v8).
\bibitem{FK04} E.M. Feichtner and D.N. Kozlov, Incidence combinatorics of resolutions, Selecta Math. (N.S.) 10, 37--60  (2004). 
\bibitem{FS05} E.M. Feichtner, B. Sturmfels, Matroid polytopes, nested sets and Bergman fans, Port. Math. (N.S.) 62, 437--468 (2005). 
\bibitem{K-Permuto} M. Kapranov, The permutoassociahedron, Mac Lane's coherence theorem and asymptotic zones for the KZ equation, Journal of Pure and Applied Algebra  85 (2), 119--142.
\bibitem{LR-P} J.-L. Loday and M. Ronco, Permutads, Journal of Combinatorial Theory Series A 05/2011; 120(2). 
\bibitem{LV-AO}  J.-L. Loday,  B. Vallette,  Algebraic operads, Springer (2012).
\bibitem{ML-categories} S. Mac Lane, Categories for the working mathematician, second edition, Springer (1978).
\bibitem{P-SISP}  Z. Petri\'c, On Stretching the Interval Simplex-Permutohedron, Journal of Algebraic Combinatorics,
39  (2014), pp. 99--125.
\bibitem{PRW08} A. Postnikov, V. Reiner,  and L. Williams, Faces of generalized permutohedra, Doc. Math. 13, 207--273 (2008). 
\bibitem {P09} A. Postnikov, Permutohedra, associahedra, and beyond, Int. Math. Res. Not. IMRN 2009, 1026--1106  (2009).
\bibitem{SS94} S. Shnider and S. Sternberg, Quantum groups: from coalgebras to Drinfeld algebras, Graduate texts in mathematical physiscs, International Press (1994).
\bibitem{S-FOPIT} J.D. Stasheff, From operads to ``physically'' inspired theories, Operads:
Proceedings of Renaissance Conferences (J.-L. Loday, J.D. Stasheff, and A.A. Voronov, eds.),
Contemporary Math., vol. 202, 1997, pp. 53--81. 
\bibitem{T-RAP}  A. Tonks, Relating the associahedron and the permutohedron, same volume as \cite{S-FOPIT}.
\bibitem{Z06} A. Zelevinsky, Nested complexes and their polyhedral realizations, Pure Appl. Math. Q. 2,  655--671 (2006).
\bibitem{Ziegler} G. Ziegler, Lectures on polytopes, second edition, Springer (1998).
\end{thebibliography}


%
%

\end{document}